\documentclass{amsart}

\usepackage{mathbbol}

\usepackage{amsthm}
\usepackage[mathscr]{eucal}
\usepackage{mathrsfs}
\usepackage{mathbbol}
\usepackage{oldgerm,units}
\usepackage{wrapfig}
\usepackage{epsfig}
  \input amssymb.sty

\usepackage{epsfig,latexsym,amsfonts,amssymb,amsmath,amscd,graphics,epic}
\usepackage{amsfonts,amssymb,amsmath,amscd,amsthm}
\usepackage[mathscr]{eucal}
\usepackage{mathrsfs}
\usepackage{color}

\newtheorem{theorem}{Theorem}[section]
\newtheorem{proposition}[theorem]{Proposition}
\newtheorem{definition}[theorem]{Definition}

\newtheorem{lemma}[theorem]{Lemma}

\newtheorem{corollary}[theorem]{Corollary}

\newtheorem{example}[theorem]{Example}
\newtheorem{remark}[theorem]{Remark}

\newtheorem{axiom}[theorem]{Axiom}

\newtheorem{properties}[theorem]{Properties}

\newcommand{\Comp}{\mathbb C}
\newcommand{\Real}{\mathbb R}

\newcommand{\Int}{\mathbb Z}
\newcommand{\Net}{\mathbb N}

\newcommand{\Fld}{\mathbb K}

\newcommand{\comp}{\mathbb C}
\newcommand{\real}{\mathbb R}

\newcommand{\one}{\mathbb{1}}
\newcommand{\zero}{\mathbb{0}}

\newcommand{\chos}[2]{{#1 \choose #2}}

\newcommand{\re}[1]{\operatorname{re} \( #1 \) }

\def \hpf{\mathcal P}
\def\hf{\mathbb H}
\def\hMult{\otimes}
\def\hFrc{\oslash}
\def\hSum{\oplus}
\def\hSub{\ominus}

\newcommand{\got}[1]{\frak{#1}}

\newcommand{\To}{\longrightarrow }

\newcommand{\om}{\omega}
\newcommand{\Om}{\Omega}

\newcommand{\ep}{\epsilon}
\newcommand{\al}{\alpha}
\newcommand{\bt}{\beta}
\newcommand{\gm}{\gamma}
\newcommand{\sig}{\sigma}

\newcommand{\dl}{\delta}

\newcommand{\lm}{\lambda}
\newcommand{\Lm}{\Lambda}

\def\({\left(}
\def\){\right)}

\newcommand{\tLim}[2]{\mathop {\lim }\limits_{#1  \to #2}}

\newcommand{\bfem}[1]{\textbf{\emph{#1}}}

\def\pSkip{\vskip 1.1mm \noindent}
\newcommand{\FigWidth}{1.7}

\def\term{term}
\def\prt{part}
\def\comp{component}

\def\vep{\varepsilon}
\newcommand{\etype}[1]{\renewcommand{\labelenumi}{(#1{enumi})}}

\def\eroman{\etype{\roman}}
\def\ealph{\etype{\alph}}
\def\earabic{\etype{\arabic}}

\def\topd{reduced}
\def\top{reduction}

\def\one{\mathbb{1}}
\def\zero{\mathbb{0}}

\def\Eq{{\operatorname{Q}}}

\def\pv{{\operatorname{re}}}
\def\ph{{\operatorname{ph}}}

\def\tail{{\operatorname{T}}}
\def\head{{\operatorname{H}}}
\def\win{{\operatorname{W}}}
\def\loss{{\operatorname{L}}}

\def\minf{-\infty}
\newcommand{\supWn}{\overline{\xi}_{W_n}}
\newcommand{\supX}{\overline{\xi}_X}
\newcommand{\infX}{\underline{\xi}_X}

\newcommand{\supWnv}[1]{\supWn(#1)}
\newcommand{\supXv}[1]{\supX(#1)}
\newcommand{\infXv}[1]{\infX(#1)}

\newcommand{\itt}[1]{\tau{(#1)}}

\newcommand{\tem}[1]{\textsl{{#1}}}
\def\prv{\tem{p.\,r.\,v.}}
\def\trv{\tem{r.\,p.\,r.\,v.}}
\def\bprv{\tem{b.\,p.\,r.\,v.}}
\def\mprv{\tem{m.\,p.\,r.\,v.}}
\def\pmf{\tem{p.\,m.\,f.}}
\def\tpmf{\tem{r.\,p.\,m.\,f.}}
\def\pdf{\tem{p.\,d.\,f.}}
\def\tpdf{\tem{r.\,p.\,d.\,f.}}
\def\cpdf{\tem{c.\,p.\,d.\,f.}}
\def\mgf{\tem{m.\,g.\,f.}}

\newcommand{\range}[1]{\operatorname{R}_{_{#1}}}
\newcommand{\domain}[1]{\operatorname{D}_{_{#1}}}

\newcommand{\cov}[1]{\operatorname{Cov}[{#1}]}
\newcommand{\covr}[1]{\operatorname{Cov}_\pv[{#1}]}
\newcommand{\covp}[1]{\operatorname{Cov}_\ph[{#1}]}

\newcommand{\var}[1]{\operatorname{Var}[{#1}]}
\newcommand{\varr}[1]{\operatorname{Var}_\pv[{#1}]}
\newcommand{\varp}[1]{\operatorname{Var}_\ph[{#1}]}
\newcommand{\mean}[1]{\operatorname{E}[{#1}]}
\newcommand{\meanr}[1]{\operatorname{E}_\pv[{#1}]}
\newcommand{\meanp}[1]{\operatorname{E}_\ph[{#1}]}

\def\div{{\operatorname{div}}}
\def\zd{Z_{\div}}
\def\zdz{Z^0_{\div}}

\def\z1{{\one}}

\def\bfz{{\bf z}}
\def\bfx{{\bf x}}
\def\bfy{{\bf y}}

\def\oset{S}

\def\wk{{\operatorname{wk}}}
\def\oleq{ \precsim_\wk}
\def\ogeq{\succsim_\wk}
\def\ole{\prec_\wk}
\def\oge{\succ_\wk}
\def\oeq{\thicksim_\wk}

\def\ptr{\Lm}
\def\cptr{\bar{\ptr}}
\def\xptr{\widetilde{\ptr}}
\def\tJ{\mathcal J}

\def\alabel{\al}
\newcommand{\af}[1]{[ {#1} ]_\alabel}
\newcommand{\aqut}[1]{ {#1}/_\alabel}
\newcommand{\asim}{\simeq_\alabel}
\newcommand{\aaa}[1]{#1_{\operatorname{\alabel}}}
\def\aleq{\aaa{\lesssim}}
\def\ageq{\aaa{\gtrsim}}
\def\ale{\aaa{<}}
\def\age{\aaa{>}}
\def\aeq{\aaa{=}}

\newcommand{\lex}[1]{#1_{\operatorname{lex}}}
\def\lxleq{\lex{\leq}}

\newcommand{\tord}[1]{#1_{\operatorname{t}}}
\def\tleq{\tord{\leq}}

\def\tle{\tord{<}}

\newcommand{\sord}[1]{#1_{\operatorname{s}}}
\def\sleq{\sord{\leq}}

\def\ablabel{| \, |}
\newcommand{\ab}[1]{#1_{\operatorname{\ablabel}}}
\def\ableq{\ab{\lesssim}}

\def\rlabel{\operatorname{re}}

\newcommand{\rrr}[1]{#1_{\operatorname{\rlabel}}}
\def\rleq{\rrr{\leq}}
\def\rgeq{\rrr{\geq}}
\def\rle{\rrr{<}}
\def\rge{\rrr{>}}
\def\req{\rrr{=}}

\def\hf{ \mathbb{PH}}
\def\hfto{\hf_{(1)}}

\def\hfn{\hf_{(n)}}
\def\hfinf{\hf_{(\infty)}}

\def\hMult{\otimes}
\def\hFrc{\oslash}
\def\hSum{\oplus}
\def\hSub{\ominus}

\def\cQ{\mathcal Q}
\def\hpf{\mathcal P}
\def\bhpf{ {\bf \mathcal P}}
\def\hpfr{\hpf_{\pv}}
\def\hpfp{\hpf_{\ph}}

\def\met{d}
\def\fn{phantom}

\def\hi{\wp \,}

\newcommand{\su}[1]{\hat{#1}}
\newcommand{\Su}[1]{\widehat{#1}}
\newcommand{\du}[1]{\bar{#1}}
\newcommand{\phf}[1]{\ph\({#1}\)}

\newcommand{\Du}[1]{\overline{#1}}
\newcommand{\com}[1]{{#1}^{\operatorname{c}}}

\newcommand{\rpff}[2]{a_{#1}(#2)}
\newcommand{\ppff}[2]{b_{#1}(#2)}

\newcommand{\norm}[1]{\| #1 \|}
\newcommand{\abs}[1]{ \left|  #1  \right|}

\def\Sig{\Sigma}

\def\gs{\zeta}
\def\gp{\got{p}}
\def\gq{\got{q}}

\def\fX{f_{_X}}

\def\fXY{f_{_{X,Y}}}

\def\gX{g_{_X}}

\def\FX{F_{_X}}

\def\tfg{{\tilde f_{_X,\gm}(t)}}
\def\tfgt{{\tilde f_{_X,\gm}}}

\def\gmX{\gm_{_X}}

\def\gmY{\gm_{_Y}}
\def\pX{p_{_X}}
\def\pY{p_{_Y}}
\def\pXY{p_{_{X,Y}}}

\def\pXr{p_{_X,\pv}}
\def\pXp{p_{_X,\ph}}

\def\SX{S_{_X}}
\def\SY{S_{_Y}}

\def\MX{M_{_X}}

\newcommand{\MXi}[1]{M_{_{X_{#1}}}}
\def\MY{M_{_Y}}
\def\MW{M_{_W}}

\def\sigX{\sig_{_X}}
\def\sigXabs{\sig_{_{\abs{X}}}}
\def\sigMnabs{\sig_{_{\abs{M_n}}}}
\def\muY{\mu_{_Y}}
\def\muX{\mu_{_X}}
\def\muXabs{\mu_{_{\abs{X}}}}
\def\muMnabs{\mu_{_{\abs{M_n}}}}

\begin{document}
\title[Phantom  Probability] {Phantom  Probability}

\author{Yehuda Izhakian}

\address{Faculty of Management, Tel Aviv University, Ramat Aviv, Tel Aviv 69978, Israel}
\email{yud@post.tau.ac.il}

\author{Zur Izhakian}\thanks{The second author acknowledges  the
support of the Chateaubriand scientific post-doctorate fellowship,
Ministry of Science, French, 2007-2008. }

\address{Department of Mathematics, Bar-Ilan University, Ramat-Gan 52900,
Israel} \address{ \vskip -6mm CNRS et Universit´e Denis Diderot
(Paris 7), 175, rue du Chevaleret 75013 Paris, France}
\email{zzur@math.biu.ac.il}


\date{October  2008}




\begin{abstract} Classical probability theory supports probability
measures, assigning a fixed positive real value to  each event,
these measures are far from  satisfactory in formulating real-life
occurrences. The main innovation of this paper is the introduction
of a new probability measure, enabling varying probabilities that
are recorded by ring elements to be assigned to events; this
measure still provides a Bayesian model, resembling the classical
probability model.

By introducing two principles for  the possible variation
  of a probability  (also known  as uncertainty, ambiguity, or imprecise probability),
together with the ``correct'' algebraic structure allowing the
framing of these principles, we present the foundations for the
theory of phantom probability, generalizing  classical probability
theory in a natural way. This generalization preserves many of the
well-known properties, as well as familiar distribution functions,
of classical probability theory: moments, covariance, moment
generating functions,    the law of large numbers, and the central
limit theorem are just a few of the instances demonstrating the
concept of phantom probability theory.

\end{abstract}

\maketitle


{\small \tableofcontents}


\section*{Introduction}
\numberwithin{equation}{section}

Over the years much effort has been invested in trying human
beings have tried to understand aspects of probability in which
the evaluations of occurrences, as well as their likelihoods of
happening, are uncertain. Although the terminology for this type
of phenomena is varied (uncertainty for physicists, ambiguity for
economists, imprecise probability for mathematicians, and phantom
for us), fundamentally, the absence of theory enabling the
formulation of such phenomena is a common problem for many fields
of study. In this paper we introduce a new approach, supported by
a novel probability measure, allowing a natural mathematical
framing of this type of problems.

 Two  main principles underlie our approach to treating probability measures associated with
 varied evaluations:
\begin{itemize}
    \item For each event, the sum of its probability and its possible distortion lies in the
    real interval
    $[0,1]$; \pSkip
    \item The overall distortions always sum up to $0$.
\end{itemize}
Having the right algebraic structure,  termed  here the ring of
\textbf{phantom numbers} that naturally records  probabilities and
their oriented variations, these principles lead to the
introduction of our new \textbf{phantom probability measure}, on
which much of the theory of classical probability can be
generalized. This generalization captures both the uncertainty of
outcomes and ambiguous likelihoods, and it is still Bayesian.

 The ring $\hf$ of {phantom
numbers} consists of elements of the form $z = a + \hi b$, each of
which is a compound of the real \term \ $a$ and the phantom \term
\ $b$ (notated, like the complex numbers, by $\hi$ instead of
$i$), and whose operations, addition and multiplication
respectively, are
\begin{equation*}
    \begin{array}{lll}
      (a_1 + \hi b_1) \hSum  (a_2 + \hi b_2) & := & (a_1+a_2) + \hi
      (b_1+b_2),
      \\[2mm]
      (a_1+ \hi b_1) \hMult  (a_2+ \hi b_2) & := & a_1 a_2 +  \hi( a_1 b_2 +  b_1  a_2 +   b_1 b_2). \\
    \end{array}
\end{equation*}
This arithmetic makes $\hf$ suitable for the purpose of carrying a
theory of probability. In many ways this ring resembles  the field
of complex numbers, but its arithmetic is different; here $\hi$ is
idempotent, i.e. $\hi^ 2= \hi$, while  $i^2 = -1$ for the
complexes.  Similar structures, though sometimes using different
terminology, have been studied in the literature, mainly from the
abstract point of view of algebra; the innovation of this paper is
the utilization in probability theory, which requires some special
setting like phantom conjugate, reduced elements, absolute value,
and norm. With these notions suitably defined, the way toward the
development of a phantom probability theory is prepared.

One of the main advantages of phantom functions $f: \hf \to \hf$,
mainly polynomial-like functions, is that they can be rewritten as
$$ f = f_\pv + \hi( \su{f} - f_\pv),$$
where $f_\pv$ and $\su{f}$ are real functions $\Real \to \Real$.
We call this property, which plays a main role in our exposition,
the \textbf{realization property} of phantoms functions.

With this realization property satisfied, most of the phantom
calculations are reduced simply to the real familiar calculations.
Moreover, for $z= a + \hi b$, the real \term \ $a$ of $z$ is the
only argument involved in $f_\pv$; this shows that when  $G$ is a
pantomization of a real function $g: \Real \to \Real$, the real
\comp \ $G_\pv$ of $G$ is just $g$. Surprisingly, the
pantomizations of all classical probability functions (moments,
variances, covariances, etc.)  admit the realization property.


Using the phantom ring structure, together with our measure
principles, we keep track of the evolution of the classical theory
of probability. The leading motif throughout  our exposition is
that restricting the theory to the real \term s of all the
 arguments involved  always leaves ones  with the well-known classical theory. Given this foundation, as well as  the
appropriate definitions, the probability insights are much clearer
and their proofs become more transparent.

The main topics covered by this paper include:
\begin{itemize}
    \item Conditional probability, independence,  and Bayes' rule;

    \pSkip
    \item Random variables (discrete, continuous, and multiple);  \pSkip

    \item Attributes  of random variables:  moments, variances,
    covariances, moment generating functions; \pSkip

    \item Inequalities (appropriately defined); \pSkip
    \item  Limit theorems.
\end{itemize}
Along our exposition we also provide many examples demonstrating
how classical results naturally  carry over to the phantom
framework. further results and applications will be appear  in our
future  papers.

The fact that the phantom probability space provides a Bayesian
probability model paves the way for developing a theory of phantom
stochastic processes and phantom Markov chains \cite{youMarkov}
with a view towards applications in dynamical systems.

We use the  notion of \textbf{imprecise probability} as a generic
term to cover all mathematical models which measure chance or
uncertainty without sharp numerical probabilities \cite{Walley91}.
The known results of past efforts  to find a theory that frames
imprecise probability give only partial or complicated answers.
For example, fuzzy probability \cite{WangFuzzy} only treats
uncertain outcomes but not varying probabilities; conversely,
complex probability provides a partial answer for deformed
probabilities but only for fixed outcomes
\cite{Bijan,Youssef94quantummechanics}. On the other hand, the
operator measure theory \cite{Takesaki} is very complicated and
not intuitive, while the min-max model \cite{ GS1989} is not
Bayesian and \cite{S1989} sometimes becomes non-additive.

These probability theories have a tremendous range of
applications, like quantum mechanics, statistics, stochastic
processes, dynamical systems, game theory, economics, mathematical
finance, or  decision making theory; to name just a few. Our
development, together with the attendant  examples, which smoothly
extend
 the known theories that have already proven to be
significant, lead one to believe that phantom theory could
contribute to these applications, and make for a better
understanding of phenomenons that arise in the real life.


\section{The phantom ground ring}


\subsection{Ground ring structure}
The central idea of our new approach is a generalization of the
field $(\Real, +, \cdot \, )$ of real numbers to a ring structure
whose binary operations are induced by the familiar addition and
multiplication of $\Real$. Focusing on application to probability
theory, to make our exposition clearer, we give the explicit
description for the certain extension of $\Real$ of order $1$,
which is suitable enough for the scope of this paper. For the sake
of completeness, and in an effort to attract audiences from
various fields of study, we recall some of the standard algebraic
definitions (see \cite{Vinberg2003}) and present the full proofs
related to the basics of the algebraic structure for the extension
of order $1$. The more general phantom framework is outlined in
the next subsection.

Set theoretically, our ground ring $\hfto(\real)$, called a
\textbf{ring with phantoms} or \textbf{phantom ring}, for short,
is the Cartesian product $\Real \times \Real$; for simplicity, we
write $a + \hi b$ for a pair $(a,b) \in \hfto(\real)$. We say that
$\hfto(\Real)$ is a phantom ring of order $1$ over the reals, and
denote it as $\hf$, for short. (The general case of order $>1$ is
spelled out later in Subsection \ref{ssec:general}.) The elements
of $\hf$ are called \textbf{phantom numbers}, usually denoted $x,
y, z$.

In what follows we use the generic notation that $a,b\in \Real$
for reals and write $z:=a + \hi b$ for a phantom number $z$; we
call $a$ the \textbf{real \term } of $z$ while $b$ is termed the
\textbf{phantom \term } of $z$. We use  the notation
$$ \pv(z):= a  \qquad \text{and} \qquad  \ph(z) := b$$
for the real \term \  and the phantom \term \ of $z= a + \hi b $,
respectively. (The reason for calling the second argument
``phantom" arises from the meaning assigned to this value in the
extension of the probability measure, as explained  in
Section~\ref{sec:phProb}.)

The set  $\hf$ is then equipped with the two binary operations,
addition and multiplication, respectively,
\begin{equation*}\label{eq:hazyOperations}
    \begin{array}{lll}
      (a_1 + \hi b_1) \hSum  (a_2 + \hi b_2) & := & (a_1+a_2) + \hi
      (b_1+b_2),
      \\[2mm]
      (a_1+ \hi b_1) \hMult  (a_2+ \hi b_2) & := & a_1 a_2 +  \hi( a_1 b_2 +   b_1 a_2 +   b_1 b_2), \\
    \end{array}
\end{equation*}
to establish the \textbf{\fn \ ring} (to be proved next), $(\hf,
\hSum, \hMult)$, with unit $1 := 1 + \hi 0$ and zero $0 := 0 + \hi
0$. We write $\hf^\times$ for $\hf \setminus \{ 0 \}$, $\hi$ for
$\hi 1$, and $a - \hi b$ for $a + \hi (-b)$.

\begin{remark} In general, similar algebraic structures (with different terminologies) are known in the literature,
mainly for graded algebras or $k$-algebras in semiring theory,
usually applied to a tensor $M \otimes_k k$, where $M$ is a module
over semiring $k$, c.f., \cite{qtheory}. However, as will be seen
immediately, in this paper we push the algebraic theory much
further for the special case where $M = k$ is a field; then,  in
this case, $k \otimes_k k$ has a much richer structure. Moreover,
some of our definitions are unique with the aim of serving
applications in probability and measure theory.
\end{remark}

Although  the multiplication of $(\hf, \hSum, \hMult)$ is somehow
reminiscent of  the multiplication of the complex numbers $\Comp$,
it is different: for the phantoms $\hi = \hi^2$ is multiplicative
idempotent, while for the complexes $i^2 = -1$ is not idempotent.

\begin{proposition}\label{prop:mGroup}
$(\hf^\times, \hMult)$ is an Abelian semigroup.
\end{proposition}

\begin{proof} Given $z_i = a_i + \hi b_i$, where $i = 1,2,3$, we have
$$
\begin{array}{lll}
& (  z_1 \hMult z_2 )  \hMult z_3   & \\ [1mm]
&=    ( a_1 a_2 + \hi  ( a_1 b_2 +   b_1 a_2 + b_1 b_2)  )  \hMult ( a_3 + \hi  b_3)   \\[1mm]
& =(  a_1 a_2 )  a_3 + \hi  (  (  a_1 a_2 )  b_3 + (  a_1 b_2 + b_1 a_2  + b_1 b_2 )  a_3 + (  a_1 b_2 + b_1 a_2  + b_1 b_2 )  b_3 )   \\[1mm]
& =a_1 (  a_2  a_3 )  + \hi   (  a_1 (  a_2  b_3)  + a_1 (  b_2 a_3 ) + b_1 ( a_2 a_3)  + b_1 ( b_2 a_3)   +  a_1(  b_2 b_3)  + b_1(  a_2 b_3)  + b_1(  b_2  b_3)   )    \\[1mm]
& =a_1 (  a_2  a_3 )  + \hi  (  a_1 (  a_2  b_3 + b_2 a_3 +b_2
b_3) + b_1 ( a_2 a_3 )  + b_1( a_2 b_3 + b_2 a_3 + b_2  b_3)  )
\\ [1mm]
& = (a_1 + \hi b_1) \hMult (  (  a_2  a_3 ) +\hi  (  a_2  b_3 +
b_2 a_3 +b_2 b_3)  )   \\ [1mm] & = z_1 \hMult ( z_2 \hMult z_3 ),
\end{array}
$$
which proves associativity. Commutativity is obtained by
$$z_1 \hMult z_2   = a_1 a_2 + \hi \(a_1 b_2 + b_1 a_2  + b_1 b_2\)
= a_2 a_1 + \hi \(a_2 b_1 + b_2 a_1 + b_2 b_1 \)  = z_2 \hMult
z_1.
$$
This shows that $(\hf, \hMult)$ is a (multiplicative) Abelian
semigroup.
\end{proof}

\begin{theorem} $(\hf, \hSum, \hMult)$ is a commutative  ring.

\end{theorem}

\begin{proof} Since $\hSum$ is defined coordinate-wise, and $(\Real, +)$ is an (additive) commutative group,
 it is clear that  $(\hf,
\hSum)$ is also a commutative group. The unique additive inverse
$-z$ of $z = a + \hi b$ is
\begin{equation*}\label{eq:addInv} -z := (-a) + \hi (-b).
\end{equation*}
The pair $(\hf, \hMult)$ is a (multiplicative) Abelian semigroup,
by Proposition \ref{prop:mGroup}, so we need to prove the
distributivity of $\hMult$ over $\hSum$:
$$
\begin{array}{lll}
  z_1 \hMult ( z_2 \hSum  z_3) & = & (  a_1 + \hi b_1 ) \hMult ( a_2 + \hi b_2 +a_3 + \hi  b_3 )
  \\[1mm]
   & = & (  a_1+ \hi b_1 ) \hMult (  ( a_2 +a_3 )+ \hi (  b_2 + b_3) ) \\[1mm]
& = & a_1 ( a_2 +a_3 ) + \hi ( a_1 (  b_2 + b_3) + b_1   ( a_2
+a_3 ) + b_1 (  b_2 + b_3) ) \\[1mm]
&= & a_1 a_2 + a_1 a_3 + \hi ( a_1 b_2  + b_1 a_2  + b_1  b_2 ) +
\hi ( a_1 b_3  + b_1 a_3  + b_1  b_3) \\[1mm]
& = & (z_1 \hMult z_2) \hSum (z_1 \hMult z_3).
\end{array}
$$
All together we have proved that $(\hf,\hSum, \hMult)$ has the
structure of a commutative ring.
\end{proof}

Note that $(\hf^\times, \hMult)$ is not a group, and thus $(\hf,
\hSum, \hMult)$ is not a field, since there are non-zero numbers
$z \in \hf^\times$ without an inverse; for example $z = 0+\hi b$.

Recalling that a nonzero ring element $z_1$ is a \textbf{zero
divisor} if there exists a nonzero element $z_2$ such that $z_1
\hMult  z_2 =0$, one observes that the phantom ring $(\hf, \hSum,
\hMult)$ is not an integral domain, i.e. it has zero divisors; for
example
$$ (0 + \hi b) \hMult (-b + \hi b) = 0 + \hi( 0 b + b(-b) + b b) = 0$$
and thus $0 + \hi b$ and $-b + \hi b$ are zero divisors.

\begin{proposition}\label{prop:zeroDiv}
All the zero divisors of $(\hf, \hSum, \hMult)$ are of the form
\begin{equation}\label{eq:zeroDiv}
 z = 0 + \hi a \qquad \text{or} \qquad z = (-a) + \hi a,
\end{equation}
for some $a \in \Real$.
\end{proposition}

\begin{proof}
Assume $z_1 = a_1 + \hi b_1$, and $ z_2 =  a_2 + \hi b_2$ are
nonzero elements such that  $z_1 \hMult z_2 = 0$, that is
$$ (a_1 + \hi b_1) \hMult( a_1 + \hi b_1)  = a_1 a_2 + \hi(a_1 b_2 + a_2 b_1 + b_1 b_2 ) =0. $$
Suppose $a_1 \neq 0$,  then by the real \term \ of the product
$a_2 = 0$. So, by the phantom \term, we should have $a_1 b_2 +
 b_1 b_2 = (a_1 +
 b_1) b_2 = 0 $. But $b_2 \neq 0$, since $z_2 \neq 0$, and thus $b_1 =
 -a_1$. This means that $z_1 = a_1 - \hi _a1 $ and $z_2 = 0 + \hi
 b_2$m  as required.
\end{proof}

 A nonzero element $z \in \hf^\times$ which is not of the form
\eqref{eq:zeroDiv} is called a \textbf{nonzero divisor}; the
collection of all zero divisors in $(\hf, \hSum, \hMult)$ is
denoted
$$ \zd(\hf) = \{ \ z \in \hf \ | \ z  \text{ is zero divisor } \}.$$
We sometimes  write $\zdz$ for the union $\zd \cup \{ 0 \}$.

\begin{definition}\label{def:du}
The \textbf{phantom conjugate} $\du{z}$  of $z = a + \hi b$ is
defined to be
$$ \du{z} :=
(a +b) - \hi b.$$
The real number $$\su{z} := a + b$$ is called the (real)
\textbf{\top } of $z$.
\end{definition}

Having the notion of (real) reduction, we can write the product of
two phantom numbers as :
 \begin{equation}\label{eq:prud}
z_1 \hMult z_2 = a_1 a_2 + \hi( \Su{z_1} \Su{z_2} - a_1 a_2 ).
\end{equation}

\begin{remark}\label{rmk:zeroDiv}
By Proposition \ref{prop:zeroDiv}, one sees that $z = a + \hi b$
in $\hf$ is a zero divisor  iff $a=0$ or $\su{z} = 0$; when both
of them are zero then $z = 0$.
Moreover, in this view, given a suitable topology on $\hf$, the
complement of  $\zdz$ in $\hf$ is dense, so we can omit the zero
divisor without detracting form the abstract theory.
\end{remark}

One can easily verify the following properties for phantom
conjugates and (real) \top s:
\begin{properties} \label{prp:dusu} For any $z = a + \hi b$
the following properties are satisfied:
\begin{enumerate}
    \item $z = a + \hi(\su{z}-a) = (\su{z} - b) + \hi  b$, \pSkip
    \item $\Du{z_1 \hSum  z_2} = \du{z}_1 \hSum \du{z}_2$, \pSkip
    \item  $\Du{(-z)} = - (\du{z})$, \pSkip
    \item $\Du{z_1 \hMult  z_2} = \du{z}_1 \hMult \du{z}_2$, \pSkip
    \item $\ph(z \hSum \du{z}) = \ph(z \hMult \du{z}) = 0$, \pSkip
    \item $\du{z} = \su{z} - \hi b$, \pSkip
    \item $\Su{z_1 \hSum z_2} = \Su{z_1} + \Su{z_2}$, \pSkip
    \item $\Su{(- z)} = -(\su{z})$, \pSkip
    \item $\Su{z_1 \hMult z_2} = \Su{z_1} \cdot \Su{z_2}$, \pSkip

\end{enumerate}

\end{properties}

\begin{remark}\label{rmk:interval} The elements of  $\hf$ can be
understood as intervals in $\Real$. This means  that is an element
$z= a + \hi b$ stands for the interval that starts at $a$ and ends
at $a+b$, i.e. the reduction of $z$. Thus,  $\hf$ can be realized
as a ring of intervals, given by:
\begin{equation}\label{eq:intervals}
    \begin{array}{lll}
 [a_1, a_1+b_1] \hSum   [a_2, a_2+b_2]
 & = &  [a_1 + a_2, a_1 +  a_2+ b_1 +b_2];
 \\[1mm]
     [a_1, a_1+b_1] \hMult  [a_2, a_2+b_2] & = &  [a_1 a_2, a_1 b_2 + a_2 b_1+ a_2 b_2]. \\
    \end{array}
\end{equation}
In order to get a canonical interval representation, $z = a + \hi
b $ is assigned to the half-open interval $[a, \su{z})$.

In this view, zero divisors are intervals with $0$ as one of their
endpoints.  This view also provides the motivation for the
definition of the conjugate: $z$ and  $\du{z}$   represent the
same interval but with switched endpoints.

\end{remark}

 In fact, $(\hf, \hSum, \hMult)$ has a much richer structure
than a standard ring;  the division is well defined for all
nonzero divisors in $\hf^\times$, and each has an inverse.
 Given a nonzero divisor $z \in \hf^\times$,  we define the multiplicative inverse of $z$ to be
\begin{equation}\label{eq:inv}
z^{-1} := \ \frac{1}{a} + \hi  \frac{(-b)}{a \( a + b \)} \ = \
\frac{1}{a} + \hi  \frac{(-b)}{a \, \su{z}} \ = \ \frac{1}{a} +
\hi  \(\frac{1}{\su{z}} - \frac{1}{a} \) ;
\end{equation}
indeed $z^{-1}$ is an inverse of $z$,
$$\begin{array}{lllll}
   z  \hMult z^{-1} & = & (a + \hi b) \hMult \(\frac{1}{a} + \hi \frac{(-b)}{a \( a + b \)} \) & &
   \\[1mm]
    & = &  a \frac 1 a + \hi \( a \frac{(-b)}{a \( a + b \)} + \frac b a + b
\frac{(-b)}{a \( a + b \)} \) &= & 1.  \\
  \end{array}
$$
One can easily verify that $z^{-1}$ is unique, and that the
 reduction of an inverse number has the form:
$$\Su{z^{-1}} \  = \ \frac{1}{a+b} \ =  \  \frac{1}{\su{z}}.$$

Since the multiplicative inverse, defined only for all $z \in \hf
\setminus \zdz(\hf)$, and the additive inverse are unique, we
define  the division and the substraction, respectively, for
$(\hf,\hSum, \hMult)$  as
$$ z_1 \hFrc z_2 \ :=  \ z_1 \hMult z_2^{-1} \qquad \text{and} \qquad z_1 \hSub z_2 \ :=  \ z_1 \hSum (-z_2), $$
where $\hFrc$ is defined only for a nonzero divisor $z_2 \notin
 \zdz$. Accordingly, we  write
\begin{equation}\label{eq:frac} z_1 \hFrc
z_2 \ = \ \frac{a_1}{a_2} + \hi \frac{ b_1 a_2 - a_1 b_2}{a_2 \(
a_2 + b_2 \)},
\end{equation}
which leads to  the following useful form:
\begin{equation}\label{eq:frac2} z_1 \hFrc
z_2 \ = \ \frac{a_1}{a_2} + \hi \frac{ \Su{z_1} a_2 - a_1
\Su{z_2}}{a_2 \Su{z_2}} \  = \  \frac{a_1}{a_2} + \hi \( \frac{
\Su{z_1}}{\Su{z_2}} - \frac{a_1}{a_2} \).
\end{equation}



\begin{definition}\label{def:pseudoPositive}
A phantom number $z= a + \hi b$ is said to be \textbf{ positive}
 if
 $a > 0$ and $b > 0$. When $a > 0$ and $\su{z} > 0 $ we say that $z$ is
\textbf{pseudo positive}.
If $a < 0$ and $b < 0$, then $z$ is said to be \textbf{negative}
and when $a < 0$ and $\su{z} < 0$ we say that $z$ is
\textbf{pseudo negative}.
When $z$ is pseudo positive or $0$ it is termed \textbf{pseudo
nonnegative}, and if is pseudo negative or $0$ is called
\textbf{pseudo nonpositive}.
\end{definition}

Clearly, any positive (negative) phantom number is also pseudo
positive (negative). In particular, if $z$ is pseudo positive, or
pseudo negative,  then $z \notin \zdz$, cf. Remark
\ref{rmk:zeroDiv}, and is multiplicatively invertible.

\begin{lemma}\label{lem:positiveSquare} Given two pseudo nonnegatives $z_1, z_2 \in \hf$
then:
\begin{enumerate} \eroman
    \item Their sum is  pseudo nonnegative, \pSkip
    \item Their product is  pseudo nonnegative, \pSkip
    \item $z^2$ is pseudo nonnegative for each $z \in \hf$, \pSkip
    \item When $z_2$ is pseudo positive, the fraction $z_1 \hFrc
z_2$ is pseudo
    nonnegative.
\end{enumerate}
\end{lemma}

\begin{proof}  $ $ \eroman
\begin{enumerate}
\item Write $z_1 + z_2 = (a_1 + a_2) + \hi(b_1 + b_2)$, since  $a_1 +
a_2$ is positive, and $b_1 \geq -a_1$ and $b_2 \geq -a_2$, the
proof is clear.
 \pSkip

    \item By Equation \eqref{eq:prud} $z_1 z_2 = a_1
a_2+\hi(\Su{z_1} \Su{z_2}-a_1 a_2)$. Then, by the hypothesis, $a_1
a_2 \geq 0 $ and $\Su{z_1 z_2 } = a_1 a_2 + \Su{z_1} \Su{z_2} -a_1
a_2 =  \Su{z_1} \Su{z_2} \geq 0$. \pSkip
    \item Use (ii) with $z = z_1 = z_2$, or write directly $z^2 = a^2 + \hi(\su{z}^2 - a^2)$, so $a^2 \geq 0$ and
    then
    $\Su{(z^2)} =  a^2+ (a+b)^2 - a^2= (a+b)^2  \geq 0$. \pSkip

    \item Writing  $z_1 \hFrc
z_2$ as in Equation \eqref{eq:frac2}, $ \frac{
\Su{z_1}}{\Su{z_2}}$ and $\frac{a_1}{a_2}$ are (real) positives,
and thus
$$\frac{\Su{z_1}}{\Su{z_2}} - \frac{a_1}{a_2} = \frac{(a_1+b_1) a_2 - a_1 (a_2
+ b_2)  }{\Su{z_2}a_2}=  \frac{b_1 a_2 -  a_1 b_2 }{\Su{z_2}a_2}
\geq  \frac{-a_1 a_2 -  a_1 b_2 }{\Su{z_2}a_2} =  -\frac{a_1
}{a_2}.$$

\end{enumerate}

\end{proof}

Next, we outline the view of our structure in the category of
rings. Categorically, we have the trivial embedding
\begin{equation*}\label{eq:embbeding}
\varphi: (\Real, +, \cdot \,) \To (\hf,\hSum, \hMult),
\end{equation*} given by
sending $\varphi: a \mapsto a + \hi 0$.  On the other hand, we
also have the onto  projection
$$ \pi : (\hf,\hSum, \hMult) \To  (\Real, +, \cdot \,),$$
given by sending  $\pi : a + \hi b \mapsto \al a + \bt b$ for some
real numbers $\al$ and $\bt$. If $\bt = 0$, the projection is
phantom forgetful, i.e. $\pi:z \mapsto \al (\pv(z))$,  while
$\pi$ is real forgetful when  $\al =0$.

\begin{remark}\label{rmk:reals} Viewing $(\hf,\hSum,\hMult)$ as an $\Real$-module,
we define the scalar multiplication $\Real \times \hf \to \hf$ as
$$ r (a + \hi b) \  := \ \varphi(r) \hMult (a + \hi b),$$
for any $r \in \Real$, which is written as  $r(a + \hi b) = (ra) +
\hi (rb)$, for simplicity. Similarly, we write $\frac{a}{r} + \hi
\frac{ b}{r}$ for $\varphi(\frac {1} {r}) \hMult (a + \hi b)$.
\end{remark}

It is easy to check that the set of real numbers forms  a subfield
in the ring of phantom numbers~$(\hf,\hSum, \hMult)$, and the
phantom numbers whose real \term \ is zero establish an ideal in
$(\hf,\hSum, \hMult)$.


\subsection{Generalization }\label{ssec:general}

In the previous subsection we  described the extension of order
$1$ of the field of real numbers. For  completeness, we present
the general definition of a phantom ring of arbitrary order.

Given a field $\Fld$ of characteristic $\neq 2$, usually the field
$\Real$ of real numbers, the \textbf{phantom ring $\hfn(\Fld)$ of
order $n$}, or \textbf{$\bf n$-phantom ring}, for short,  is built
over the product $\Fld \times \cdots \times \Fld$ of $n+1$ copies
of $\Fld$ indexed $0,1, \dots,n$. Accordingly, the elements of
$\hfn(\Fld)$ are just $(n+1)$-tuples $(x_0, x_1, \dots, x_n)$ and
$(y_0,y_1, \dots, y_n)$ denoted, respectively, as $\bfx$ and
$\bfy$. $\hfn(\Fld)$ is then equipped with the following binary
operations, addition and multiplication, respectively:
\begin{equation}\label{eq:genArit}
\begin{array}{lll}
   \bfx \hSum \bfy & := & (x_0+ y_0, x_1 + y_1, \ \dots \ , x_n + y_n),  \\[1mm]
 \bfx \hMult \bfy & := & (x_0 y_0, 
  \ \dots \ ,
 x_i \bar y_{i-1} + y_i \bar x_{i-1} + x_i y_i,\  \dots \ ,
  x_n \bar y _{n-1} + y_n \bar x_{i-1} + x_n y_n),
  \end{array}
\end{equation}
where $\bar x_i = \sum_{j=0}^{i}x_j$ and $\bar y_i = \sum_{j=0}^i
y_j$.

(Note that the notation here is different from that used in the
previous subsection, in particular the $x_i$ and the $y_i$, $i
 \geq 1$  stand for the phantom \term s for the respective level.)

Numbering the copies of $\Fld$ sequentially, the first copy
$\Fld_0$ is considered as the real \prt \ of $\hfn(\Fld)$ while
$\Fld_i$, $i \geq 1$, is said to be the \textbf{phantom of level}
$i$ of $\hfn(\Fld)$. Note that $\hf_{(0)}(\Fld)$ is just $\Fld$,
which is a subfield of $\hfn(\Fld)$.

Having the operations rigorously defined for any $n \in \Net$,
using the arithmetic  defined in \eqref{eq:genArit}, we can push
$n$ to infinity and also define the $\infty$-phantom ring
$\hfinf(\Fld)$.

 In the sequel, for simplicity,
we apply  our development only to $\hfto(\Real)$, which as we have
said is denoted $\hf$, though we note that extends smoothly to any
$\hfn(\Fld)$, with $n
>1$, defined over a suitable field $\Fld$. Generalizing the future
definitions suitably to $n$, the $n$-phantom ring $\hfn(\Fld)$
carres also the same properties as $\hfto(\Fld)$, to be described
in the next sections.

\pSkip \pSkip \bfem{Notations:} For the rest of this paper,
assuming that the reader is familiar with the arithmetical
nuances, we write $z_1 + z_2$ for $z_1 \hSum z_2$, $z_1 - z_2$ for
$z_1 \hSub z_2$, $z_1 z_2$ for the product $z_1 \hMult z_2$,
$\frac{\ z_1 \ }{ \ z_2 \ }$ for the division $z_1 \hFrc z_2$, and
$z^n$ for $z \hMult \cdots \hMult z$ repeated $n$ times. The
phantom ring $(\hf,\hSum, \hMult)$ is denoted $\hf$, for short.

\subsection{Relations and orders} In the sequel,
mainly for the development of  phantom probability theory, we need
some relations that help to utilize the structure of $\hf$.

 To
make our paper reasonably self-contained, let us recall the
property of a binary relation on a set for being an order:

\begin{definition}\label{def:oset} A binary relation $\oleq$  is a \textbf{weak
order}  on a set $\oset$ if the following properties hold:
%

\begin{enumerate}  \eroman
    \item

Reflexivity: $s \oleq s$  for all $s \in \oset$; \pSkip


\item
 Transitivity: $s_1 \oleq s_2$ and  $s_2 \oleq s_3$ implies $s_1 \oleq s_3$; \pSkip

\item  Comparability (trichotomy law): for any $s_1, s_2 \in \oset$, either $s_1 \oleq s_2$  or $s_2 \oleq s_1$. \pSkip
\end{enumerate}

(When $s_1 \oleq s_2$ and $s_2 \oleq s_1$ we write $s_1 \oeq
s_2$.) A \textbf{weakly ordered set} is a pair $(\oset,\oleq)$
where $\oset$ is a set and $\oleq$ is a weak order on $S$. When
$S$ consists of phantom numbers we say that $(\oset,\oleq)$ is a
\textbf{phantom weakly ordered set}.

Adding the extra axiom:
\begin{enumerate} \setcounter{enumi}{2}  \eroman
\item
 Antisymmetry: $s_1 \oleq s_2$ and $s_2 \oleq s_1$  implies $s_1= s_2$; \pSkip
\end{enumerate}
the order  $\oleq$ is then a \textbf{total order}, or
\textbf{order}, for short, and is denoted as $\leq$.
\end{definition}
Clearly, $\oeq$ induces an  equivalent relation on $\hf$, the
classes of which are $\hf /_{\oeq}$,   and when $\oeq$ is a total
order $\oeq$ is replaced by full equality $=$.
 We use the notation $\square_\wk$ to distinguish this order,
mainly when writing $\ole$,  from the other orders used in the
sequel. Therefore, the symbol $<$ and $\leq$  always denote the
usual order of the real numbers.

Although, in general, many relations may serve as weak order on
$\hf$, in this paper we require  the weak order $\oleq$ to have
the following properties:
\begin{properties}\label{prp:orderProperties} $ $
\begin{enumerate} \eroman
    \item Compatibility  with the standard order of the reals, that is
$$ z_1 \oleq z_2 \ \iff  \ a_1 \leq a_2,$$
for any $z_1 = a_1 + \hi 0$ and  $z_2 =  a_2 + \hi 0$.
    \item  Compatibility  with the arithmetic  operations of
    $\hf$:
\begin{enumerate} \ealph
    \item if $ z_1 \oleq z_2$  then  $ z_1 + z_3 \oleq z_2 + z_3$,
    for any $z_3 \in \hf$; \pSkip
    \item if $ z_1 \oleq z_2$  then $ z_1  z_3 \oleq z_2  z_3$,
    for any pseudo positive $z_3 \in \hf$; \pSkip
    \item if $ z_1 z_3 \oleq z_2$  then  $ z_1 \oleq \frac{ z_2}{  z_3}$,
    for any pseudo positive $z_3 \in \hf$.
\end{enumerate}

\end{enumerate}
\end{properties}
\noindent Viewing $\hf$ as an Euclidian space, we usually assume
that all the elements that are $\oeq$ form a connectable  set.

\begin{example}
Viewing $\hf$ as $\Real \times \Real$, our main example for a
total order  on $\hf$ is the \textbf{lexicographic order} $\lxleq$
defined as
\begin{equation}\label{eq:lexorder}
   a_1 + \hi b_1 \ \lxleq  \ a_2 + \hi b_2  \ \iff \ \left\{%
\begin{array}{ll}
    a_1 <   a_2 , & { a_1 \neq a_2;} \\[1mm]
    b_1 \leq b_2  , & {a_1 = a_2;} \\
\end{array}%
\right.
\end{equation}
which is a total order satisfying the above conditions.
\end{example}

In the continuation, when writing  $\oleq$, we assume the weak
order $\oleq$ is provided with the set structure. The reader
should keep in mind that one  interpretation for an order which is
also total is the lexicographic order $\lxleq$.

\begin{remark}
Note that our definition of pseudo positivity, cf. Definition
\ref{def:pseudoPositive}, is independent of the given order
$\oleq$ on $\hf$.

%
\end{remark}


In the sequel, mainly for probability theory, we also use the
notation $\rleq$ for the (real) relation
 \begin{equation}\label{eq:realOrder} z_1 \rle z_2 \ \iff \
\re{z_1} < \re{z_2},
\end{equation} the other real relations $\req$,
$\rleq$, $\rge$, and $\rgeq$ are defined similarly.

We define the real-valued function $\af{ \ } : \hf \to \Real $,
with a real positive parameter $\al \in \real$, given by
\begin{equation}\label{eq:af}
\af{ \ }: a + \hi b \ \longmapsto \ a + \frac{b}{\al},
\end{equation}
and write $\af{z}$ for the image of $z \in \hf$ in $\Real$. Then,
$\af{ \ }$ determines the equivalence relation on $\hf$ given by
$$ z_1 \asim z_2 \ \iff \ \af{z_1} = \af{z_2},$$
and written $z_1 \aeq z_2$. (Note that in the special case when
$\al = 2$, by this definition, we always have $z \aeq \du{z}$.)
The quotient ring of $\hf$, taken with respect to $\af{ \ }$, is
denoted as $\aqut{\hf}$; clearly $\aqut{\hf} \cong \Real$.

In the same way, $\af{ \ }$ induces a weak order on $\hf$,
provided as
\begin{equation}\label{eq:aorder} z_1 \ale z_2 \ \iff
\ \af{z_1} < \af{z_2};
\end{equation}
and satisfying Properties \ref{prp:orderProperties}; the relations
$\aleq$, $\age$, and $\ageq$ are determined similarly.

These relations are very important for  advanced topics in phantom
probability theory and their applications, mainly discussed  in
the sequel papers \cite{youDecision,youMarkov}.


\subsection{Powers and exponents}
Writing $z^n$, with $n \in \Net$,  for the product $z \cdots
 z$ with $z$ repeated $n$ times,  for any $z = a + \hi
b$ we have
\begin{equation*}
    z^n = a^n + \hi  \sum_{i=1}^{n} \chos{n}{i} a^{n-i}b^{i};
\end{equation*}
as usual, $z^0$ is identified with the unit $1$. This form leads
to the following friendly formula:
\begin{equation}\label{eq:pwo}
    z^n =
    a^n + \hi (  (a + b)^{n} - a^n) =
    a^n + \hi (  \su{z}^{n} - a^n).
\end{equation}
Following accepted standards, we write  $z^{-n}$ for
$\frac{1}{z^n}$, and therefore get the extension to integral
powers of phantom numbers.

Equation \eqref{eq:pwo} plays a main role throughout our
development and, together with Equation \eqref{eq:prud}, leads to
the next important formula, which is used frequently in the
sequel: 
\begin{equation}\label{eq:pwo2}
  z_1^n \, z_2^m  =  a^n_1 a_2^m + \hi \(  {\Su{z_1}}^{n}{\Su{z_2}}^{m}-a^n_1
    a_2^m \).
\end{equation}
(To verify this equality, combine Equation  \eqref{eq:pwo} and
Equation \ref{eq:prud}.)

\begin{properties} Given  a phantom number $z \in \hf$, then:
\begin{enumerate}
    \item $z^i z^j = z^{i+j}$, \pSkip
    \item $\frac{z^i}{z^j} = z^{i-j}$,  \pSkip
    \item$\(z^i\)^j = \(z^j\)^i = z^{ij}$, \pSkip
\end{enumerate}
for any $i,j \in \Int$.
\end{properties}

Of course, one can take an arbitrary finite number of
multiplicands, $z_1, z_2, \dots, z_n$, and get recursively
\begin{equation}\label{eq:pwo3}
\begin{array}{lll}
  z_1^{i_1} z_2^{i_2} \cdots  z_n^{i_n}  & = & a_1^{i_1} a_2^{i_2} \cdots
  a_n^{i_n} + \hi \(  {\Su{z_1}}^{i_1} {\Su{z_2}}^{i_2} \cdots
  {\Su{z_n}}^{i_n} -a _1^{i_1} a_2^{i_2} \cdots   a_n^{i_n} \),
\end{array}
\end{equation}
for any $i_1, i_2, \dots, i_n \in \Int$.
\begin{definition}\label{def:relaizarion} When a phantom equation $\Eq$ can
be written in terms of two real equations, $\Eq_ \pv$ and
$\Su{\Eq}$, as
$$ \Eq = \Eq_\pv + \hi(\Su{\Eq} - \Eq_\pv),  $$
we say that $\Eq$ has a \textbf{realization form}, or
equivalently, that it admits the \textbf{realization property}.
\end{definition}
\noindent For that matter an equation might be an arithmetic
expression or a function, where $\Eq_\pv$ and  $\Su{\Eq}$ stand
respectively for the real and the reduction of each argument
involved in $\Eq$.

For example,  Equations \eqref{eq:pwo}, \eqref{eq:pwo2}, and
\eqref{eq:pwo3} above admit the  realization property. In the
sequel, we will see that many other familiar equations admit this
nice property. Having this property, as spelled out later for
probability theory, phantom results are induced by known results
for reals, which makes the development much easier.

Since the realization property us satisfied for each $z \in \hf$
and any natural power $n\in \Net$, cf. Equation \eqref{eq:pwo}, it
easy is to determine the \textbf{n'th root}, if it exists, of a
 phantom number $z = a + \hi b$ as:
\begin{equation}\label{eq:nroot}
\begin{array}{lll}
  \sqrt[n]{z}  & = &  \sqrt[n]{a} + \hi (\sqrt[n]{\su{z}}  -\sqrt[n]{a}), \\[1mm]
\end{array}
\end{equation}
where $\sqrt[n]{a}$ and $\sqrt[n]{\su{z}}$ are, respectively, the
real $n$'th roots of $a$ and $\su{z}$, and $n \in \Net$ is a real
positive number. Clearly, when $n$ is even, both $a$ and $\su{z}$
must be nonnegative.

In the usual way, we sometimes write $z^{\frac{1}{n}}$  for
$\sqrt[n]{z}$, and have the properties:
\begin{properties} Given  pseudo
nonnegative phantom numbers $z$, $z_1$, and  $z_2$   then:
\begin{enumerate}
    \item $\sqrt[n]{z_1} \sqrt[n]{z_2} = \sqrt[n]{z_1 z_2}$,
    \pSkip
    \item $\sqrt[n]{\frac{z_1}{z_2}} =
    \frac{\sqrt[n]{z_1}}{\sqrt[n]{z_2}}$ , for pseudo positive
    $z_2$, \pSkip
    \item $\sqrt[n]{z^m} = \( \sqrt[n]{z}\)^m  = \(
    {z^{\frac{1}{n}}}\)^m = z^{\frac{m}{n}}$,
\end{enumerate}
for any positive  $m,n \in \Net$.
\end{properties}

In the specific case when  $n=2$, clearly, each  pseudo
nonnegative phantom number $z = a + \hi b \in \hf$ has a
\textbf{square root}
\begin{equation}\label{eq:root}
\begin{array}{lll}
  \sqrt{a  + \hi  b}  & = &  \sqrt{a} + \hi (  \sqrt{a  +  b} - \sqrt{a}
     ) \\[1mm]
    &  = &  \sqrt{a} + \hi ( \sqrt{\su{z}} - \sqrt{a}  ). \\
\end{array}
\end{equation}
Actually, $\sqrt{a}$ in the equation stands for $\pm \sqrt{a}$;
therefore there always exists a \textbf{nonnegative square root}
of $a + \hi b $, i.e. a root whose real  and phantom \term s are
both nonnegative.

In the standard way, we define the \textbf{exponent} of an element
$z \in  \hf$ to be the infinite phantom sum
$$e^z = 1 +  z   +  \frac{z^{2}}{2!} +  \frac{z ^{3}}{3!} + \frac{z^{4}}{4!} + \  \cdots $$
\begin{proposition}\label{prop:exp}
Given $z$, $z_1$ and $z_2$ in $\hf$ then: \pSkip
\begin{enumerate}
     \item $e^0$ = 1, \pSkip
    \item $e^z = e^a + \hi \(e^{a+b}-e^a\) $, \pSkip

    \item $e^{z_1}  \  e^{z_1} = e^{z_1 + z_2}$, \pSkip
      \item $e^{z_1} /  \ e^{z_1} = e^{z_1 - z_2}$.
\end{enumerate}
\end{proposition}

\begin{proof} (1) is by definition.
(2) Expand $e^z$ and use Equation  \eqref{eq:pwo}, i.e.
$$
\begin{array}{lll}
e^{a+ \hi b} & = & 1 +  {(a+ \hi b)}  +  \frac{{(a+ \hi b)} ^{
2}}{2!} +
\frac{{(a+ \hi b)} ^{ 3}}{3!} + \  \cdots  \\[1mm]
& = &1+a +\frac{a^2}{2!}+\frac{a^3}{3!}+ \cdots  \\
& & \ \ + \hi \( \(1+{(a+ b)} +\frac{(a+ b)^2}{2!}+\frac{(a+ b)^ 3}{3!} + \cdots \) - \( 1+a +\frac{a^2}{2!}+\frac{a^3}{3!}+ \cdots \)\)
\\[1mm]
& = & e^a + \hi \(e^{a+b}-e^a\).
\end{array}
$$
(3) Using the identity in (2), write
$$
\begin{array}{llll}
e^{z_1} \ e^{z_2} & = & \( e^{a_1} + \hi \(e^{a_1+b_1}-e^{a_1}\) \) \  \( e^{a_2} + \hi \(e^{a_2+b_2}-e^{a_2}\) \)
\\[1mm]
& = & e^{a_1}e^{a_2} + \hi \( e^{a_1}\(e^{a_2+b_2}-e^{a_2}\)+ e^{a_2}\(e^{a_1+b_1}-e^{a_1}\) +\(e^{a_1+b_1}-e^{a_1}\)\(e^{a_2+b_2}-e^{a_2}\) \)   \\[1mm]
& = & e^{(a_1+a_2 +\hi(b_1 + b_2)}  \\[1mm]
& = & e^{z_1 + z_2}.
\end{array}
$$
(4) Straightforward from $(3)$ by taking $e^{z_1 + z_2} = e^{z_1 +
(-z_2)}$.
\end{proof}
Proposition \ref{prop:exp} (2) yields  the following  convenient
form,  i.e. the realization form, for the  phantom exponent:
\begin{equation}\label{eq:exp} e^z = e^a + \hi \(e^{a+b}-e^a\) =
e^a + \hi \(e^{\su{z}}-e^a\),
\end{equation}
often used  in   phantom probability theory.

Analogously to classical theory, for any pseudo positive $z \in
\hf$ we define the \textbf{logarithm} as
$$\log(z) = (z - \z1)  - \frac{(z - \z1)^{2}}{2} + \frac{(z - \z1)^{3}}{3} -
\frac{(z - \z1)^{4}}{4} + \ \cdots, $$
where $\z1 = 1 + \hi 1$, and prove that
$$
\begin{array}{lll}
  \log(z) & = & \log(a) + \hi(\log( a+b) - \log(a)) \\ [1mm]
   &
  = & \log(a) +
\hi(\log( \su{z}) - \log(a)), \\
\end{array}
$$
that is, the realization property for phantom logarithm.

\subsection{Phantom spaces}

Modules over the phantom ring, called \textbf{phantom modules},
are  just like  standard modules over rings \cite{Row1994}. For
the reader's convenience we state this explicitly:

\begin{definition}\label{def:module0} A phantom  \textbf{$\hf$-module}
$V$ is an additive group  $(V,\hSum ,\zero_V)$ together with a
scalar multiplication $ \hf \times V \to V$ satisfying the
following properties for all $z \in \hf $ and $v, w  \in V$:
\eroman
\begin{enumerate}
    \item $z (v \hSum w) =  z v \hSum z w;$ \pSkip
    \item $(z_1 \hSum z_2)v = z_1v \hSum z _2 v;$ \pSkip
    \item $(z_1  z_2)v = z_1 (z_2 v);$ \pSkip
    \item $1  v =v;$ \pSkip
    \item $0 v =\zero_V = z  \zero_V .$
  \end{enumerate}
 \end{definition}

The direct sum $\bigoplus_{j \in \tJ} \hf$ of  copies
    (indexed by $\tJ$) of the phantom ring $\hf$  is denoted as~$\hf^{(\tJ)},$
    with zero element ${ \bf 0 }  = (0, \dots, 0),$
   and is called the \textbf{phantom space}. When $\tJ = \{ 1, \dots, n \}$,
then the phantom space $\hf^{(\tJ)}$ is denoted as $\hf^{(n)}$ and
we say that  $\hf^{(n)}$ is an  $n$-phantom space. As element of
$\hf^{(n)}$ is just an $n$-tuple $(z_1, \dots, z_n)$ and is
denoted as $\bfz$.



Denoting the nonnegative real numbers as $\Real_+$, we recall the
standard definition of a norm, formulated for the $n$-phantom
space:

\begin{definition}\label{def:norm}
A \textbf{norm} on $\hf^{(n)}$ is a real-valued function $\norm{
\, \,}: \hf^{(n)} \to \Real_+$ that satisfies:
 \eroman
\begin{enumerate}
    \item $0 \leq \norm{\bfz} \in \Real$ and $\norm{\bfz} = 0$ iff
    $\bfz = \bf 0$, \pSkip
    \item $\norm{r  \bfz}= |r| \norm{\bfz}$ for each $r \in
    \Real$, \pSkip
    \item $\norm{\bfz' \hSum \bfz''} \leq \norm{\bfz'} +
    \norm{\bfz''}$,
\end{enumerate}
for any $\bfz, \bfz', \bfz'' \in \hf^{(n)}$.
\end{definition}

In what follows we use the \textbf{absolute value}, also called a
\textbf{modulus}, $\abs{ \ }: \hf \To \Real_+$  given by
\begin{equation}\label{eq:abs}
\abs{ a + \hi b}  = \sqrt{\(a+\frac{b}{2}\)^2 +  \(\frac{b}{2}
\)^2}.
\end{equation}
(When $z$ is only a real \term, i.e. $z = a + \hi 0$.  This
definition coincides with the familiar  absolute value of the
reals.)

\begin{proposition}\label{prop:abs} The absolute value
$\abs{ \ \, }$ as defined in Equation \eqref{eq:abs} is a norm on
$\hf$.
\end{proposition}

\begin{proof}
(i) and (ii) are immediate by definitions. To prove (iii), we show
that $\abs{z_1+z_2} ^2 \leq (\abs{z_1}+\abs{z_2})^2$. Expanding
both sides of this form, and letting $\al_i = a_i+\frac{b_i}{2}$,
$\bt_i = \frac {b_i}{ 2}$ for $i = 1,2$, we have
$$\(\al_1 + \al_2\)^2 + \(\bt_1 + \bt_2 \)^2 \leq  \al_1^2 +  \bt_1^2 +
2 \abs{z_1} \abs{z_2} + \al_2 ^2 + \bt_2^2.$$
Discard similar components on both sides and write $\abs{z_1}$ and
$\abs{z_2}$ explicitly to get
$$2\al_1 \al_2 + 2
\bt_1 \bt_2 \leq 2 \sqrt{\al_1^2 + \bt_1^2} \sqrt{\al_2^2 +
\bt_2^2}.$$
Canceling  the common multipliers and taking squares,  we have
$\(\al_1 \al_2 + \bt_1 \bt_2\)^2 \leq\(\al_1^2 + \bt_1^2\)
\(\al_2^2 + \bt_2^2\),$
and thus
$$2\al_1 \al_2 \bt_1 \bt_2 \leq
\al_1^2\bt_2 ^2  + \bt_1^2\al_2^2,$$
which implies $ 0 \leq \al_1^2\bt_2 ^2 - 2\al_1 \al_2 \bt_1 \bt_2
+ \bt_1^2\al_2^2 = \(\al_1\bt_2  - \bt_1\al_2\)^2$. This proves
property (iii) of Definition~\ref{def:norm}.
\end{proof}

Using the reduced form of phantom numbers,  and Properties
\ref{prp:dusu} (1), $\abs{z}^2$ can be written also as
$$\abs{z}^2 = {\( a+ \frac{
\su{z}-a}{2}\)^2 +   \(\frac{\su{z}-a}{2}\)^2} =
 \frac{a^2 + \su{z}^2}{2}.$$
Having a weak order satisfying Properties
\ref{prp:orderProperties} (i), one  also has $\abs{z} \ogeq 0$ for
any $z \in \hf$.

\begin{figure} \qquad

\begin{minipage}{0.45\textwidth}
\begin{picture}(10,170)(0,0)
\includegraphics[width=\FigWidth in]{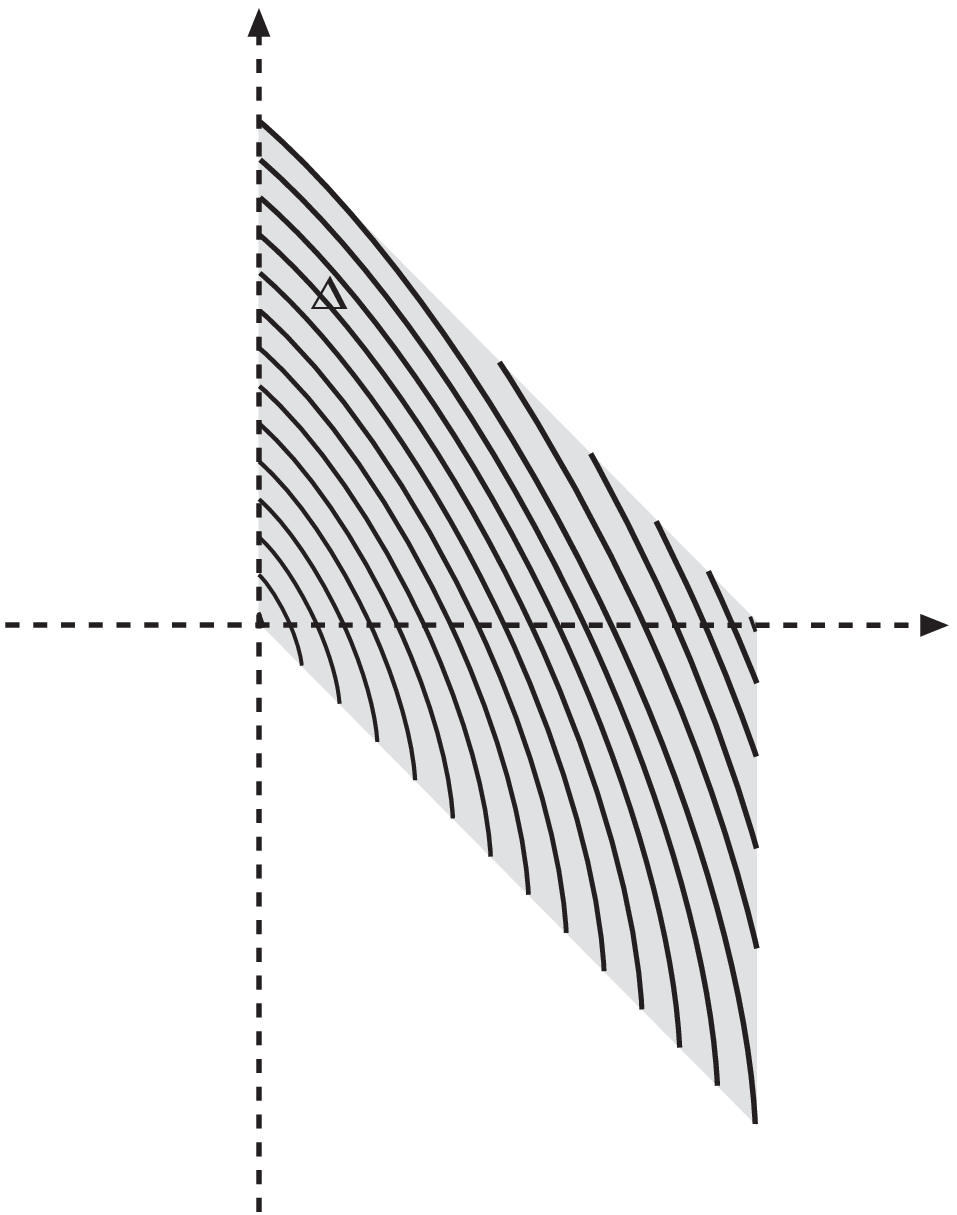}
\put(-95,160){$\hi$} \put(2,75){$\pv$}
\end{picture}
\center{(a)}
\end{minipage}\hfil
\begin{minipage}{0.4\textwidth}
\begin{picture}(10,170)(0,0)
\includegraphics[width=\FigWidth in]{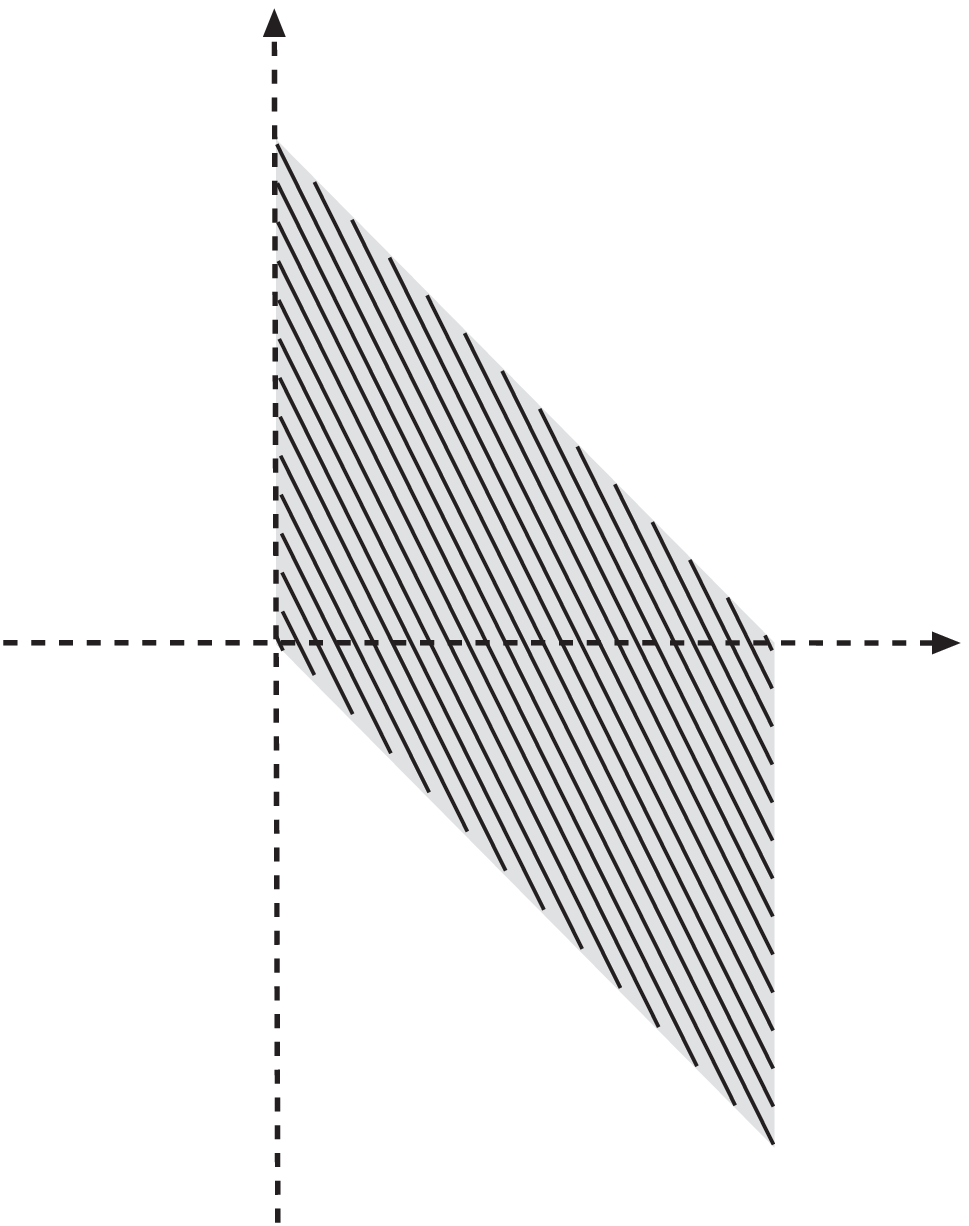}
\put(-95,160){$\hi$} \put(2,75){$\pv$}
\end{picture}
\center{(b)}
\end{minipage}
\caption{ (a) The iso-norm points on the compactification  $\cptr$
of $\hf$. (b) The $\af{ \ }$~-equivalent points, for $\al := 2$,
on $\cptr$. }\label{fig:norm}
\end{figure}

\begin{remark}\label{rmk:abs} There are several main reasons for
defining the  absolute value on $\hf$ as it has been defined in
Equation \eqref{eq:abs}:

\begin{enumerate} \eroman

    \item  $\abs{z} = \abs{\du{z}}$, for each
$z \in \hf$; indeed, to verify this identity, we have
$$
\begin{array}{lllll }
  \abs{\bar{z}} = \abs{(a + b) - \hi b } & = &  \sqrt{\(a+b -
\frac{b}{2}\)^2 +  \(\frac{b}{2}\)^2} &&\\[1mm]
   & = & \sqrt{\(a+ \frac{b}{2}\)^2
+  \(-\frac{b}{2}\)^2}  &= & \abs{z}. \\
\end{array}
$$

\item Considering $\hf$ as $\Real \times \Real$, the
two point compactification of each copy of $\Real$  is viewed as a
parallelogram with vertices $(0,0)$, $(0,1)$, $(1,0)$, and
$(1,-1)$ together with  the following correspondences:
$$ (\minf,\minf) \mapsto (0,0), \quad  (\minf,\infty) \mapsto (0,1),
\quad   (\infty,\minf) \mapsto (1,-1), \quad  (\infty, \infty)
\mapsto (1,0).$$ (See Figure \ref{fig:norm} (a).)

Accordingly, a $0 + \hi 0$ is the unique  point having absolute
value $0$, and  $1 + \hi 0$ is the unique point having  absolute
value $1$. The significance of this property become apparent later
in the discussion on the phantom probability measure.

\item The same view of $(ii)$, applied for  $\af{ \ }$ with $\al = 2 $, cf. Equation \eqref{eq:af},
shows that the classes $\af{ 0 }$ and $\af{ 1 }$  in $\hf/_{\al}$
are singletons. (See Figure \ref{fig:norm} (b).)
\end{enumerate}

These properties are very important for applications in phantom
probability theory.
\end{remark}

Having the norm $\abs{ \ }$ on $\hf$, we equipped $\hf$ with the
following relation:
\begin{equation}\label{eq:order}
    z_1 \ableq z_2 \iff \abs{z_1} <  \abs{z_2}.
\end{equation}
Accordingly, $0$ is the unique minimal element  in $\hf$.
(Clearly, this relation is also a weak order  on $\hf$; however,
since $\ableq$ ignores signs, it does not satisfy Properties
\ref{prp:orderProperties}.)

%
%

Basing on Equation \eqref{eq:abs}  we defined the \textbf{norm}
on the $n$-phantom space $\hf^{(n)}$ as:
\begin{equation*}\label{eq:norm}
\norm{ \bfz }  = \sqrt{ \abs{z_1}^2 + \cdots + \abs{z_n} ^2}
\end{equation*}
where $\bfz = (z_1,\dots, z_n)$.
\begin{proposition}\label{prop:norm} $\norm{ \ }$ is norm.

\end{proposition}

\begin{proof}
Straightforward from $\abs{ \ }$ being a norm.
\end{proof}

We use to $\norm{ \ }$ to define the map $\met: \hf^{(n)} \times
\hf^{(n)} \to \Real_+$ given by
 \begin{equation}\label{eq:metric}
\met: \bfz_1 \times \bfz_2 \longmapsto \norm{\bfz_1 - \bfz_2},
\end{equation}
where the substraction is taken coordinate-wise.
 We write
$\met(\bfz_1, \bfz_2)$ for the image of $ \bfz_1 \times \bfz_2$
under $\met$.
\begin{proposition}
$\met$ is a metric on $\hf^{(n)}$.
\end{proposition}

\begin{proof} By Proposition \ref{prop:norm}, we have $\met(\bfz_1, \bfz_2) \geq
0$, for any $\bfz_1$ and $\bfz_2$, and equals $0$ iff  $\bfz_1 =
\bfz_2$. Symmetry is clear. $\met(\bfz_1,\bfz_3) \leq
\met(\bfz_1,\bfz_2) + \met(\bfz_2,\bfz_3)  $ is derived from the
triangular law satisfied by $\norm{ \ }$.
\end{proof}

Note that using the metric \eqref{eq:metric} we always have
$$ \met(z, \du{z}) = \abs{a-(a+b) + \hi (b-(-b)) } =  \abs{-b + \hi 2b } = b,$$
for every $z = a + \hi b$ in $\hf$.

\begin{remark}\label{rmk:Borel}
The fact that $\hf^{(n)}$ is metric space allows us to define a
Borel $\sig$-algebra over $\hf^{(n)}$ in the usual way.
\end{remark}

%

\subsection{Polynomials} Polynomials over the phantom ring, called \textbf{phantom polynomials},
are defined just as formal polynomials over rings \cite{Row2006}.
As usual polynomials, say in $n$ phantom variables
$\lm_1,\dots,\lm_n$, form a ring which is denoted as
$\hf[\lm_1,\dots,\lm_n]$; these polynomials can also be viewed as
sums of polynomials in $2n$ real variables.

\begin{remark}\label{rmk:polynomials}
Given a polynomial $f  = \sum_i \al_i  \lm_1^{i_1}  \cdots
\lm_n^{i_n}$ in $\hf[\lm_1,\dots,\lm_n]$, it can be written as $f
= f_\pv + \hi f_\ph$, where  $f_\pv$ and $f_\ph$ are real
polynomials.

Suppose $\lm_i = u_i + \hi v_i$, $i = 1,\dots,n$, is a sum of two
variables $u$ and $v$ that take real values, and let
$$f_\pv(u_1, \dots, u_n)  =
\sum_i \re{\al_i}   u_1^{i_1}  \cdots  u_n^{i_n} \quad \text {and}
\quad
 \su{f}(\Su{\lm}_1, \dots, \Su{\lm}_n)  = \sum_i
\Su{\al}_i  \Su{ \lm}_1^{i_1}  \cdots  \Su{ \lm}_n^{i_n}$$
be two polynomials over the reals. (Note that, since $\Su{\lm}_i =
u_i + v_i$, $\su{f}$ is considered as a real polynomial in $2n$
variables.) Then, by Equation \eqref{eq:pwo2}, $f$ is written as
$$ f(\lm_1, \dots, \lm_n) = f_\pv(u_1, \dots, u_n) + \hi\(\su{f}(\Su{\lm}_1, \dots, \Su{\lm}_n) - f_\pv(u_1, \dots, u_n)\).$$
Therefore, phantom polynomials also admit the {realization
property}, in this case in the sense of functions.

\end{remark}

The \textbf{conjugate polynomial} $\du{f}$ of  $f  = \sum_i \al_i
\lm_1^{i_1} \cdots \lm_n^{i_n}$ is defined as $$\du{f}  = \sum_i
\Du{\al}_i \lm_1^{i_1} \cdots \lm_n^{i_n}.$$

\begin{proposition}\label{prop:congPoly} $\Du{f(z_1, \dots,
z_n )} = \du{f}(\du{z}_1, \dots, \du{z}_n)$ for any $ f \in
\hf[\lm_1,\dots,\lm_n]$ and each $(z_1, \dots, z_n) \in
\hf^{(n)}$.

\end{proposition}

\begin{proof}Straightforward by Properties \ref{prp:dusu}.
\end{proof}

\subsection{Basic analysis} Finally, we provide the necessary
notions for basic analysis over the phantoms; we present only the
general tools needed for our exposition. Most of these notions are
the phantom analogues to those in complex analysis; in general we
adopt the philosophy of analysis over the complexes.

\begin{definition}
Let $z_1, z_2, \dots$  be an infinite sequence of phantom numbers,
and let $z$ be another phantom number. We say that the sequence
$z_n$ \textbf{converges} to $z$, written $ \tLim{n}{\infty} z_n =
z $, if for every real $\ep > 0$ there exists some $n_0$ such that
$\abs{z_n - z} < \ep$, for all $n > n_0$.
\end{definition}

\begin{lemma}\label{lem:converge}
A sequence $z_i = a_i + \hi b_i$, $i= 1,2,\dots   $ converges to
$z = a + \hi b$ iff $ \tLim{n}{\infty} a_n = a$ and $
\tLim{n}{\infty} b_n = b$  as real sequences.
\end{lemma}

\begin{proof} $(\Leftarrow)$ Clear by definition, cf. Equation
\eqref{eq:abs}.

 $(\Rightarrow)$  Write $ \tLim{n}{\infty} \sqrt{\(a_n -a +\frac{b_n - b }{2}\)^2 +  \(\frac{b_n -b}{2}
\)^2} = 0$. Each, $\(\frac{b_n -b}{2} \)^2$  and $\(a_n -a
+\frac{b_n - b }{2}\)^2$ is positive and converges to $0$. Thus,
by the latter component, $b_n \to b$. Then, by the first
component, $a_n \to a$.
\end{proof}

A function $f: D \to \hf $, whose domain is a subset $D \subset\hf
^{(n)}$, is termed a  \textbf{phantom function}, while a function
$g: \Real^{(n)} \to \real$ is called a \textbf{real function}. We
say that a function is a phantom-valued function if its range lies
in $\hf$; similarly a function whose range lies in $\Real$ is
called real-valued,

\begin{definition}
Given a phantom function $f: D \to \hf $, we say that $w_0 \in
\hf$ is the \textbf{limit} of $f$ when $z \to z_0 \in D$ if for
any real $\ep> 0$ there exists a real $\dl
> 0$ such that for any $z$ with $\abs{z - z_0} < \dl$ we have
$\abs{f(z) - w_0} < \ep$. In such a case we write $\tLim{z}{z_0}
f(z) = w_0$.

\end{definition}

A function $f $ is \textbf{continuous} at $z_0 \in D$ if
$\tLim{z}{z_0} f(z) = f(z_0)$, and is said to be continuous on $D$
if it is continuous at each $z_0 \in D$.

Suppose $f : D \to  \hf$ is a phantom function, where $ D\subset
\hf$ is a set, and $z_0$ is an interior point of $D$. The
\textbf{derivative} of $f$ at $z_0$ is defined as
$$f'(z_0) = \tLim{z}{z_0} \frac{ f(z) - f(z_0)}{z - z_0}
,$$
provided this limit exists (depending also on $z - z_0$ being a
nonzero divisor). In this case, $f$ is called
\textbf{differentiable} at $z_0$. If $f$ is differentiable for all
points in an open disk centered at $z_0$ then $f$ is called
\textbf{analytic} at $z_0$. The phantom function $f$ is analytic
on the open set $D \subset \hf$  if it is differentiable (and
hence analytic) at every point in $D$. (The familiar properties of
derivation are also satisfied for phantom derivation.)

\begin{example}\label{exp:drPoly} The derivative of a polynomial $f  = \sum_i \al_i
\lm^i$ at $z_0$, where $\lm = u + \hi v $, written as $ f(\lm) =
f_\pv(u) + \hi(\su{f}(\Su{\lm}) - f_\pv(u))$ by Remark
\ref{rmk:polynomials}, is provided by using Equation
\eqref{eq:frac} as:
$$
\begin{array}{lll}
   f'(z_0) & = &  \tLim{z}{z_0} \frac{f_\pv(a)- f_\pv(a_0)}{a-a_0} +
\hi \frac{ (\su{f}(\Su{z})-f_\pv(a) - \su{f}(\Su{z_0}) +
f_\pv(a_0))(a-a_0) - (f_\pv(a)- f_\pv(a_0))(b-b_0) }
{(a-a_0)(a-a_0+ b - b_0) }  \\[1mm]
   & = &  f'_\pv(a_0) + \hi \tLim{z}{z_0} \frac{
(\su{f}(\Su{z})-\su{f}(\Su{z_0}))(a-a_0)}{(a-a_0)(\su{z} -
\Su{z_0}) } - \frac{(f_\pv(a)- f_\pv(a_0))(\su{z} - \Su{z_0})}
{(a-a_0)(\su{z} - \Su{z_0}) } \\[1mm]
& = & f'_\pv(a_0) + \hi\( \su{f}'(\Su{z_0}) - f'_\pv(a_0)\). \\
\end{array}
$$
\end{example}

When a phantom function has the realization property, i.e. it is
of the form  $f_\pv(t) + \hi f_\ph(t)$, where both $f_\pv$ and
$f_\ph$ are real functions with $t \in \Real$, the derivative of
$f$ is given as
\begin{equation}\label{eq:dirfrac}
f' = f'_\pv(t) + \hi f'_\ph(t).
\end{equation}
Indeed, write
$$ f' = \tLim{t}{t_0} \frac{ (f_\pv(t) + \hi f_\ph(t)) -  (f_\pv(t_0) + \hi f_\ph(t_0))}
{t - t_0}$$ which by Equation \eqref{eq:frac} is
\begin{equation*} f' = \tLim{t}{t_0}
\frac{f_\pv(t) - f_\pv(t_0)}{t- t_0} + \hi \tLim{t}{t_0} \frac{
 (f_\ph(t) - f_\ph(t_0)) }{\( t -t_0\)} .
\end{equation*}

Phantom \textbf{integration} is not really anything different from
real integration over pathes. For a continuous phantom-valued
function $ \phi(t): [a, b] \in  \Real \to \hf$, where $\phi =
\phi_\pv + \hi \phi_\ph$, we define
\begin{equation}\label{eq:integral1}
\int_a^b \phi(t) dt = \int_a^b \phi_\pv(t) dt + \hi \int_a^b
\phi_\ph(t) dt .
\end{equation}
For a function which takes phantom numbers as arguments, we
integrate over a path  $\gm$ (instead of a real interval) in $\hf$
realized as $\Real \times \Real$.  If one thinks about the
substitution rule for real integrals, the following definition,
which is based on Equation \eqref{eq:integral1} should come as no
surprise.

\begin{definition}\label{def:int}
Suppose $\gm$ is a smooth path parameterized by $\gm(t) : [a,b]
\to \hf$, $a \leq t \leq b,$ where $t, a,b \in \Real$, and $f :
\hf \to \hf$ is a phantom function which is continuous on $\gm$.
Then we define the \textbf{integral} of $f$ on $\gm$ as
\begin{equation}\label{eq:integral2}
\int_ \gm f(z) dz  = \int_ \gm f(\gm(t)) \gm'(t) dt.
\end{equation}
\end{definition}

This is simply the path integral of $f$  along the path $\gm$.
This integral can be defined analogously to the Riemann integral
as the limit of sums of the form $ \sum (f \circ \gm)(\tau_k)  (
t_k - t_{k-1})$, so is the Riemann–-Stieltjes integral of $f \circ
\gm$
 with respect to $\tau$. Using this definition, the integral can be extended to rectifiable paths,
i.e. ones for which $\gm$
 is only of bounded variation.

\begin{properties} Suppose $\gm$ is a smooth path, $f$ and $g$ are phantom functions which are continuous
on $\gm$, and $w \in \hf$ is constant.
\begin{enumerate}
    \item $\int_\gm (f +   w g) dz = \int_\gm f dz +  w \int_\gm g dz
    $. \pSkip
    \item If $\gm$  is parameterized by
$\gm(t)$, $a \geq t \geq b$, define the path $-\gm$
 through $-\gm(t) = (a + b - t)$,  $a \geq t \geq b$. Then $\int_\gm f dz = - \int_\gm
 f(z)
 dz$. \pSkip
    \item If $\gm_1$ and $\gm_2$ are paths so that $\gm_2$ starts where $\gm_1$ ends then
define the curve $\gm_1 \gm_2$ by following $\gm_1$ to its end,
and then continuing on $\gm_2$ to its end. Then $\int_{\gm_1\gm_
2} f(z) dz = \int_{\gm_1} f(z) dz + \int_{\gm_ 2}  f(z) dz$.
\end{enumerate}
\end{properties}

Assume $f$ is given as $f_\pv + \hi f_\ph$, and is defined along a
smooth path $\gm$,  given in a parametric form
\begin{equation*}\label{eq:parGm}
 \gm = \{ z = z(t) \ : \ \al \leq t \leq \bt, \  t \in \Real  \},
\end{equation*}
for some real $\al, \bt$,  for $z = z_\pv + \hi z_\ph $; we also
write $z_\pv = a(t)$ and $z_\ph= b(t)$. Then, using the familiar
line integral from calculus, Equation \eqref{eq:integral2} is
written in the following useful form:
\begin{equation}\label{eq:integral3}
\int_ \gm f(z) dz  = \int_\al^\bt f_\pv(t)a'(t) dt + \hi
\int_\al^\bt f_\pv(t)b'(t) +  f_\ph(t)(a'(t)+ b'(t))  dt,
\end{equation}
where $ f_\pv(t)$ and   $ f_\ph(t)$ stand respectively for
$f_\pv(a(t), b(t))$ and $f_\ph(a(t), b(t))$.

\section{Phantom  probability spaces}\label{sec:phProb}

\subsection{Phantom  probability laws}
We first recall the necessary basics of  standard measure theory,
then we further extend these basics to obtain the phantom setting
that generalizes the familiar classical probability framework. We
use \cite{pCHU01a,pFEL68a,pFEL71a} as general references for
classical probability theory.

A \textbf{measure space} is a triple $(\Om,\Sig,\mu)$, where
$\Sig$ is a $\sig$-algebra of subsets over a set $\Om$ and  $\mu:
\Sig \to [0,\infty]$ is a real-valued function, called a
\textbf{measure}, that satisfies the properties:
\begin{enumerate} \eroman
    \item $\mu(\emptyset) =  0;$ \pSkip
    \item $\mu(\bigcup_{i=1}^\infty A_i) = \sum_{i=1}^\infty \mu
    (A_i)$ for any countable sequence $A_1, A_2, A_3, \dots$ of pairwise disjoint sets
    in~$\Sig$.
\end{enumerate}
A measure $\mu$ is  \textbf{monotonic} if $\mu(A_1) \leq \mu(A_2)$
for each  $A_1 \subseteq A_2$.

A \textbf{probability measure} is a measure with total measure one
(i.e. $\mu(\Om) = 1$), cf. \cite{Billingsley}; a
\textbf{probability space} $(\Om,\Sig,P)$ is a measure space with
a probability measure $\mu : = P$ that satisfies the additional
probability axiom
$$ P(A) \geq 0, \qquad \forall A \in \Sig.$$
When $(\Om,\Sig,P)$ is a probability space, $P(A)$ is said to be
the \textbf{probability} of $A$, $\Om$ is called the
\textbf{sample space} and its elements are called
\textbf{outcomes},  usually denoted as $\om_1, \om_2, \dots  $.  A
collection of possible outcomes is called an \textbf{event}. In
the sequel, mainly in the examples, we use the letter $P$ to
denote a standard probability measure, i.e. $P : \Sig \to [0,1]
\subset \Real$.

\pSkip \bfem{Terminology:} In what follows, when using the term
``standard'', or ``real standard'', we refer to  the known
classical results, based on the above (real) probability measure,
appearing in the literature on probability theory,
\cite{pCHU01a,pDUR96a}. \pSkip

Roughly speaking, our aim is to generalize the probability measure
$P:\Sig \to \real$ to a phantom-valued function $\hpf: \Sig \to
\hf$, whose  real \comp \ is a standard probability measure while
its phantom \comp \  satisfies  an extra axiom. One way to realize
this extra axiom, enforced only on the phantom \comp, is to
understand the phantom as a signed distortion (either positive or
negative) assigned to each evaluation of the
 probability measure. Therefore, given a fixed event, its probability
together with an arbitrary distortion should still be  positive
(in the standard sense) and should  not exceed $1$.

\begin{remark}
In the continuation  the sample space  $\Om$ need not be a
standard sample space, and is also generalized to a
\textbf{phantom sample space} -- a sample space consisting of
phantom elements, called \textbf{phantom outcomes}. In what
follows, the notation $\Om$ is also used  for a phantom sample
space, and we  use  the standard terminology of outcomes and
events, respectively, for elements and subsets of $\Om$.
\end{remark}

Recall that $\hf$ is assumed to be equipped with a weak order
$\oleq$, coinciding with the standard order on $\Real$, usually a
total order.  However, to ensure that  our formalism is abstract
enough, we formulate our setting in terms of a general phantom
weak order $\oleq$ on $\hf$.

A phantom-valued function
$$\hpf: \Sig \To \hf,$$
 is called a
\textbf{phantom probability measure} if it satisfies the following
axioms. We denote the real \comp \ and the phantom \comp \ of
$\hpf$ as $\hpfr$ and $\hpfp$, respectively, each being a
real-valued function, and write $\hpf = \hpfr + \hi \hpfp$:
\begin{axiom}[\bfem{Phantom probability
measure}]\label{axm:phantomProb} $ $ \eroman
\begin{enumerate}
\item Nonnegativity:  $0 \leq \hpfr(A) \leq 1$ for each $A \in \Sig$,  \pSkip

\item Normalization:  $\hpf (\Om)= 1,$ \pSkip

\item Additivity: $\hpf(A \cup B) = \hpf(A)+\hpf(B)$ for any pair of
disjoint events  $A$ and $B$ in $\Sig$, \pSkip

\item Phantomization: $-\hpfr(A) \leq \hpfp(A) \leq 1 -\hpfr (A)$ for
each $A \in \Sig$.
\end{enumerate}
(The order $\leq$ is the standard order of the real numbers.)
\end{axiom}
As one can see,  conditions (i)-(iii) are none other than  the
well known classical probability axioms,  referring to real \comp
\ of $\hpf$ (condition  (iii) is also imposed on $\hpf_\ph$),
while the extra axiom (iv) is enforced on the phantom \comp. (This
axiom can be equivalently written as $0 \leq \hpfr(A) + \hpfp(A)
\leq 1$.) Therefore, the real \comp \ $\hpf_\pv$ of any phantom
probability measure $\hpf$ is always a standard (real) probability
measure. These axioms  properly frame our earlier probability
principles.

 Let $\cptr \subset \hf$ be the set
\begin{equation}\label{eq:teta}
 \cptr= \{ z \in \hf \   | \ a \in [0,1], \  -a \leq  b \leq 1-a \},
\end{equation}
each of whose points has a real \term \ belonging  to real
interval $\cptr_{\pv} = [0,1] \subset \Real$ and a phantom \term \
limited to the interval $[-a , 1-a]$, conditional on the real
\term \ of the points.  The set $\cptr$ is  called the
\textbf{phantom probability zone},  all of whose elements are
pseudo positive.

\begin{remark}\label{rmk:delOrder} In order to define  our
probability theory appropriately, we need to enforce the following
requirement on the weak order provided with $\hf$:
\begin{equation}\label{eq:delOrder}
    0 \ \ole \ z \ \ole \  1,  \qquad \text{ for each } z \in \ptr.
\end{equation}

For example, $\aleq$ with $\al > 1$ (cf. Equation
\eqref{eq:aorder}), or $\ableq$ (cf. Equation \eqref{eq:order}),
are weak orders that satisfy this condition. The total order
$\lxleq$ (cf. Equation \eqref{eq:lexorder}) also admits this
property.
\end{remark}

For any phantom probability measure $\hpf$, one sees that always
$$ \hpf: \Sig \To \cptr,$$
cf. Axiom \ref{axm:phantomProb} (iv).  Therefore, the target of a
phantom probability measure always lies in $\cptr$.  Given a fixed
event $A$, we use  the Gothic letter $\gp$ to denote the image
$\hpf(A) = p + \hi q$ of $A$, indicating that the corresponding
phantom number belongs to $\cptr$, and thus stands for a phantom
probability value.

Note that a phantom probability measure is notated by a
calligraphic letter, while a standard measure is notated by a
capital letter.

To avoid nonzero annihilators, in the sequel exposition,  we
usually  restrict the target of the  phantom probability measure
to the set
$$ \ptr = \{ \ \gp \in \cptr \ : \ \gp \text{ is not a zero divisor } \},$$
which we call the \textbf{restricted phantom probability zone}.
(Note that $0 \in \ptr$ and the $\cptr$ is the topological closure
of $\ptr$.)
 In the remainder, unless
otherwise  specified, we always assume the probability values are
in $\ptr$, i.e. we exclude all the possible zero divisors in
$\cptr$. (In  view of Remark \ref{rmk:abs} (ii), $\cptr$ is
realized as the compactification of $\hf$, while the elements of
$\ptr$ are all pseudo nonnegative, cf. Definition
\ref{def:pseudoPositive}.)

\begin{lemma}\label{lem:converse} Suppose $\gp, \gp' \in \ptr$, then:   \eroman
\begin{enumerate}
    \item $(1 - \gp) \in \cptr$, \pSkip
    \item $\gp \gp' \in \cptr$. \pSkip
\end{enumerate}

\end{lemma}
\begin{proof}

    (i) Let $\gp = p + \hi q$, then $1- \gp = (1-p) + \hi(-q)$. Clearly,
as $p \in [0,1]$,  the real \term \ $(1-p) \in [0,1]$.  The
phantom \term \ should satisfy $$ - (1-p) \leq  -q \leq 1 - (1-p)
\ \ ( \ = p),$$ i.e. $ -p \leq  q \leq  (1- p)$, but is  given by
the assumption that $\gp \in \ptr$. \pSkip

    (ii) Let $\gp = p + \hi q$ and $\gp' =
p'+ \hi q'$; then  $\gp \gp' = p p'  + \hi (p q' + q p' + q q')$,
clearly $p p' \in \cptr_\pv$. Using Axiom \ref{axm:phantomProb}
(iv), write
$$
\begin{array}{lll}
  p (-p') + (-p) p' + (-p)( -p')  & \leq  &  p q' + q
p' + q q'\\[1mm]
   &  \leq &  p (1 -p') + (1-p) p' + (1 - p)(1-
p'), \\
\end{array}
$$
and expand to get
$$p (-p')  \ \leq  \
p q' + q p' + q q' \  \leq \ 1 - p p',$$
as desired.
%
%
\end{proof}

Given an element $\gp \in \cptr$, the element $(1-\gp)$, also in
$\cptr$, is regarded as the \textbf{phantom complement} of $\gp$
in~$\cptr$.

\begin{lemma}\label{lem:teta} The image of a phantom probability measure $\hpf$ is
well defined for phantom addition and multiplication, that is
$\hpf(A)+\hpf(B)$,  $\hpf(A) \hpf(B)$, and $\hpf(\com{A})$  are in
$\ptr$ for any $A, B \in \Sig$.

\end{lemma}
\begin{proof}
The addition is axiomatic, since $A \cup B \subseteq \Sig$,
implies $\hpf(A)+ \hpf(B)\in \ptr$; cf.  Axiom
\ref{axm:phantomProb} (iv).

For the multiplication, take  $\gp = \hpf(A)$ and $\gq = \hpf(B)$
and apply Lemma \ref{lem:converse} (ii).  Since $\hpf({A}) \in
\ptr$, by Lemma \ref{lem:converse} (i), $\hpf(\com{A}) = 1 - \gp$
is in $\ptr$.
\end{proof}

\begin{definition}\label{def:panProbSpace}
A triple $(\Om,\Sig,\hpf)$, where $\Sig$ is a $\sig$-algebra of
subsets of $\, \Om $ and $\hpf$ is a phantom probability measure,
is called a \textbf{phantom probability space}. Given an event $A
\in \Sig$, $\hpf(A)$ is said to be the \textbf{phantom
probability} of $A$.
\end{definition}

 As mentioned earlier, in the context of phantom probability
spaces, the phantom \term \ should be realized as a signed bounded
distortion, with respect to each event,  dispersed non-uniformly
over the probability space and it has total sum $0$.
%
Accordingly, the phantom probability measure can be understand as
a family of  real probability measures $P_t : \Sig \to [0,1]$,
each satisfying
\begin{equation}\label{eq:intrvalProb}
    P_t(A) \in [\hpf_\pv(A), \hpf_\pv(A) + \hpf_\ph(A)], \qquad
    \text{for any } A \in \Sig,
\end{equation}
(or $P_t(A) \in [\hpf_\pv(A)+ \hpf_\ph(A), \hpf_\pv(A) ]$, when $
\hpf_\ph(A)$ is negative).

We say that a phantom probability measure $\hpf' : \Sig \to \hf$,
\textbf{agree with} $\hpf : \Sig \to \hf$ if it satisfies,
\begin{equation*}
     [\hpf'_\pv(A), \hpf'_\pv(A) + \hpf'_\ph(A)] \subseteq [\hpf_\pv(A), \hpf_\pv(A) + \hpf_\ph(A)], \qquad
    \text{for any } A \in \Sig.
\end{equation*}
A real phantom probability measure $P: \Sig \to \Real$, is said to
\textbf{agrees with} $\hpf : \Sig \to \hf$ if it satisfies
Equation ~\eqref{eq:intrvalProb}.

 This argument provides the
basis for Axiom \ref{axm:phantomProb} (iv): the sum of a
probability and a distortion can not exceed the probability of a
whole $\Om$ and is never negative, since otherwise it would
violate the standard laws of probability.

\begin{remark}\label{rmk:zerProb}

In  light of the previous  paragraph, the former notions obtain
the following special meaning:
\begin{enumerate}
    \item A probability $\hpf(A) = p + \hi q$ and its conjugate $(p-q) - \hi
    q$ resemble the same likelihoods (in the usual sense) lying between $p$ and
    $p+q$. This is one of the reasons for specifying a norm in
    which $\abs{\hpf(A) } = |\Du{\hpf(A)} | $.

    When one wants to dismiss this similarity and to have a unique canonical representative,
     he can use the correspondence
$ q \rightsquigarrow [0, q)$,
    and therefore $  p + \hi q$ is understood as $\{ P \in [p, p+q ) \}.$
    (The same setting can be used for sample spaces as well.)

    \item Zero divisors in the image of $\hpf$, if they  exist, correspond
    to events whose likelihood might be equal $0$, i.e. $\hpf(A) = 0 + \hi
    q$ or $\hpf(A) = p - \hi p $, with $p,q \geq 0$; cf.
    Proposition \ref{prop:zeroDiv}.

   An exclusive case  is when nothing is known about
    the likelihood of an event $A$; this scenario is recorded by  $\hpf(A) = 0 + \hi 1$.

    \item Fixing $\hpf_\ph := 0$ for the
    phantom \comp \ of $\hpf$, one gets the standard probability  model.
    \pSkip

    \item Cases in which two probabilities  $\hpf(A)$ and $\hpf(B)$ are both
    phantom numbers but their sum is real might happen and mean  that the probability of $A \cup B$ is
    fixed, but the probability of the interior subdivision is
    uncertain.
\end{enumerate}

\end{remark}

\begin{remark}\label{rmk:topProb}
Any  \fn \ probability measure $\hpf: \Sig \to \hf$ is associated
 neutrally  with a (real) \textbf{\topd \  probability measure}
$$\Su{\hpf}: \Sig \To [0,1],$$ given by sending each $A \in \Sig$
to $\Su{\hpf}(A)$. Since  $ \hpf(\com{A}) = 1 - \hpf(A)$ for each
$A \in \Sig$ and $\hpf(\Om) =1$, it is easy to verify that
$\Su{\hpf}$ is a proper
 standard real probability measure.

 We  recall that the real \comp, $\hpf_\pv: \Sig \to [0,1],$  of the phantom measure $\hpf$ is a proper
 standard real probability measure as well.
\end{remark}
As will be seen in the sequel,
 this \topd \  probability measure plays a crucial role in our
 future development.


\subsection{Digression}


In  view of Subsection \ref{ssec:general}, the phantom probability
measure $\hpf: \Sig \to \hf$, is a certain case of a phantom
measure with  probability zone $\cptr \subset \hf_{(1)}(\Real)$ of
order $1$. The general case is given with the phantom probability
measure $$\hpf: \Sig \To \hf_{(n)}(\Real),$$ of order $n$,  and
the following generalization of Axiom \ref{axm:phantomProb} (iv):
\begin{enumerate} \eroman  \setcounter{enumi}{3}
    \item  $-\hpfr(A) \leq \sum_{\ell=1}^{i}\hpf_{\ell}(A) \leq 1 -\hpfr (A)$ for
each $A \in \Sig$ and $i = 1,\dots,n$,
\end{enumerate}
where $\hpf_{\ell}$ denotes the phantom \comp \ of $\hpf$ of level
$\ell$.

Generalizing our  definitions appropriately, most of the following
theory extends  smoothly to phantom measures of order $n$.

\subsection{Elementary properties of phantom
probability}

\begin{proposition}[\bfem{Basuc properties of phantom
probability I}]\label{prp:pamtonProb} Given a phantom probability
measure $\hpf$, the following properties are satisfied for each
$A$ and $B$ in $\Sig$: \earabic
\begin{enumerate}
 \item $  \hpf_\ph(\Om) = 0$, \pSkip
 \item $ \hpf(\emptyset)=0$, i.e. $ \hpf_\ph(\emptyset) =  {\hpf_\pv(\emptyset)} = 0$, \pSkip
 \item $ -1 \leq  {\hpf_\ph(A)} \leq 1$ for for each $A \in \Sig$,
 \pSkip
 \item $ 0 \leq \abs{\hpf(A)}  \leq 1$, \pSkip
 \item $ \hpf(\com{A}) = 1 - \hpf(A)$, \pSkip
 \item $ \hpf(A \cup B ) = \hpf(A) + \hpf(B) - \hpf(A \cap
 B),$ \pSkip
  \item $\hpf_\pv (A) \leq  \hpf_\pv(B) $  if $A \subseteq B \subseteq \Om$, \pSkip

 \item $ \hpf_\pv (A \cup B) \leq  \hpf_\pv (A) + \hpf_\pv(B).$

\end{enumerate}
\end{proposition}

\begin{proof} $ $
\begin{enumerate}
    \item By definition, cf. Axiom \ref{axm:phantomProb} (ii).
    \pSkip
    \item By Axiom \ref{axm:phantomProb} (iii),
    $\hpf(\Om \cup \emptyset) = \hpf(\Om) + \hpf( \emptyset)$. Thus, by Axiom \ref{axm:phantomProb} (i),
      $\hpf_\pv( \Om) + \hpf_\pv( \emptyset) =1 + \hpf_\pv( \emptyset) \leq 1
    $, namely $\hpf_\pv( \emptyset) = 0$. On the other hand, by
    property (1), $\hpf_\ph( \Om) + \hpf_\ph( \emptyset) = 0 + \hpf_\ph(
    \emptyset) = 0$, so $\hpf_\ph(
     \emptyset) = 0$. \pSkip

    \item Immediate by Axiom \ref{axm:phantomProb} (i) and Axiom \ref{axm:phantomProb}
    (iii). \pSkip

    \item Let $\hpf(A) = p + \hi q$, then $| p + \hi q|^2 = p^2 + pq +
    \frac{q^2}{2}$. Thus, since $q \leq 1-p$ and  $0 \leq  p  \leq
    1$,
    $$p^2 + pq + \frac{q^2}{2}  \leq  p^2 + p(1-p) + \frac{(1-p)^2}{2} = \frac{1 + p^2}{2} \leq \frac{1+1}{2}.$$
On the other hand, since $q \geq -p $,  $p^2 + pq +
    \frac{q^2}{2} \geq  p^2 + p(-p) + \frac{(-p)^2}{2} =  \frac{ p^2}{2} \geq 0$. \pSkip

\item Straightforward from property (1). \pSkip

\item Write $A \cup B = A \cup (A^c \cap B)$  and $B  =  (A \cap B ) \cup (A^c \cap
B)$. The additivity axiom yields $$P(A \cup B) = P(A) + P(A^c \cap
B) \qquad \text{and} \qquad P(B) = P(A \cap B) + P(A^c \cap B).$$
 Subtracting the second equality
from the first and rearranging terms, we obtain the required.

\item{ and (8)} are precisely the well known relations  for the
real probability measure $\hpf_\pv$.
\end{enumerate}
\end{proof}

One of the reasons for defining the absolute value as in
 Equation \eqref{eq:abs} is that $0 \in \ptr$ is the unique element
 with $\abs{z} = 0$ and $1$ is the unique element with $\abs{z} =
 1$. The equiv-norm elements in $\ptr$ are a restriction of ellipses
 centered around the origin,  to the first quadrant( see Figure
 \ref{fig:norm}).

The same reason also led  to defining  the relation $\ale$ as in
Equation \eqref{eq:aorder},  since then for each $z \in \ptr$ we
have
 $0 \leq \af{z} \leq 1$, where $0$ and $1$  are obtained uniquely,
 i.e. $\af{1} =  1$ and $\af{0} = 0$. Thus, their equivalent classes
 in $\ptr$ are singletons, and they are all the singletons in $\aqut{\ptr}$. All
 the equivalent classes are parallel line segments  having a slope $=
 -\al$.

Next, we plug in the given weak order $\oleq$ on $\hf$; note that
this weak order assumes satisfying the condition of Remark
\ref{rmk:delOrder}.

\begin{proposition}[\bfem{Elementary properties of phantom
probability II}]\label{prp:pamtonProb2} For any phantom
probability measure $\hpf$, the following properties are satisfied
for each $A$ and $B$ in $\Sig$: \earabic
\begin{enumerate}
 \item $ 0 \oleq  {\hpf(A)} \oleq 1$ for  each $A \in \Sig$,
 \pSkip

 \item $\hpf(A) \oleq  \hpf(B) $  whenever $A \subseteq B \subseteq \Om$, \pSkip

 \item $ \hpf(A) \oleq \hpf(A \cup B),$ for any pair of disjoint
 events
$A$ and $B$ in $\Sig$, \pSkip
 \item $ \hpf(A \cup B ) \oleq  \hpf(A) + \hpf(B).$

\end{enumerate}
\end{proposition}

\begin{proof} $ $
\begin{enumerate}
    \item $ {\hpf(A)} \in \cptr$ in which  $ 0 \oleq  z \oleq 1$ for each $z
    \in \cptr$, cf. Remark \ref{rmk:delOrder}.  \pSkip
    \item Write $B = A \cup C$. So, by Axiom \ref{axm:phantomProb} (iii),
      $\hpf(A\cup C) = \hpf(A) + \hpf(C)$. But $\hpf(C) \ogeq 0$ by property (1),  and hence
      $ \hpf(A) + \hpf(C) \ogeq \hpf(A) $.  \pSkip
    \item Immediate by property (2).\pSkip
     \item Combine property (1) and Proposition \ref{prp:pamtonProb} (6).
\end{enumerate}

\end{proof}

\begin{proposition}[\bfem{Compound phantom probability measure}] Suppose $\hpf_i$, $i = 1,\dots,m$, are phantom
probability measures and $z_i \in \cptr$ are phantom numbers  $z_i
= a_i + \hi b_i$ such that $\sum_i z_i = 1$.  Then $\hpf = \sum_i
z_i \hpf_i $ is also a legitimate phantom probability measure.
\end{proposition}

\begin{proof} We need to verify the axioms of being a phantom probability
measure; cf. Axiom \ref{axm:phantomProb}:
\begin{enumerate} \eroman
    \item For a fixed $A \in \Sig$, let $\hpf_{\min, \pv}(A) = \min\{ \hpf_{i,\pv}(A) \}$
    and let  $\hpf_{\max,\pv}(A) = \max \{ \hpf_{i,\pv}(A) \}$.
    Then,
    $$ 0 \leq  \ \hpf_{\min,\pv}(A) \sum_i a_i \ \leq \ \sum_i a_i \hpf_{i,\pv}(A)
     \  \leq \hpf_{\max,\pv}(A) \sum_i a_i \leq 1. $$
The fact that each  $z_i$ is  in $\ptr$ insures that $0 \leq  a_i
\leq 1$, and thus $a_i \hpf_i(A) \in [0,1] $. \pSkip

    \item $\hpf(\Om) = \sum_i z_i \hpf_i(\Om)  =  \sum_i z_i 1 = 1$.
    \pSkip
    \item By the additivity of each $\hpf_i$, since $A$ and $B$ are assumed to be disjoint, write
    $\hpf(A \cup B) =  \sum_i z_i \hpf_i(A \cup B) = \sum_i z_i \(\hpf_i(A) + \hpf_i(B)\) =
    \sum_i z_i \hpf_i(A) + \sum_i z_i  \hpf_i(B) = \hpf(A) + \hpf(B)
    $. \pSkip

    \item To prove that $ \hpf_\ph(A) \leq 1 -  \hpf_\pv(A) $, we
    expand
    $$ \begin{array}{lll}
           \hpf_\ph(A) & = &  \sum_i \( (a_i + b_i) \hpf_{i,\ph}(A) + b_i \hpf_{i,\pv}(A)\)
           \\[1mm]
                &  \leq  & \sum_i \( (a_i + b_i)(1-
    \hpf_{i,\pv}(A)) + b_i \hpf_{i,\pv}(A)\)\\ [1mm]
    & = & \sum_i (a_i + b_i) -  \sum_i a_i \hpf_{i,\pv}(A) \\ [1mm]
    & = & 1 - \hpf_{\pv}(A).
             \end{array}
$$
 The same argument shows
    that $- \hpf_\pv(A) \leq \hpf_\ph(A)$.
\end{enumerate}
\end{proof}

\begin{corollary} The space of phantom probability measures on a $\sig$-algebra
 $\Sig$ is closed under an action of  probability measures.
That is,
 given a family of  phantom probability measures $\hpf_i: \Sig \to \cptr$,
 $i=1,\dots,m$,
 and a phantom probability measure
 $\cQ: \{ \{ 1 \},\dots, \{m \} \}  \to
 \cptr$, then $\hpf$,
 defined as
 $$\hpf = \sum_i \cQ(i) \hpf_i,  $$
is also a phantom probability measure.
\end{corollary}

\subsection{Initial examples}
The following examples are presented mainly to  demonstrate  how
nonstandard problems are formulated naturally using phantom
probability models. Later we show the phantom analogues to
well-known probability distributions.

We start  with an example whose sample space is also a sample
space in the usual sense.
\begin{example} Consider an unfair coin whose probability $P_t(\head)$ to get
head  (in a single  experiment) is unfixed, but belongs to the
interval $[0.4, 0.6]$.  Accordingly, for any possibility of
$P_t(\head)$,  the probability $P_t(\tail)$ to get tail must
satisfy $P_t(\head)+ P_t(\tail) = 1$, and thus is also restricted
to the interval $[0.4, 0.6]$.

This situation is formulated phantomly by letting
$$ \hpf(\head) = 0.4 + \hi 0.2 \qquad \text{and} \qquad  \hpf(\tail) = 0.6 - \hi 0.2.$$
In this view, the real \term \ of $\hpf(\head)+ \hpf(\tail)$ is
constantly $1$, while the distortion, which is at most $0.2$, is
encoded in the phantom \term s of  $\hpf(\head)$ and
$\hpf(\tail)$.
\end{example}

In  classical probability theory the uniform probability is
defined by assigning an identical probability to each event $A$ in
$\Sig$, which recall is formulated as $  P(A_i) = \frac{ 1 }{ k }$
for a discrete model with $\Sig = \{  A_1,\dots ,A_k  \} $. This
trivial formulation becomes meaningless in the phantom framework,
since by Axiom \ref{axm:phantomProb} (iv) the phantom \term \ must
be identically $0$ for each $\hpf(A_i)$. But, in  view of Remark
\ref{rmk:zerProb} (1), one can alternate between phantom
probabilities and their  conjugates, unless $k$ is even, to have
the sum of  phantom \term s equal $0$.

Next we  consider an example with a phantom probability space.
\begin{example}
Assume a financial investment with an expected profit of $5M \$ $
up to $10 M \$ $ in the case of success, which is estimated to
have $40 \% - 60 \%$ likelihood and $0$ profit otherwise. Using a
phantom probability space we formulate this investment with $\Sig
= \{ \{ 5M + \hi 5M \}, \{ 0 \} \}$ where the phantom probability
measure $\hpf : \Sig \to \hf$ is given by
$$ \hpf:  \{ 5M + \hi 5M  \} \mapsto 0.4 + \hi 0.2, \qquad   \hpf: \{ 0 \} \mapsto 0.6 - \hi 0.2.$$
\end{example}

An exclusive case, very difficult to formulate using  classical
probability theory, is the following:
\begin{example}
Assume a gambler who  knows nothing about the chances of  winning
in a new roulette game. Denoting the event of wining and losing
respectively by $\win$ and $\loss$,  we define the phantom
probability measure:
$$ \hpf: W  \mapsto 0 + \hi 1, \qquad   \hpf: V \mapsto 1 - \hi 1.$$
Recall that these two phantom numbers are zero divisors in $\hf$.
\end{example}




\subsection{Conditional phantom probability and
Bayes' rule}

Analogously to  classical theory, the \textbf{conditional phantom
probability} of $A$, given a fixed conditioning event $B$, is
denoted $\hpf(A | B)$ and defined as
$$ \hpf(A | B) \ = \  \frac{\hpf(A  \cap  B) }{\hpf(B)},  $$
where  $\hpf(B)$ assumed nonzero and a nonzero divisor. As a
consequence, given two disjoint events $A$ and $B$, where $\hpf(B)
\notin \zdz$, we have the equality:
\begin{equation}\label{eq:conDisProb} \hpf(A | B) \, {\hpf( B)}
\ = \ {\hpf(A  \cap  B) }.
\end{equation}

 It can be verified that for a fixed event $A$, the conditional
phantom probability forms a legitimate phantom probability law
that satisfies Axiom \ref{axm:phantomProb}.  The fact that the
phantom \term \  of $\hpf(A | A)$  is    $0$ can be seen using
Equation \eqref{eq:frac}, that is
$$ \hpf(A | A) \ = \  \frac { p} {p} + \hi \frac{ p q  - q p}{p \( p +
q \)} \  = \  1 + \hi 0,$$
for  $\hpf(A) = p + \hi q$.


Assuming  all the conditioning events have probabilities that are
nonzero divisors and are $\neq 0$, applying Equation
\eqref{eq:conDisProb} recursively  we have
\begin{equation*}\label{eq:multRule}
 \hpf( \cap_{i = n}^n A_i) =  \hpf(A_1) \  \hpf(A_2 | A_1) \
 \hpf(A_3 | A_1 \cap A_2) \ \cdots  \  \hpf(A_n | \cap_{i=1}^{n-1}
 A_i).
\end{equation*}

A sequence of events $A_1,\dots,A_n \in \Sig$  is said to be a
\textbf{partition} of $\Om$ if each possible outcome  is included
in one and only one of the events $A_1, \dots, A_n$. That is, the
sample space $\Om$ is the disjoint union of the events $A_1,
\dots, A_n$.

\begin{theorem}[\bfem{Total probability theorem}]\label{thm:total}
Let $A_1, \dots , A_n$ be disjoint events that form a partition of
the sample space $\Om$  and assume that $\hpf(A_i) \neq 0$ is not
a zero divisor, for all $i = 1, \dots, n$. Then, for any event
$B$, we have
$$ \hpf(B) \  = \ \sum_i \hpf(B \cap A_i) \ = \ \sum_i \hpf(A_i)\ \hpf(B|A_i).$$
\end{theorem}
\begin{proof}
The events $A_1, \dots , A_n$ form a partition of the sample space
$\Om$, so the event $B$ can be decomposed into the disjoint union
of its intersections $A_i \cap B$ with the sets $A_i$. Using the
additivity axiom, Axiom \ref{axm:phantomProb} (iii), it follows
that $ \hpf(B) = \sum_i \hpf(B \cap A_i)$. The proof is completed
by the definition of conditional probability, i.e. $\hpf(B \cap
A_i) = \hpf(A_i)\hpf(B|A_i)$.
\end{proof}

\begin{theorem}[\bfem{Bayes' rule}]\label{thm:Bays}
Let $A_1, \dots , A_n$ be disjoint events that form a partition of
the sample space $\Om$  and assume that $\hpf(A_i) \neq 0$  is not
a zero divisor, for all $i = 1, \dots, n$. Then, for any event
$B$, we have
$$ \hpf(A_i | B) = \frac{\hpf(A_i) \ \hpf(B | A_i)}{\hpf(B)} =
 \frac{\hpf(A_i) \ \hpf(B | A_i)}{\hpf(A_1) \ \hpf(B|A_1) + \ \cdots \ + \hpf(A_n) \ \hpf(B|A_n)} .$$
\end{theorem}

\begin{proof}
To verify Bayes' rule, note that $\hpf(A_i)\hpf(B |A_i)$ and
$\hpf(A_i |B)\hpf(B)$ are equal, because they are both equal to
$\hpf(A_i \cap B)$. This yields the first equality. The second
equality follows from the first by using the total probability
theorem to rewrite $\hpf(B)$.
\end{proof}

\begin{remark} In  light of Theorems \ref{thm:total} and \ref{thm:Bays},
the phantom probability space,
as introduced in Definition \ref{def:panProbSpace} provides a
 \textbf{Baysian probability model}.
 Being a Baysian probability model is a crucial
 property of a theory of stochastic processes and Markov chains.
 The absence of this property is one of the
 main deficiencies in some alternative models that has been
 suggested in the past, cf.
 \cite{Bijan,GS1989,Youssef94quantummechanics}.
\end{remark}


\subsection{Independence}
When the equality
$$ \hpf(A | B)  \ =  \ \hpf(A)$$
holds, we say that the event $A$ is (phantomly)
\textbf{independent} of the event  $B$. Note that by the
definition $ \hpf(A|B) = \hpf(A \cap B)/ \ \hpf(B)$, this is
equivalent to
$$\hpf(A \cap B) = \hpf(A) \hpf(B).$$ We adopt this latter
relation as the definition of independence because it can be used
even if $\hpf(B)$ is a zero divisor or $0$, in which case $
\hpf(A|B)$ is undefined.

The symmetry of this relation also implies that independence is a
symmetric property; that is, if $A$ is independent of $B$, then B
is independent of $A$, and we can unambiguously say that $A$ and
$B$ are \textbf{independent events}.

We noted earlier that the conditional phantom probabilities of
events, conditioned on a particular event, form a legitimate
probability law. Thus, we can talk about phantom independence of
various events with respect to this conditional law. In
particular, given an event $C$, the events $A$ and $B$ are called
\textbf{conditionally independent} if
\begin{equation}\label{eq:cond3}
\hpf(A \cap B |C) \ = \ \hpf(A|C) \ \hpf(B |C).
\end{equation}
The definition of  conditional probability and the multiplication
rule yield
$$
\begin{array}{lllll}
  \hpf(A \cap B |C)  & = & \frac{ \hpf(A \cap B \cap C) } {\hpf(C)} & &
  \\ [1mm]
   & = & \frac{ \hpf( B \cap C)  \hpf(A | B \cap C) } {\hpf(C)} & &
  \\ [1mm]
   & = & \frac{ \hpf(C) \ \hpf( B | C) \  \hpf(A | B \cap C) }{\hpf(C)}
   & = & \hpf( B | C) \  \hpf(A | B \cap C),\\
\end{array}
  $$
and thus, using Equation \ref{eq:cond3},
$$ \hpf( B | C) \  \hpf(A | B \cap C) \ = \ \hpf(A|C) \ \hpf(B |C)$$

 After canceling the factor $\hpf(B |C)$, assumed nonzero
divisor and  $\neq 0$, we see that conditional independence is the
same as the condition $$\hpf(A|B \cap C) = \hpf(A|C).$$ In other
words, this relation states that if $C$ is known to have occurred,
the additional knowledge that $B$ also occurred does not change
the probability of $A$, even though it is a phantom probability
(understood as a varied probability in the classical sense).
Interestingly, like in the classical theory, independence of two
events $A$ and $B$ with respect to the unconditional probability
law, does not imply conditional independence, and vice versa.

We generalize the definition of phantom independence to finitely
many events, and say that the events $A_1,A_2,$ $\dots , A_n$ are
independent if
$$ \hpf(\cap_{i \in S} A_i) = \prod_{i\in S} \hpf(A_i),$$ for every subset $S$ of $\{1, 2, \dots , n\}$.

\section{Phantom random variables }
In many probabilistic models, a random variable is a real-valued
function $X: \Om \to \Real$ of the outcomes of an experiment,
which means that each outcome is assigned with a fixed single
(real) numerical value. In  real life this is far from being
satisfactory; for example consider that these numerical values
correspond to instrument readings or stock prices. Our module
allows the assignment of a varied numerical value, recorded as a
phantom number,  to each outcome.

Consider a random experiment with a sample space $\Om$. A
\textbf{phantom random variable}, written \textbf{\prv} for short,
$$X:
\Om \To \hf$$
is a single-phantom-valued  function of the form
\begin{equation}\label{eq:pRrad}
 X: \om \ \longmapsto \ x(t) = \rpff{x}{t} + \hi  \ppff{x}{t}, \qquad t \in \Real,
\end{equation}
that assigns a phantom number $x = X(\om)$,  called the
\textbf{value} of $X$, to each sample element $\om \in \Om$. We
write $X_\pv$ and $X_\ph$, respectively, for  the real and the
phantom \comp s of $X$. In this realization, $X$ is parameterized
by real numbers, denoted by $t$;  later we shall see that this
parametrization is either discrete or continuous.

Note that the terminology which used here is the traditional
terminology of probability theory, and for this reason we use the
letter $x$, which stands for $a_x + \hi b_x$,  to denote the
phantom evaluation of $X$ at $\om$, while $z$ stands for an
arbitrary element in $\hf$. (Clearly a \prv \ is not a variable at
all in the usual sense, but a function.)

The sample space $\Om$ is called the \textbf{domain} of the \prv \
$X$ and is denoted  $\domain{X}$.  The collection of all phantom
values of $X(\om)$, where $\om \in \Om$, is termed the
\textbf{phantom range}, or just \textbf{range}, for short, of the
\prv \ $X$ and is denoted by   $\range{X}$.  Thus, the range
$\range{X}$ of \prv \ $X(\Om)$ is a certain subset of the set of
all phantom numbers, usually assumed without zero divisors.

Note that two or more different sample elements might give the
same value of $X(\om)$, but two different numbers in the range
cannot be assigned to the same sample point.

\begin{remark}\label{rmk:torder} When the real parametrization in Map \eqref{eq:pRrad} is
one-to-one, i.e. $x(t_1) \neq x(t_2)$ for any $t_1 \neq t_2$, then
the parametrization induces a total order on $\range{X}$. We
denote this order as $\tleq$.

\end{remark}

 Clearly, any function $g: \range{X} \to \hf$ of
a \prv, i.e. a function whose domain  contains the range of $X$,
defines another \prv.

To any \prv \ $X$ we associate the \textbf{\topd \ phantom random
variable}, written \textbf{\trv} for short,
$$\Su{X}: \Om \To \Real$$
given by
\begin{equation}\label{eq:tpRed}
 \Su{X}: \om \longmapsto X_\pv(\om) + X_\ph(\om),
\end{equation}
which in  view of  Map \eqref{eq:pRrad} is $\Su{X}(\om) =
\rpff{x}{t} + \ppff{x}{t}$, with $t \in \Real$.
\begin{remark}\label{rmk:rRange}
Along our next development we use the weak order $\oleq$ on $\hf$,
assumed to satisfy condition \eqref{eq:delOrder}. Accordingly the
range of a \prv \ is well ordered. (Recall that  the main examples
are  $\aleq$ and $\ableq$ for a weak order, cf. Remark
\ref{rmk:delOrder}, and the lexicographic order, cf. Equation
\eqref{eq:lexorder} for a  total order.)

 We also remark
that the range $\range{X}$ of any given \prv \ $X$ can be embedded
in $\Real$, cf. Map \eqref{eq:pRrad}, and therefore, as pointed
out  earlier, can be parameterized by the real numbers.
\end{remark}

If $X$ is a \prv \ and $z \in \hf$ is a fixed phantom number, not
necessarily in $\range{X}$, we define the event $(X  = z)$ as the
preimage of $z$, i.e.
$$(X =  z) = \{ \om \in \Om: X(\om) =  z\},$$
which has probability $\hpf(X = z)$. Note that when $z \notin
\range{X}$, we set $\hpf(X = z) = 0$.

Similarly, for fixed numbers $z$, $z_1$, and $z_2$ in $\hf$, we
define the following events:
\begin{equation*}\label{eq:disFun}
\begin{array}{rll}
 (X \oeq z)  & = & \{ \om \in \Om: X(\om) \oeq z\} \\[1mm]
  (X \oleq z)  & = & \{ \om \in \Om: X(\om) \oleq z\} \\[1mm]
  (X \oge z)  & = & \{ \om \in \Om: X(\om) \oge z\} \\[1mm]
   (z_1 \ole X \oleq z_2)  & = & \{ \om \in \Om: z_1 \ole X(\om) \oleq z_2\}
   \\[1mm]

\end{array}\end{equation*}
which have respectively the phantom probabilities $\hpf(X \oeq
z)$, $\hpf(X \oleq z)$, $\hpf(X \oge z)$, and $\hpf(z_1 \ole X
\oleq z_2)$. (We emphasize that these values need not be real
numbers.)

Note that when $z \notin \range{X}$, we can still have $\hpf(X
 \oeq z) \neq 0$. Of course this can only happen for a weak order; for a total order $\hpf(X
= z) = 0$ for each $z \notin \range{X}$.

Given an arbitrary phantom number $z$ and a \prv \ $X$, with
one-to-one real  parameterizations $t$,   we define the function
  $\supX{  }:\hf \to \range{X}$ by
\begin{equation}\label{eq:supX}
    \supXv{z} \ := \ \max_{\tleq} \{ \max _ {\oleq} \{ x \in \range{X} \ :  \  x \oleq z\} \},
\end{equation}
and sometimes write $x^z$ for $\supXv{z}  \in \range{X}$. This
function is well defined unless $\range{X} = \emptyset$. Note that
the interior $\max$ provides a set of elements in $X$ which are in
the same equivalence class, determined by $\oeq$,   while the
exterior pick the maximal $t$-element in this class.

In the same way we define the function $\infX: \hf \to \range{X}$
as
\begin{equation}\label{eq:infX}
    \infXv{z} \ := \ \min_{\tleq}\{ \min _{\oleq} \{ x \in \range{X} \ :  \  x \ogeq z\} \},
\end{equation}
and write $x_z$ for $\infXv{z}  \in \range{X}$. As before, we set
$ \infXv{x} = x$ for each $x \in \range{X} $. We emphasize thar
the interior ``$\min$" and the ``$\max$" above are taken with
respect to $\oleq$ -- the weak order on $\hf$, and the exteriors
are taken with respect to $\tleq$ -- the total order on
$\range{X}$. The use of these functions is mainly for continuous
\prv \, as will be seen later.



\subsection{Discrete random variables}
\begin{definition}
A \prv \   $X$  is called \textbf{discrete} if its range
$\range{X}$, i.e. the set of values that it can take, is finite or
at most countably infinite. \end{definition} \noindent When  $X$
is discrete,  we sometimes denote the values of $X$ as $x_1, x_2,
\dots, x_k, \dots$, indicating that it is parameterized by $t \in
\Net$.

 The most important way
to characterize a random variable is through the (phantom)
probabilities of the values that it can take. For a discrete
random variable $X$, these values  are captured by the
\textbf{probability mass function} of $X$, written \textbf{\pmf}
for short, and  denoted $\pX$. In particular, if $x$ is any
possible value of $X$, the \textbf{probability mass} of $x$,
denoted $\pX(x)$, is the phantom probability of the event $\{X
\oeq  x\}$ consisting of all outcomes that give rise to a value of
$X$ equal to $x$. That is
\begin{equation*}\label{eq:mass}
 \pX(x) \ =  \  \hpf(X =   x),
\end{equation*}
and therefore  $ 0 \oleq \pX(x) \oleq 1$, cf. Proposition
\ref{prp:pamtonProb2} (1). The mass function of $X$ is extended to
the whole $\hf$ by setting $\pX(z) := 0$ for each $z \notin
\range{X}$, and thus
\begin{equation*}\label{eq:mass2}
 \pX(z) \ =  \ \left\{%
\begin{array}{ll}
    \hpf(X =   z), &   z \in \range{X} ;  \\ [1mm]
    0, &   z \notin \range{X}. \\
\end{array}%
\right.
\end{equation*}

By the axioms of the phantom probability
 measure, we therefore have
$$ \sum_{x \in X} \pX(x) = 1,   $$
where in the summation above, $x$ ranges over all the possible
numerical phantom values of $X$. This follows from the additivity
and normalization axioms, because the events $\{X =  x\}$ are
disjoint and form a partition of the sample space $\Om$, as $x$
ranges over all possible values of $X$. By a similar argument, for
any set $S$ of phantom numbers, we also have
$$ \hpf(X \in  S) =  \sum_{x \in S} \pX(x),   $$
where $X \in  S$ means the values of $X$ which are contained in
$S$. (This notation it is  a bit misleading, but is the
traditional notation.)

Let us summarize the properties of phantom mass functions:
\begin{properties} Properties of a \pmf \ $\pX$: \pSkip
\begin{enumerate}
    \item $0  \oleq  \pX(z)  \oleq 1$, \pSkip
    \item $ \pX(z)  \in \ptr$, and thus $0 \leq  \abs{\pX(z)} \leq 1 $, \pSkip
    \item $ \pX (z) = 0$ if $z \neq x_1, x_2, \dots,  $ for all  $x_i \in \range{X}$, \pSkip
    \item $ \sum_k \pX (x_k)  =1$.
\end{enumerate}
\end{properties}

In  view of Remark \ref{rmk:topProb}, with any $\pmf$ \  $\pX$ we
 associate the  \textbf{\topd \ probability mass function}, written \textbf{\tpmf} \ for
 short,
$$ \Su{\hpf}(X \in S) =  \sum_{x \in S} \Su{\pX}(x). $$

Given a phantom number, in particular a phantom probability $\gp
\in \ptr$, by Definition \ref{eq:abs},  we always have
\begin{equation}\label{eq:hatp}
    \Su{1-\gp} = 1 - \su{\gp}.
\end{equation}
Note that for each phantom number $\gp \in \ptr$ we have the
inclusions $\su{\gp} \in [0,1]$ and $1 - \su{\gp} \in [0,1]$  in
the real interval $[0,1]$. We use this property in the forthcoming
examples.

 The following examples demonstrate  how well-known probability mass
 functions are generalized naturally in the phantom framework. Moreover, these
 examples show that most of these probabilities have the realization property.
  We  keep the traditional notation and denote the parameter as $\lm$,
 though, here it takes phantom values.

\begin{example}[\bfem{The binomial \prv}]\label{emp:binomial}
A biased coin with ambiguous  probability is tossed $n$ times. At
each toss, the coin comes up a head with phantom probability $\gp
= p + \hi q$, and a tail with phantom probability $ \gq = 1-\gp$,
independently of prior tosses.

Let $X$ be the number of heads in the $n$-toss sequence. We refer
to $X$ as a binomial \prv \  with parameters $n \in \Net$ and $\gp
\in \hf$. The \pmf \ of $X$ consists of the binomial
probabilities:

$$ \pX(k) = \hpf(X = k) = \chos{n}{k}\gp^k\gq^{n-k}, \qquad k= 0, 1, \dots , n. $$ (Note
that here and elsewhere, we simplify notation and use $k$, instead
of $x_k$, to denote the discrete  values of integer-valued random
variables.)

The normalization property   $ \sum_x \pX(x) = 1$, specialized to
the binomial random variable, is written
$$\sum_{k=0}^n \chos{n}{k}\gp^k(1-\gp)^{n-k} =1.$$
To see that this property is satisfied, apply  Equation
\eqref{eq:pwo2} for $\gp^k(1-\gp)^{n-k}$, that is
$$\gp^k(1-\gp)^{n-k}= p^k(1-p)^{n-k} + \hi\(
 (\su{\gp})^k
(1-\su{\gp})^{n-k} -p^k(1-p)^{n-k} \).$$
Recalling that $ \sum_{k=0}^n \chos{n}{k} p^k(1-p)^{(n-k)} = 1$
for any real $p \in [0,1]$, cf. \cite{pCHU01a}, and $\su{\gp} \in
[0,1]$, we take the sum
$$ \sum_{k=0}^n \chos{n}{k} p^k(1-p)^{n-k} + \hi  \sum_{k=0}^n
\chos{n}{k}\( \su{\gp}^k (1-\su{\gp})^{n-k} -p^k(1-p)^{n-k}\)  = 1
+ \hi\(1-1\)$$ to obtain the desired.
\end{example}

\begin{example}[\bfem{The geometric \prv}]\label{emp:geoemetic} The geometric \prv \  is the number
$X$ of trials, each with phantom probability $\gp$, needed for
success the first time. Its \pmf \ is given by
$$ \pX(k) = (1 - \gp)^{k-1} \gp, \qquad k = 1, 2, \dots. $$
This is a legitimate \pmf . Indeed, use Equation \eqref{eq:pwo2}
to write
$$\pX(k) = ((1 - p) - \hi q)^{k-1} (p + \hi q) =  (1 - p)^{k-1}p + \hi\((1- \su{\gp})^{k-1}(\su{\gp})-(1
- p)^{k-1}p \).$$
Recalling that $\sum_{k =1}^{\infty} (1 - p)^{k-1}p =1$ for any
real $p \in [0,1]$,  we take the sum
$$ \sum_{k =1}^{\infty} \pX(k) =
\sum_{k =1}^{\infty} (1 - p)^{k-1}p + \hi \sum_{k =1}^{\infty}
\((1- \su{\gp})^{k-1}(\su{\gp})-(1 - p)^{k-1}p \) = 1 +
\hi(1-1),$$ to get the required.
\end{example}

\begin{example}[\bfem{The Poisson \prv}]\label{emp:Poisson}
A Poisson \prv \  takes nonnegative integer values. Its \pmf \ is
given by
$$\pX(k) \ = \ e^{-\lm} \frac{\lm^k}{ k!}, \qquad k= 0, 1, 2,  \dots ,$$
where $\lm$ is a pseudo positive phantom parameter characterizing
the
\pmf. It is a legitimate \pmf \ because %
$$\sum_{k=0}^{\infty} e^{-\lm} \frac{\lm^k}{ k!} \  = \
e^{-\lm}\( 1 + \lm + \frac{\lm^2}{2!} + \frac{\lm^3}{3!} + \dots
\) \ = \ e^{-\lm} e^{\lm} \  = \  1.$$
The latter equality is by Proposition \ref{prop:exp} (3).
\end{example}

\subsection{Continuous random variables} The case when  phantom random
variables are continuous is much delicate than the discrete case,
especially since paths are involved and their parametrization
needs to be included carefully in our formulation.

\begin{definition}\label{def:contDens}
A \prv \  $X$ is called \textbf{continuous} if its probability law
can be described in terms of a piecewise continuous phantom
function $\fX : \hf \to \hf$, called the \textbf{probability
density function} of $X$, written \textbf{\pdf} \ for short, whose
real \comp \ is nonnegative and which satisfies
\begin{equation}\label{eq:dInt}
 \hpf(X \in S) =  \int_{S} \fX (x)dx,
\end{equation}
for every subset $S$ of
\begin{equation}\label{eq:gmX}
\gmX = \{ X = \rpff{x}{t} + \hi \ppff{x}{t} \ | \ t \in \Real \},
\end{equation}
where $\rpff{x}{t}$ and  $\ppff{x}{t}$ are real piecewise
differentiable  functions.
\end{definition}
  In fact, we  care only about $\fX $ restricted to $\gmX$, which is
  the range of $X$,
on which $\fX $ is piecewise continuous. The set $\gmX$ is
realized as a path interval in $\hf = \Real \times \Real$,
isomorphic to an interval in $\Real$, and it plays  a main role in
our exposition.


Note that $S$ does not need not be continuous. In such a case,
assuming $S$ compounds of countably many continuous subsets $S_i$,
the integral is translated to the sum of integrals over the $S_i$,
i.e.
\begin{equation}\label{eq:intDecompose}
 \int_{S} \fX (x)dx = \sum_{S_i}
\int_{S_i} \fX (x)dx
\end{equation}
 where each $S_i$ is continuous and the $S_i$'s are pairwise disjoint. In
the sequel, for simplicity, we assume $S$ is  a continuous subset
of  $X$, otherwise  we apply the same consideration as
\eqref{eq:intDecompose}.

\begin{figure}[!h] \qquad

\begin{minipage}{0.4\textwidth}
\begin{picture}(10,150)(0,0)
\includegraphics[width=\FigWidth in]{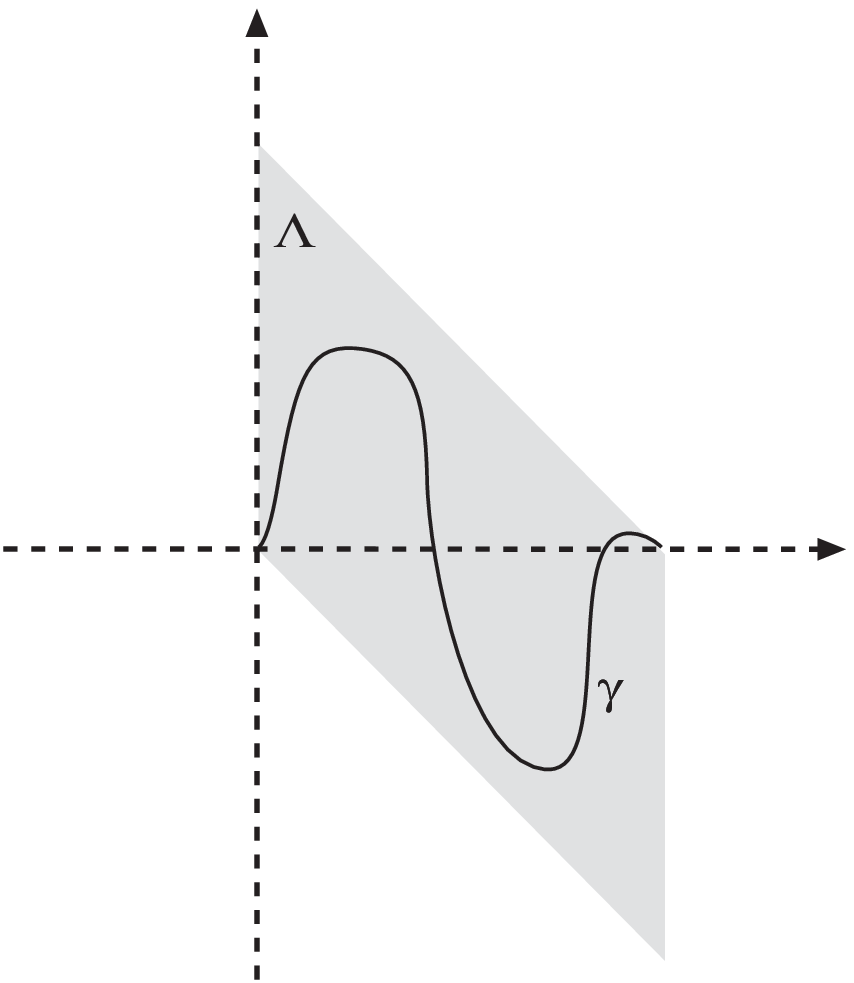}
\put(-90,145){$\hi$} \put(2,62){$\pv$}
\end{picture}
\center{(a)}
\end{minipage}\hfil
\begin{minipage}{0.4\textwidth}
\begin{picture}(10,170)(0,0)
\includegraphics[width=\FigWidth in]{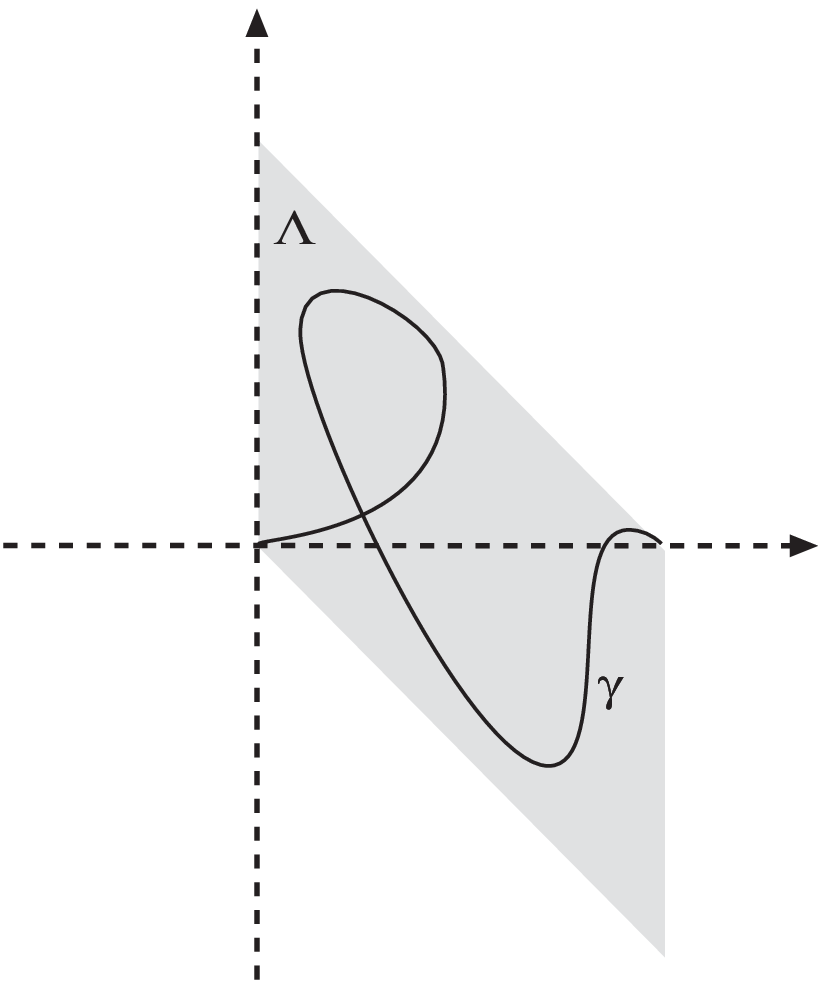}
\put(-90,145){$\hi$} \put(2,62){$\pv$}
\end{picture}
\center{(b)}
\end{minipage}
\caption{(a) Illustration of the compactification  of a path used
in most applications. (b) The compactification of a path $\gmX$
having self intersection points. }\label{fig:path}
\end{figure}

\begin{remark}\label{rmk:pathes}
\begin{enumerate} \eroman
    \item
 Since $\gmX$ is a parameterized
path in  $\hf$, i.e. $x(t) = \rpff{x}{t} + \hi \ppff{x}{t}$ for
each $x \in \gmX$, the map $\Real \to \gmX$ is not necessarily
one-to-one, and several reals may have the same image. In other
words $\gmX$ might have self-intersection points, see for example
Figure \ref{fig:path} (a). These cases require  special treatment
that is beyond  the scope of this paper. Therefore, in the rest of
this paper, when dealing with paths, they are always assumed to
have no self-intersections. \pSkip

\item
Having this assumption, i.e. $\gmX$ has no self-intersection
points, given a point $x = \rpff{x}{t} + \hi \ppff{x}{t}$ of
$\gmX$, for notional convenience, we write $\itt{x}$ for the real
$t$-value, determined  by the parametrization of $\gmX$,  whose
image is $x$.  Moreover, as mentioned before (cf. Remark
\ref{rmk:torder}), the given parametrization also determines an
order on $\gmX$, which we denote as $\tleq$, and write $x_1 \tleq
x_2$ when $\itt{x_1} \leq \itt{x_2}$. \pSkip

\item
In  view of Remark \ref{rmk:abs} (ii), in most applications the
image of a path $\gmX$ in the compactification of $\hf$ has the
endpoints $(0,0)$ and $(1,0)$, as illustrated by Figure
\ref{fig:path} (b), and usually does not have points with the same
real \term . However, we do not limit ourself to this type of
path. \pSkip

\item
Abusing  the notation, as in the complex convention,  in order to
address the situation that a continuous \prv \ $X$ is provided  as
a path $\gmX$, parameterized by $t\in \Real$, we sometime write
$\gmX(t)$ for $x(t)$ and $\gmX'(t)$ for the derivative $x'(t)$ of
$x(t)$ with respect to $t$, which by Equation \eqref{eq:dirfrac}
is just $a'(t) + \hi b'(t)$.
\end{enumerate}
\end{remark}

Since the integration of \pdf 's  is performed along paths whose
order does not need to  be compatible with the weak order on
$\hf$, we distinguish between cases in which the integration scope
is determined by points that belong to the path and  cases when
these points are arbitrary phantom points. We start with the
former case, then we extend it the latter case.

First, since we are dealing with an integral along a parameterized
path, cf. Equation \eqref{eq:dInt}, using  Equation
\eqref{eq:integral2} this integral can be written as
$$ \hpf(X \in S) = \int_{S} \fX (x)dx = \int_{S} \fX (\gmX(t))\gmX'(t)dt, $$
for any $S \subset \gmX$, assumed to be continuous.
 In particular, the probability that the value
$x$ of $X$ falls within a path interval of $\gmX$, whose endpoints
are $x_1$ and $x_2$,  is
\begin{equation}\label{eq:pdf1} \hpf(x_1 \tleq X \tleq x_2 ) =
\int_{\itt{x_1}}^{\itt{x_2}} \fX (\gmX(t))\gmX'(t)dt, \qquad \text
{for } x_1,
 x_2 \text{ are on  }  \gmX.
\end{equation}
Recall that by Equation \eqref{eq:integral3}, this integral
decomposes into two real integrals, one for the real \comp \ and
the second for the phantom \comp. The evaluation of each \comp \
can be interpreted as the areas confined between $\gmX$ and the
graphs of the corresponding function.

To simplify the  notation, we define
\begin{equation}\label{eq:edf} \tfg := \fX (\gmX(t))\gmX' (t)
\end{equation}
and rewrite Equation \eqref{eq:pdf1} as
\begin{equation}\label{eq:pdf2} \hpf(x_1 \tleq X \tleq x_2 ) = \int_{\itt{x_1}}^{\itt{x_2}} \tfg dt.
\end{equation}
(In fact $\gmX$ is determined by $X$, but we use this notation to
indicate that $\tfg$ depends on the path $\gmX$ in $\hf$.)

As in the discreet case, a \pdf \ $\tfgt$ is associated with the
\textbf{\topd \  probability density function}, written
\textbf{\tpdf} \ for short, defined as
\begin{equation}\label{eq:pdf2} \Su{\hpf}(x_1 \tleq X \tleq x_2 ) =
\int_{\itt{x_1}}^{\itt{x_2}} \Su{\tfgt}(t) dt,
\end{equation}
by taking the integral over the sum of the real \comp \ and the
phantom \comp \ of $\tfgt$.

 As in the standard theory, for any single value $x$ we have
$$ \hpf(X = x ) = \int_{\itt{x}}^{\itt{x}}
\tfg dt = 0.$$ Therefore, including or excluding the endpoints of
an interval in $\gmX$ has no effect on its probability:
$$ \hpf(x_1 \tleq X \tleq x_2 ) = \hpf(x_1 \tle X \tleq x_2 ) = \hpf(x_1 \tleq X \tle x_2 ) =
\hpf(x_1 \tle X \tle x_2 ),$$
for any $x_1, x_2 \in \gmX$.

As usual, to be qualified  as a \pdf, $\fX $ must  satisfy the
normalization property
$$ \int_{\gmX}
\fX  dx  \ = \ \int_{-\infty }^{\infty} \tfg dt \ = \ 1, $$
where its real \comp \ must take only nonnegative values, i.e.,
$\re{\fX (x)} \geq 0$ for every $x \in \gmX$. Accordingly, the
function $\tfgt$ satisfies
$$ \int_{-\infty }^{\infty} \re{ \tfg } dt = 1 \qquad \text{and} \qquad
\int_{-\infty }^{\infty} \phf{ \tfg }dt =  0; $$ that is,
normalization and vanishing of phantoms, respectively.

Now, we turn to  the cases in which the scopes of random variables
are determined  by arbitrary phantom values. To extend $\fX$ to
the whole $\hf$, we fix $\fX(z) := 0$ for each $z \notin \gmX$ and
define
\begin{equation}\label{eq:pdfz1} \hpf(z_1 \oleq X \oleq z_2 ) =
\int_S \tfg dt, \qquad S = \{ x \in \gmX \ : \ z_1 \oleq x \oleq
z_2 \}
\end{equation}
%
%
for any $z_1, z_2 \in \hf$. In the case that $S$ is not
continuous, the integral is decomposed into the  sum of countably
many integrals along the  path intervals as in
\eqref{eq:intDecompose}, taken with respect to $\tleq$ in the
positive direction.

 When $S$ is continuous,  Equation \eqref{eq:pdfz1} has the form
\begin{equation}\label{eq:pdfz2} \hpf(z_1 \oleq X \oleq z_2 ) =
\int_{\itt{\infXv{z_1}}}^{\itt{\supXv{z_2}}} \tfg dt,
\end{equation}
which is just a line integral along a path interval of $\gmX$.
Recall that $\supXv{z_1}$ and $\infXv{z_1}$, cf. Equations
\eqref{eq:supX} and \eqref{eq:infX}, provide the ``top'' point and
the ``bottom'' point  on $\gmX$ closest to $z_2$ and $z_1$,
respectively; these points are unique since $\tleq$ is a total
order applied on the weak order $\oleq$. In this sense, we capture
all the elements in $\range{X}$ which are less or
$\oeq$-equivalent to $z_2$ and greater or $\oeq$-equivalent to
$z_1$.

Accordingly, we also have,
\begin{equation}\label{eq:pdfz3} \hpf(X \oeq z ) =
\int_{\itt{\infXv{z}}}^{\itt{\supXv{z}}} \tfg dt,
\end{equation}
which is the integration along all elements in $\range{X}$,
assumed continuous, that are $\oeq$-equivalent to $z$. In this
view it is easy to see that we might have $\hpf(X \oeq z ) \neq
0$.

Let us outline some basic properties of probability density
functions:
\begin{properties} Properties of \pdf \ $\fX $:

\begin{enumerate}
    \item $0  \oleq  \fX(z)  \oleq 1$,  for each $z \in \hf$,  \pSkip
    \item $ \fX(z)  \in \ptr$, and thus $0 \leq  \abs{\fX(z)} \leq 1 $, \pSkip
    \item $ \fX (z) = 0$ if  $z \notin \gmX$. 
\end{enumerate}
\end{properties}

In most applications, the  random variable  is either discrete or
continuous, but if the \pdf \ of a \prv \  $X$ possesses features
of both discrete and continuous random variable's, then the random
variable $X$ is called a \textbf{mixed phantom random variable}.

\subsection{Cumulative phantom distribution function}

The  \textbf{cumulative phantom distribution function}, written
\textbf{\cpdf }  for short,  $\FX$ of a \prv \  $X$  provides the
probability $\hpf(X \oleq z)$,  i.e.
\begin{equation}\label{eq:CDF}
    \FX(z) = \hpf(X \oleq z) = \left\{%
\begin{array}{ll}
    \sum_{x_k \oleq z} \pX(x_k), & \hbox{$X$  discrete;} \\[2mm]
    \int_{S} \tfg dt, & \hbox{$X$  continuous.} \\
\end{array}%
\right.
\end{equation}
 for every $z$ in $\hf$. Here, $S$ is defined as $\{ x \in \gmX \ : \ x \oleq z \}$,
 cf. Equation \eqref{eq:pdfz1} with $z_1$ tending  to $-\infty$, where $-\infty$ stands for    $-\infty(1 +\hi)$.

Loosely speaking, the \cpdf \ $\FX(z)$ ``accumulates" phantom
probability "up to" the phantom value $z$. As in  classical
theory, most of the information about a random experiment
described by the \prv \ $X$ is recorded by the behavior of
$\FX(z)$.

Using Equation \eqref{eq:pdfz1}, if $X$ is a continuous random
variable, then
$$ \hpf( z_1 \ole X \oleq z_2  ) = \int_{S} \tfg dt = \FX(z_2) - \FX(z_1). $$
(The discrete analogue is  obvious.) To emphasize, although $\FX$
gets an argument that is a phantom number, it also depends on the
parametrization of $X$.

\begin{properties}\label{prp:comDenFunc} Writing $z \to \infty$, for $z = a + \hi b$
with $a \to \infty$ and $b \to \infty$, we have the following
properties:

\begin{enumerate}
    \item  $0 \oleq \FX(z) \oleq 1$, \pSkip

    \item  $0 \leq \abs{\FX(z)} \leq 1$, \pSkip
    \item $\FX(z_1) \oleq \FX(z_2)$ if $z_1 \oleq z_2$, \pSkip
    \item $\tLim{z}{\infty}  \FX(z) =   \FX(\infty) =  1 $, \pSkip
    \item $\tLim{x}{-\infty}  \FX(z) = \FX(-\infty) =  0 $, \pSkip

    \item $\tLim{z}{z_0^+}  \FX(z) =  \FX(z_0^+)$, where
     $z_0^+ =  \tLim{0 < \abs{\vep}}{0} z_0^+  + \vep
    $. \pSkip
\end{enumerate}
\end{properties}
\noindent (The verification  of these properties is
straightforward.)

Accordingly, having Properties \ref{prp:comDenFunc}, one can also
compute other probabilities, such as
\begin{equation}\label{eq:disFun2}
\begin{array}{rll}
  \hpf(X \oge z_1)  & = & 1 - \FX(z_1), \\[1mm]
   \hpf(X \ole z_2)  & = &  \FX(z_2^-), \qquad \text{where } z_2^-  =  \tLim{{0 {\oleq} \ep}}{0} (z_2 -
   \ep).
   \\[1mm]
\end{array}\end{equation}

\subsection{Moments and variance}

In the sequel, we use the notation $\mean{ \ }$, $\var{ \ }$, and
$\cov{ \ }$ respectively for moments, variance, and covariance.
These notations are used for both the phantom sense and the
standard sense (applied for real numbers) of the respective
probability functions, where the meaning is understood from the
context. We also point out that for these functions, and others,
the standard form ia always captured in the  real \comp \ of the
function. As will be seen, this attribute is provided for free by
the phantom structure, and one should keep in mind that the
generalization to the phantom framework is performed only through
the phantom \term s of the arguments. This is the leading idea of
our forthcoming exposition.

\begin{definition}
The $n$'th \textbf{moment} of a \prv \  $X$  is defined by
\begin{equation}\label{eq:moment}
    \mean{X^n} \ = \  \left\{%
\begin{array}{ll}
    \sum_x x^n \pX(x), & \hbox{$X$  discrete;} \\[2mm]
    \int_{\gmX} (\gmX(t))^n \tfg dt, & \hbox{$X$ continuous.} \\
\end{array}%
\right.
\end{equation}
The first moment $\mean{X^1}$ is called the \textbf{mean}, or the
\textbf{expected value},  of $X$ and is denoted by $\muX $.

\end{definition}

The $n$'th \textbf{\topd \  moment} $\mean{(\Su{X})^n}$ of a \prv
\ $X$ is the standard moment for reals, applied to $\Su{X}$ with
$\Su{\pX}$, or  $\Su{\tfgt}$, for  a discrete or a continuous $X$,
respectively. In the same way we define the $n$'th
\textbf{conjugate moment} $\mean{(\Du{X})^n}$ of $X$, by taking
$\Du{X}$, computed with respect to $\Du{\pX}$, or  $\Du{\tfgt}$;
clearly this is a phantom function.

We write $\meanr{X^n} $  and $\meanp{X^n}$, respectively,  for the
real \term \ and the phantom \term \ of $\mean{X^n}$ and therefore
have
$$\meanr{X^n} \ = \ \mean{{X_\pv}^n},$$
where $\mean{{X_\pv}^n}$ is taken with respect to the real \comp \
of ${\pX}$ or ${\tfgt}$. (This relation is not satisfied for the
phantom \comp, i.e. $\meanp{X^n} \neq \mean{{X_\ph}^n}$, since it
also involves  the real \term \ of the arguments.)

When $X$ is discrete, using Equation \eqref{eq:pwo2}, we write
$$
\begin{array}{lll}
  \mean{X^n} & =  & \sum_x x^n \pX(x) \\[1mm]
  & =  & \sum_x a_x^n p_x + \hi( (a_x+b_x)^n(p_x+q_x) - a_x^n
  p_x), \\[1mm]
   & =  & \sum_x a_x^n p_x + \hi( \su{x}^n\su{\gp}_x - a_x^n p_x), \\
\end{array}
$$
for  $x= a_x+\hi b_x$ and $\pX(x) = \gp_x =  p_x + \hi q_x$; a
similar form is also obtained  for a continuous $X$. Accordingly,
the phantom moment
 satisfies the realization property, and Equation \eqref{eq:moment}
gets the following friendly form:
\begin{equation}\label{eq:moment2}
   \mean{X^n} \ = \ \mean{{X_\pv}^n} + \hi(\mean{\Su{X}^n} - \mean{{X_\pv}^n}),
\end{equation}
where $\mean{{X_\pv}^n}$ stands for  $\sum p_x a_x^n $ and
$\mean{\Su{X}^n}$ stands for  $\sum \su{\gp} \su{x}^n$.

\begin{proposition}\label{prop:conMean}
     $\Du{\mean{X^n}} = \mean{\Du{X}^n}$, where $\mean{\Du{X}^n}$
    is computed with respect to the conjugates of $X$ and the
probability measure.

\end{proposition}
\begin{proof} Straightforward by the additivity and the
multiplicativity  of the phantom conjugate, cf.
Properties~\ref{prp:dusu}.
\end{proof}

The \textbf{variance} of a \prv \  $X$, denoted by $\sigX ^2$ or
$\var{X}$, is defined as
\begin{equation}\label{eq:var1} \sigX ^2 \ = \ \mean{ (X - \muX )^2 }  \end{equation}
and thus
\begin{equation}\label{eq:var2}
    \sigX ^2 = \left\{%
\begin{array}{ll}
    \sum_x (x - \muX )^2 \pX(x), & \hbox{$X$ is discrete;} \\[2mm]
    \int_{\gmX} (\gmX(t)-\muX )^2\tfg dt, & \hbox{$X$ is continuous}.  \\
\end{array}%
\right.
\end{equation}

As usual, $\varr{X}$  and $\varp{X}$ denote respectively the real
\term \ and the phantom \term \ of $\var{X}$, where $\varr{X} =
\var{X_\pv}$ are taken with respect to the real \comp \ of the
probability  measure.  Thus, we always have $ \varr{X} \geq 0$,
since it is just a standard (real) variance. Moreover, we have the
following property:
\begin{proposition}\label{prop:var}
The variance $\var{X}$ is pseudo nonnegative, for any \prv \ $X$.
\end{proposition}
\begin{proof} In  view of Equation \eqref{eq:var2}, since
the square of a phantom number and probabilities are pseudo
nonnegative, then the proof is completed by Lemma
\ref{lem:positiveSquare} applied for the sum of the products.
\end{proof}

The \textbf{\topd \ variance}, i.e. a real-valued  function, is
defined
as $$\var{\Su{X}} = \mean{ (\Su{X} - \mean{\Su{X}})^2 }. $$ Since
this is a standard (real) variance, then we always have
$\var{\Su{X}} \geq 0$. The \textbf{conjugate variance} is defined
similarly as   $\var{\Du{X}} = \mean{ (\Du{X} - \mean{\Du{X}})^2 }
$, i.e. taken with respect to the conjugates of $X$ and the
probability measure.

Expanding the right-hand side of Equation \eqref{eq:var1}, we
obtain the following relation:
\begin{equation}\label{eq:var3} \var{X} = \mean{X^2}- \mean{X}^2,
\end{equation}
which is a useful formula for determining the variance. Plugging
Equation \eqref{eq:moment2} in to  this form, and simplifying, one
obtains   the realization property for variance:
\begin{equation}\label{eq:var4} \var{X} = \var{X_\pv} + \hi \(
\var{\Su{X}} - \var{X_\pv}\).
\end{equation}
We recall that the notation $\var{ X }$ is used for both the
phantom and the standard variance, where the meaning is understood
from the context.

\begin{remark} Both $\mean{X^n}$  and $\var{X}$ are phantom
functions whose real \comp s satisfy the familiar properties of
$n$'th moment and  variance.
\end{remark}

Assuming $g(X)$ is a phantom function of a \prv \ $X$, the
expected value of the \prv \ $g(X)$ is given by
\begin{equation*}\label{eq:meanFun}
    \mean{g(X)} \ =  \ \left\{%
\begin{array}{ll}
    \sum_x g(x) \pX(x), & \hbox{$X$  discrete;} \\[2mm]
    \int_{\gmX} g(\gmX(t)) \tfg dt, & \hbox{$X$ continuous.} \\
\end{array}%
\right.
\end{equation*}
It is straightforward to verify that when $g$ is a  linear phantom
function, say $g(X) = \al X + \bt$, with $\al, \bt \in \hf$, then
 \begin{equation}\label{eq:linearity}
 \mean{g(X)} = \al  \mean{X}+ \bt \qquad \text{and} \qquad \var{g(X)} = \al^2  \var{X},
\end{equation}
the latter formula is obtained by \eqref{eq:var2}.

\begin{proposition}
     $\Du{\var{X}} = \var{\Du{X}}$, where $\var{\Du{X}^n}$
     is computed with respect to the conjugates of $X$ and the
probability measure.

\end{proposition}
\begin{proof}  Use Equation \eqref{eq:var3} and  Properties \ref{prp:dusu}  to write
    $\Du{\var{X}} =  \Du{\mean{X^2}}- \Du{\mean{X}}^2$,
which by Proposition \ref{prop:conMean} is $\mean{\Du{X^2}}-\mean{
\Du{X}}^2$ and again by Properties \ref{prp:dusu}
$\mean{\Du{X}^2}-\mean{ \Du{X}}^2$, that is    $ \var{\Du{X}}$.
\end{proof}

Having this  property of Proposition \ref{prop:var}, we attained
the following additional phantom analog:
\begin{definition}
 The \textbf{standard phantom deviation} $\sigX $ of a \prv \ $X$,
is defined to be the maximal  nonnegative phantom square root of
$\var{X}$, cf. Equation \eqref{eq:root}, i.e. the  root with the
nonnegative real \term \ and the maximal nonnegative phantom
\term.
\end{definition}

Since each phantom number has a nonnegative square root, cf.
Equation \eqref{eq:root}, and $\var{X}$ is pseudo nonnegative,
then the standard phantom deviation is well defined for any \prv .
Using Equation \eqref{eq:root} it is easy to see that the
standard phantom deviation also admits the realization property,
i.e.
$$ \sigX  = \sig_{_{X_\pv}} + \hi ( \sig_{_{\Su{X}}} -  \sig_{_{X_\pv}}).$$

\subsection{Special  examples} The following examples show how the

classical mean and variance carry naturally on to the phantom
framework. When a \prv \ is discrete, we can retain  the exacrt
classical setting, while the continuous cases require a
modification of definitions which involves the parametrization of
$X$, i.e. that of $\gmX$. Yet, the standard (real) distributions
are received as the private cases for the phantom ones.

\begin{example}[\bfem{Mean and variance of the
Bernoulli}]\label{exmp:BernoulliMean}
Consider the experiment of tossing a biased coin, which comes up a
head with phantom probability $\gp$ and a tail with probability $1
- \gp$, and the Bernoulli \prv \ $X$ with \pmf
$$ \pX(k) \ = \ \left\{%
\begin{array}{ll}
    \gp, & \hbox{$k = 1$;} \\[1mm]
    1 -\gp , & \hbox{$k = 0$.} \\
\end{array}%
\right.    $$
Then $\mean{X} = 1 \gp + 0 (1 -\gp) = \gp $,  $ \mean{X^2} =  1^2
\gp + 0^2 (1 -\gp) = \gp$ and thus $\var{X} = \mean{X^2} -
\mean{X}^2 = \gp - \gp^2$.
\end{example}
\begin{example}
[\bfem{The mean of the Poisson}] The mean of the Poisson \pmf \
with pseudo positive parameter $\lm \in \hf$
$$\pX(k) \ = \ e^{-\lm} \frac{\lm^k} {k!}, \qquad k = 0, 1, 2,
\cdots ,$$
can be calculated as follows:
$$\mean{X}  \ = \ \sum_{k=0}^{\infty} k e^{-\lm} \frac{\lm^k} {k!}
\ \overset{*}{=} \ \sum_{k=1}^{\infty} k e^{-\lm} \frac{\lm^k}
{k!} \ = \ \lm \sum_{k=0}^{\infty} k e^{-\lm} \frac{\lm^{k-1}}
{k!} \ = \ \lm.
$$
($*$ the component indexed $k=0$ is zero.) A similar calculation
shows that the phantom variance of a Poisson random variable is
also $\lm$.
\end{example}

\begin{example}[\bfem{The phantom exponential
\prv}]\label{exmp:exoMean} Let $x(t) = a(t) + \hi b(t) $. We write
$a := a(t)$, $b:= b(t)$, and $x: = x(t)$, for short. The notation
$x'$ stands for the derivative of $x$ with respect to $t$, and
thus $x' = a' + \hi b'$.

 A
phantom exponential \prv \ has a \pdf \ with the form
$$ \fX (x) = \left\{%
\begin{array}{ll}
    \frac{\lm}{x'} e^{-\lm x}, & x \hbox{ is pseudo positive;} \\[2mm]
    0, & \hbox{otherwise ,} \\
\end{array}%
\right.     $$ where $\lm\in \hf$ is a pseudo positive parameter.
In particular  $\lm \neq 0$   and is not a zero divisor in $\hf$.

One observes that when $X$ is a real random variable, i.e $x(t) :
= t$, then $x' = 1$ and the phantom exponential $\fX$ collapses to
the known exponential random variable.

Using Equation \eqref{eq:prud} and Equation \eqref{eq:dirfrac},
respectively, we have
$$x \lm = a\lm_\pv +
  \hi(\su{x}\su{\lm}  -  a\lm_\pv), \qquad \text{and} \qquad
    \frac{\lm}{x'} \  =  \  \frac{\lm_\pv}{a'} + \hi\(
\frac{\su{\lm}}{\Su{x'}} - \frac{\lm_\pv}{a'} \).
$$ Then, Equation
  \eqref{eq:exp} yields
  $$
\begin{array}{lllll}
 \fX & = &   \frac{\lm}{x'} e^{-\lm x} &  =  &
 \( \frac{ \lm_\pv}{a'} + \hi\( \frac{ \su{\lm}}{\Su{x'}} -
\frac{\lm_\pv}{a'} \) \) e^{-a\lm_\pv -
  \hi(\su{x}\su{\lm}  -  a\lm_\pv)}  \\[1mm]
  &&  &  =  & \( \frac{ \lm_\pv}{a'} + \hi\( \frac{ \su{\lm}}{\Su{x'}} -
\frac{\lm_\pv}{a'} \) \)\( e^{-a\lm_\pv} + \hi(
  e^{-\su{x}\su{\lm}}  -  e^{-a\lm_\pv})\)  \\[1mm]
 && & = & \frac{\lm_\pv}{a'} e^{-a \lm_\pv}  + \hi \( \frac{ \su{\lm}}{\Su{x'}} e^{-\su{x} \su{\lm}}
   - \frac{\lm_\pv}{a'} e^{-a \lm_\pv} \).
\end{array}
$$
Now, for $\tfgt =  \fX \gmX'$ we get
$$
\begin{array}{lllll}
\tfgt & = &   \frac{\lm}{x'} e^{-\lm x} x' &   = &
 {\lm_\pv} e^{-a \lm_\pv}  + \hi \( { \su{\lm}} e^{-\su{x} \su{\lm}}
   - {\lm_\pv} e^{-a \lm_\pv} \),
\end{array}
$$
which shows that that $\tfgt$ satisfies the normalization
property. This because each component  is by itself a real
exponential, and thus the real \comp \ is $1$ and the phantom
\comp \ is summed up to $0$.

A similar computation as before shows that
$$ x\tfgt = {a\lm_\pv} e^{-a \lm_\pv}  + \hi \( {\su{x} \su{\lm}} e^{-\su{x} \su{\lm}}
   - {a\lm_\pv} e^{-a \lm_\pv} \) $$
 Recalling that for a real exponential random
variable, $\mean{X} = 1/\lm$ and $\var{X} = 1/\lm^2$, cf.
\cite{pCHU01a},  taking the integral
$$ E(X) \ =  \ \int_{0}^{\infty}
x(t)\tfg dt \ = \ \frac{1}{\lm_\pv} + \hi \(\frac{1}{\su{\lm}} -
\frac{1}{\lm_\pv} \)$$
we get the mean of exponential \prv \ in terms of the phantom
parameter $\lm$ as
$$ E(X)  \ = \  \frac{1}{\lm_\pv} + \hi \(\frac{1}{\su{\lm}} - \frac{1}{\lm_\pv} \)
\ =  \ \frac{1}{\lm_\pv} + \hi \(\frac{-\lm_\ph}{\lm_\pv(\lm_\pv +
\lm_\ph)}\) \ = \ \frac{1}{\lm},
$$
cf. Equation \eqref{eq:inv}.

In fact we could also have obtained this relation  in a shorter
way by using Equation \eqref{eq:moment2}, but, for the matter of
validation, we have presented the detailed computation.
\end{example}


\subsection{Normal random variables}
A continuous \prv \  $X$ is said to be \textbf{phantom normal}, or
\textbf{phantom Gaussian}, if it has a \pdf \  of the form
\begin{equation*}\label{eq:normal}
    \fX (x) = \frac{1}{\sig x'  \sqrt{ 2 \pi}}  e^{- (x - \mu)^2/2\sig
    ^2},
\end{equation*}
where $\mu$ and $\sig$ are two phantom scalar  parameters
characterizing the $\pdf$, with $\sig$ assumed pseudo positive.
For simplicity, we also assume $x'$ differentiable, and write $x'$
for the derivative of $x := x(t)= a(t) + \hi b(t)$ with respect to
$t$.

Note that in comparison to the classical case the normal \prv \
includes the extra argument $x'$ in the denominator. Yet, as we
had for the exponential \prv , when $X$ is assumed to take only
real values, the phantom normal density function collapses to the
classical normal density function.

\begin{proposition}\label{prop:Normal1}
 $\fX (x)$ satisfies the  the normalization property
\begin{equation}\label{eq:normal2}
 \frac{1}{\sig \sqrt{ 2 \pi}} \int_{\gmX} \frac{1}{x'} e^{- (x - \mu)^2/2\sig
    ^2}  dx \ = \
     1,\end{equation}
where $\gmX$ is parameterized by  $t \in \Real $, assumed
differentiable.
\end{proposition}
\begin{proof}
Let $w = x - \mu$, and therefore  $w' = x'$. Then, $(x -
\mu)^2/\sig ^2 = w^2/\sig^2$, and using simple computation one can
verify that
\begin{equation*}
 w^2/\sig
    ^2 \ = \  w_\pv ^2/ \sig_\pv^2 + \hi(\su{w} ^2/ \su{\sig}^2 - w_\pv ^2/ \sig_\pv^2).
    \end{equation*}
Plugging  this into $e^{-(x -\mu)^2/2\sig^2} $ and using Equation
\eqref{eq:exp}, we  have
$$e^{-(x - \mu)^2/2\sig^2}  \ =  \ e^{-w^2/2\sig^2}  \ =  \ e^{-w_\pv ^2/2 \sig_\pv^2} + \hi\(
e^{-\su{w} ^2/ 2\su{\sig}^2} -e^{-w_\pv ^2/2 \sig_\pv^2} \). $$
Thus,
$$\frac{1}{x'}e^{-w^2/2\sig^2}  \ =  \
\frac{1}{a'} e^{-w_\pv ^2/2 \sig_\pv^2} + \hi\(
\frac{1}{\Su{x'}}e^{-\su{w} ^2/ 2\su{\sig}^2} -
\frac{1}{a'}e^{-w_\pv ^2/2 \sig_\pv^2} \). $$
But then, $\frac{1}{x'}e^{-w^2/2\sig^2} \gmX'$ is just $ \
e^{-w_\pv ^2/2 \sig_\pv^2} + \hi\( e^{-\su{w} ^2/ 2\su{\sig}^2}
-e^{-w_\pv ^2/2 \sig_\pv^2} \)$.

Recalling that $\frac{1}{\sig} = \frac{1}{\sig_\pv} - \hi
\frac{\sig_\ph}{\sig_\pv \su{\sig}} = \frac{1}{\sig_\pv} + \hi(
\frac{1}{ \su{\sig}}- \frac{1}{\sig_\pv})$ and integrating this in
\eqref{eq:normal2}, given in parametric form, we get
\begin{equation}\label{eq:normDomp}
 \frac{1}{\sig_\pv \sqrt{ 2 \pi}} \int  e^{-w_\pv ^2/2 \sig_\pv^2} dt +
\hi\(   \frac{1}{\su{\sig}\sqrt{ 2 \pi}} \int e^{-\su{w} ^2/
2\su{\sig}^2} dt -  \frac{1}{\sig_\pv\sqrt{ 2 \pi}} \int e^{-w_\pv
^2/2 \sig_\pv^2} dt \) \  = \ 1 + \hi(1-1),
\end{equation} since the phantom \comp \ is
the sum of two standard normal distributions, each equal to $1$.
\end{proof}
Equation \eqref{eq:normDomp} shows that the phantom normal \pdf \
admits the realization property.

%
%

\begin{proposition}\label{prop:normal2}
    The mean and the variance of a normal \prv \ $X$ with phantom
    parameters $\mu$ and $\sig$ are
$$\mean{X} = \mu \qquad \text{and} \qquad \var{X}= \sig^2.$$
\end{proposition}
\begin{proof} Consider the realization property of  Equation \eqref{eq:normDomp} combined
respectively, with Equation \eqref{eq:moment2} and Equation
\eqref{eq:var4}.
\end{proof}

\begin{theorem} Normality is preserved under linear transformations.  If $X$ is
a normal \prv \  with mean $\mu$ and variance $\sig^2$, and if
$\al, \bt \in \hf$ are phantom scalars, then the \prv \ $Y =  \al
X + \bt$ is also normal, with mean and variance $$\mean{Y} = \al
\mu + \bt, \qquad \var{Y} = \al^2 \sig^2.$$
\end{theorem}
\begin{proof} Immediate by Proposition \ref{prop:normal2} and Equation
\eqref{eq:linearity}.
\end{proof}

A normal random variable $Y$ with zero mean and unit variance is
said to be a \textbf{standard phantom normal}. Its \cpdf , denoted
as $\Phi$, is given by
\begin{equation}\label{def:standardNormal}
\Phi(z) \ = \ \hpf(Y  \oleq z) \ = \ \frac{1}{ \sqrt{ 2 \pi}}
\int_S \frac{1}{y'} \, e^{- y^2/2} dy,
\end{equation}
where $S = \{ y \in \gmY  : \ y \oleq z\}$, assumed continuous and
differentiable. Clearly, this integral can also be written in the
parametric form as given in Equation \eqref{eq:pdfz2}.

Let $X$ be a normal \prv \ with mean $\muX$  and variance
$\sigX^2$. We ``standardize" $X$ by defining a new random variable
$Y$ given by
$$Y = \frac{X- \muX}{\sigX}.$$
 Since $Y$ is a linear transformation
of $X$, it is normal. Furthermore,
$$\mean{Y} \ = \  \frac{\mean{X} - \muX}{\sigX} \ = \ 0, \qquad
  \var{Y} \ = \ \frac{ \var{X}}{ \sigX } \  = \  1.$$
   Thus, $Y$ is a
standard normal \prv . This fact allows us to calculate the
probability of any event defined in terms of $X$: we redefine the
event in terms of $Y$, and then use the standard normal \prv .

The (classical) normal random variable plays an important role in
a broad range of probabilistic models. The main reason is that,
generally speaking, it models well the additive effect of many
independent factors, in a variety of engineering, physical, and
statistical contexts. As we have shown the normal \prv \ preserves
this property and generalizes the classical one in a natural way.

Mathematically, the key fact is that the sum of a large number of
independent and identically distributed (not necessarily normal)
phantom random variables has an approximately normal \cpdf,
regardless of the \cpdf \ of the individual random variables. This
property is captured in the celebrated central limit theorem,
extended to the phantom framework, which will be discussed in
Section~\ref{sec:limitThm}.
\section{Multiple random variable}
Consider a random experience having the sample space $\Om$.  A
\textbf{multiple phantom random variable}, written \textbf{\mprv}
for short,
 is a multiple-phantom-valued function
\begin{equation*}
    (X_1, \dots, X_n) : \Om \ \To \ \hf^{(n)},
\end{equation*}
given by
\begin{equation}\label{eq:mpvr}
    (X_1, \dots, X_n) : \om \ \longmapsto \ (X_1(\om), \dots, X_n(\om)),
\end{equation}
with each $X_i$ a \prv \ on $\Om$ as in Equation \eqref{eq:pRrad}.
The phantom range of the \mprv \ $(X_1, \dots, X_n)$  is denoted
by $\range{X_1, \dots,X_n}$, and defined by $$\range{X_1,
\dots,X_n} = \{ (x_1, \dots, x_n) \  :  \  \om \in \Om, \  x_1 =
X_1(\om),  \ \dots \ , \  x_n = X_n(\om) \}  \ .
$$

If the $X_i$'s are each, by themselves, discrete $\prv$'s, then $
(X_1, \dots, X_n)$ is called a discrete $\mprv$. Similarly, if the
$X_i$'s are each, by themselves, continuous \prv's, then $(X_1,
\dots, X_n)$ is called a continuous \mprv. When $n=2$ we write
$(X, Y)$ for $(X_1,X_2)$ and call it a
  \textbf{bivariate phantom random variable}, written
  \textbf{\bprv} for short.
In the remainder of this section, to make the exposition clearer,
we present the case of \bprv;  the extension to \mprv \ is
straightforward.

Consider two discrete \prv's $X$ and $Y$ associated with the same
experiment. The \textbf{joint phantom mass function} of $X$ and
$Y$ is defined by
$$\pXY (x, y) \ =  \ \hpf(X = x, Y = y)$$
  for all pairs
of phantom numerical values $(x, y)$ that $X$ and $Y$ can take;
otherwise it equals zero.  (Here and elsewhere, we will use the
abbreviated notation $\hpf(X = x, Y = y)$ instead of the more
precise notation $\hpf(\{X = x\} \cap \{ Y = y\})$.)

The joint \pmf \ determines the probability of any event that can
be specified in terms of the \prv's $X$ and $Y$. For example, if
$A$ is the set of all pairs $(x, y)$ that have a certain property,
then
$$\bhpf((x,y) \in A) \ = \  \sum_{(x,y) \in A} \pXY (x, y). $$
In fact, as in  classical theory, we can calculate the \pmf's of
$X$ and $Y$ by using the formulas
$$\pX(x) = \sum_{y} \pXY (x, y), \qquad \pY(y) = \sum_{x} \pXY (x, y), $$
where $x$ and $y$ range respectively over all the phantom  values
of $X$ and $Y$.

\begin{definition}
Two discrete \prv 's  $X$ and $Y$ are  said to be
\textbf{independent} if
$$\pXY (x, y) \ = \ \pX(x) \,  \pY (y), \qquad \text{ for all }  x,
y.$$
\end{definition}

We say that two continuous \prv 's  associated with a common
experiment are \textbf{jointly continuous}, and can be described
in terms of a joint \pdf \ $\fXY$, if $\fXY$ is a continuous
function whose real \comp \ is nonnegative and that satisfies
\begin{equation}\label{eq:JC1}
\bhpf((X,Y) \in B) \ = \ \underset{B \ \ }{\int \int} \fXY(x, y)
dx fy
\end{equation}
for every subset $B$ of
$$\gmX \times \gmY  \ = \ \{ (x,y ) \ : \ x \in \gmX, \ y \in \gmY
\},
$$
where $\gmX$ and $\gmY$ are as in Equation \eqref{eq:gmX}.

Accordingly, given $z_1,z_2, w_1,w_2 \in \hf$, we define
$$ \SX \ = \ \{ x \in \gmX \ : \ z_1 \oleq x \oleq z_2\},
 \qquad  \SY \ = \ \{ y \in \gmX \ : \ w_1 \oleq y\oleq w_2\}, $$
and have
$$\hpf(z_1 \oleq X \oleq z_2, \ w_1 \oleq Y \oleq w_2 ) \ = \
 \int_{\SY} \int_{\SX}\fXY(x, y) dx fy .
$$
(Note that when $\SX$ or $\SY$ are not connected, the integral
decomposes into a sum of  integrals, assumed finitely many.)

Furthermore, by letting $B$ in Equation \eqref{eq:JC1} be the
entire set $\gmX \times \gmY$, we obtain the normalization
property
$$\bhpf((X,Y) \in B) \ = \ \int_{\gmY} \int_{\gmX} f_{X,Y}(x, y) dx fy \ = \ 1.$$
As before, when $\SX$ and $\SY$ are subpaths of $\gmX$ and $\gmY$,
respectively, we can use the orders $\tleq$ and $\sleq$  induced
on $\gmX$ and $\gmY$ by their parametrization, respectively,
together with the integral form  \eqref{eq:integral2}, and write:
$$\hpf(x_1 \tleq X \tleq x_2, \ y_1 \sleq Y \sleq y_2 ) \ = \
 \int_{\itt{y_1}}^{\itt{y_2}} \int_{\itt{x_1}}^{\itt{x_2}}
\fXY (\gmX(t),\gmY(s))\gmX'(t)\gmY'(s)dtds .
$$
for $x_1, x_2 \in \gmX$,  $y_1, y_2 \in \gmY$.

The \textbf{marginal} $\pdf$ 's  $\fX $  and $f_Y$ of  $X$  and
$Y$, respectively, are given by:
$$ \fX (x) =  \int_{\gmY}  \fXY (x, y) dy, \quad \text{and}  \quad
f_Y(y) = \int_{\gmX} \fXY (x, y) dx.$$

In full analogy with the discrete case, we say that two continuous
\prv 's $X$ and $Y$ are \textbf{independent} if their joint \pdf \
is the product of their  marginal \pdf's:
$$\fXY (x, y) \ = \  \fX
(x)f_Y (y), \qquad \text{for all } x, y.$$

By simple computation one can verify that:

\begin{properties}
 If $X$ and $Y$ are
independent \prv 's  then:
\begin{enumerate}
    \item The \prv 's $g(X)$ and $h(Y)$ are independent, for any functions $g$ and $h$, \pSkip
    \item $\mean{z_1 X + z_2 Y + z_3} = z_1 \mean{X} + z_2 \mean{Y} + z_3
    $, \
    for $z_i \in \hf$, \pSkip
    \item  $\mean{XY} =
\mean{X}\mean{Y}$,  and more generally $\mean{g(X) h(Y)} =
\mean{g(X)}\mean{h(Y)}$, \pSkip
    \item $ \var{X+Y} = \var{X} +
\var{Y}. $
\end{enumerate}
\end{properties}

These properties can be verified easily by direct computation or
by using the relaxation property admitted by $\mean{ \ }$ and
$\var{ \ }$ and validation of these property for the standard
(real) cases.

\subsection{Covariance and correlation}

The \textbf{covariance},  denoted by $\cov{X, Y }$, of two \prv 's
$X$ and $Y$ is defined as
\begin{equation}\label{eq:cov}
\cov{X, Y } \ = \ \mean{\ (X - \mean{X})( Y - \mean{Y }) \ }.
\end{equation}
The \prv 's  $X$ and $Y$ are said to be \textbf{uncorrelated} if
 $\cov{X, Y } = 0$.  When $\covr{X, Y } = 0$ we say that
$X$ and $Y$ are \textbf{real uncorrelated}, and if $\covp{X, Y }=
0$ we say that $X$ and $Y$ are \textbf{phantomly uncorrelated}.

We let  $$\cov{ \Su{X}, \Su{Y} } \ = \  \mean{ \ (\Su{X} -
\mean{\Su{X}})( \Su{Y} - \mean{\Su{Y} }) \ },$$ and call it the
(real) \textbf{\topd \ covariance} of $X$ and $Y$, where
expectations are computed with respect to the reductions of $X$
and the probability measure.

Let $U = X - \mean{X} $ and $V = Y - \mean{Y }$, then
$$\cov{X,Y} = \mean{U_\pv V_\pv + \hi (\Su{U} \Su{V} - U_\pv
V_\pv)},$$
 cf. Equation \eqref{eq:pwo2}, which  by Equation
\eqref{eq:moment2} is
$$\cov{X,Y} = \mean{U_\pv  V_\pv}  + \hi (\mean{\Su{U} \Su{V}} -
\mean{U_\pv  V_\pv}).$$ In other words,
\begin{equation}\label{eq:cov2}
\cov{X, Y } = \cov{X_\pv, Y_\pv}  + \hi (\cov{\Su{X}, \Su{Y}} -
\cov{X_\pv, Y_\pv}),
\end{equation}
that is the phantom covariance admits the realization property.

Roughly speaking,  positive or negative parts (cf. Definition
\ref{def:pseudoPositive}) of covariance indicate that the values
of $X - \mean{X}$ and $Y - \mean{Y}$ obtained in a single
experiment "tend" to have the same or the opposite sign,
respectively. Thus, the signs of the real and the phantom \term \
of the covariance provide an important qualitative indicator of
the relation between the real \comp s and the phantom \comp s of
$X$ and $Y$. If $X$ and $Y$ are independent, then
\begin{equation*}\label{eq:cov1}
\begin{array}{lllll}
  \cov{X, Y} & = & \mean{ \ (X - \mean{X})( Y - \mean{Y }) \ } &&\\ [1mm]
   & = & \mean{X -
\mean{X}} \ \mean{ Y - \mean{Y }} & = & 0. \\
\end{array}
\end{equation*}
Therefore, if $X$ and $Y$ are independent, they are also
uncorrelated. However,  as in  classical theory, the reverse is
not true.

The \textbf{correlation coefficient} $\rho(X,Y)$ of two \prv's $X$
and $Y$, whose variances are nonzero divisors, is defined as
\begin{equation}\label{eq:corr}
\rho(X,Y) \ = \ \frac{\cov{X, Y}}{\sqrt{\var{X}\var{Y }}} .
\end{equation}
This maybe viewed as a normalized version of the phantom
covariance $\cov{X, Y}$, and as the computation below shows, the
real \term \ of $\rho(X,Y)$ ranges from $-1$ to $1$.

Using Equation \eqref{eq:var4} and Equation \eqref{eq:cov2},
together with the realization properties of the square root
\eqref{eq:root}, Equation \eqref{eq:corr} receives the familiar
form:
\begin{equation*}
\begin{array}{lll}
  \rho(X,Y) &  =  & \frac{\cov{X_\pv Y_\pv}  + \hi (\cov{\Su{X} \Su{Y}} -
\cov{X_\pv Y_\pv})}{\sqrt{(\var{X_\pv} + \hi (\var{\Su{X}}
-\var{X_\pv}))}\sqrt{(\var{Y_\pv} + \hi( \var{\Su{Y}} -
\var{Y_\pv}))}}
\\[2mm]
  &  =  & \frac{\cov{X_\pv Y_\pv}  + \hi (\cov{\Su{X} \Su{Y}} -
\cov{X_\pv Y_\pv})}{(\sqrt{\var{X_\pv}} + \hi (\sqrt{\var{\Su{X}}}
-\sqrt{\var{X_\pv}}))(\sqrt{\var{Y_\pv}} + \hi(
\sqrt{\var{\Su{Y}}} - \sqrt{ \var{Y_\pv}}))}
\\[2mm]
  & =  &   \frac{\cov{X_\pv Y_\pv}  + \hi (\cov{\Su{X} \Su{Y}} -
\cov{X_\pv Y_\pv})}{\sqrt{\var{X_\pv}\var{Y_\pv}} + \hi(
\sqrt{\var{\Su{X}} \var{\Su{Y}} } - \sqrt{\var{X_\pv} \var{Y_\pv}})} \\
[2mm]
& = &  \frac{\cov{X_\pv Y_\pv}}{\sqrt{\var{X_\pv}\var{Y_\pv}}}  +
  \hi \frac{\cov{\Su{X} \Su{Y}} \sqrt{\var{X_\pv}\var{Y_\pv}} - \cov{X_\pv Y_\pv}  \sqrt{\var{\Su{X}}
\var{\Su{Y}}} }{\sqrt{\var{X_\pv}\var{Y_\pv}} \sqrt{\var{\Su{X}}
\var{\Su{Y}}}}.
\end{array}
\end{equation*}
%
%
Therefore,
\begin{equation}\label{eq:cor2}
  \rho(X,Y) = \ \rho(X_\pv, Y_\pv)
   + \hi \(  \rho(\Su{X}, \Su{Y})  -\rho(X_\pv, Y_\pv)  \), \hskip
   21mm
\end{equation}
which is the realization property for phantom covariance.

 Let $\xptr$ be the set
\begin{equation}\label{eq:teta}
 \xptr= \{ z \in \hf \   | \ a \in [-1,1], \  -(1+a) \leq  b \leq 1-a \},
\end{equation}
i.e. it is  the pointwise product  $2 \ptr - 1$. We write
$\xptr_{(+,+)}$ for the subset of $\cptr$ consisting of all
phantom points whose real and phantom \term s are positive;
$\xptr_{(+,-)}$, $\xptr_{(-,+)}$, and $\xptr_{(-,-)}$ are defined
respectively according to the positivity signs of the real \term \
and the phantom \term \ of their points.

\begin{proposition} Given any two \prv's $X$ and
$Y$, then $ \rho(X, Y) \in \xptr$ as defined in Equation
\eqref{eq:teta}.
\end{proposition}
\begin{proof}
Using  classical theory, both  $\rho(\Su{X}, \Su{Y})$ and $
\rho(X_\pv, Y_\pv)$ range  from $-1$ to $1$; the proof is
completed by Equation \eqref{eq:cor2}.
\end{proof}

If $\rho \in \xptr_{(+,+)} $ (or $\rho \in \xptr_{(-,-)}$), then
the real and the phantom values of $x - \mean{X}$ and $y -
\mean{Y}$ ``tend'' to have the same (or opposite, respectively)
sign, and the size of $|\rho|$ provides a normalized measure of
the extent to which this is true. In fact, always assuming that
$X$ and $Y$ have positive variances, it can be shown that $\rho =
1$ (or $\rho = -1$) if and only if there exists a constant
positive phantom number $\al$, or negative, respectively, such
that
$$y - \mean{Y} = \al (x - \mean{X}), \quad \text{for all possible
numerical values } (x, y).$$
When $\rho \in \xptr_{(+,-)} $, or $\rho \in \xptr_{(-,+)} $, then
the real \term s of $x - \mean{X}$ and $y - \mean{Y}$ ``tend'' to
have the same (or opposite, respectively) sign opposite to that of
their phantom  \term s.

\section{Moment generating functions}


The \textbf{moment generating function}, written \textbf{\mgf} \
for short, of the distribution function of a \prv \  $X$ (also
referred to as the \textbf{transform} of $X$) is a phantom
function $M_X(\gs)$ of a free phantom parameter $\gs \in \hf$,
defined by
$$\MX(\gs )  \ =  \ \mean{e^{\gs X}}.$$
%
In more detail, the corresponding transform of $X$ is given by:
\begin{equation}\label{eq:mgf}
    \MX(\gs) \ =  \ \left\{%
\begin{array}{ll}
    \sum_x e^{\gs x} \pX(x), & \hbox{$X$  discrete;} \\[2mm]
    \int_{\gmX} e^{\gs x} \fX (x) dx, & \hbox{$X$ continuous.} \\
\end{array}%
\right.
\end{equation}

Let $\gs = \gs_\pv +\hi \gs_\ph $ and use Equation \eqref{eq:exp}
to write
\begin{equation*}\label{eq:mgf1}
\begin{array}{lll }
   \mean{e^{\gs X}} & = &  \mean{e^{\gs_\pv X_\pv
+ \hi(\su{\gs}\Su{X}-  \gs_\pv X_\pv)}}\\[1mm]
   & = & \mean{e^{\gs_\pv X_\pv} + \hi(e^{
\su{\gs}\Su{X}}- e^{\gs_\pv X_\pv})}. \\
\end{array}
\end{equation*}
Then, by Equation \eqref{eq:moment2}, one has the realization
property for \mgf
\begin{equation*}\label{eq:mgf2} \mean{e^{\gs X}} \ = \  \mean{e^{\gs_\pv X_\pv}} + \hi \(\mean{e^{
\su{\gs}\Su{X}}}- \mean{e^{\gs_\pv X_\pv}}\),
\end{equation*}
and thus
\begin{equation}\label{eq:mgf3} \MX(\gs) \ = \ \MX(\gs_\pv) +
\hi(\MX(\su{\gs})-\MX(\gs_\pv)),
\end{equation}
where $\MX(\gs_\pv))$ and $\MX(\su{\gs})$ are standard (real)
moment generating functions.

\begin{theorem}[\bfem{Inversion property}] The \mgf \ $\MX(\gs)$ completely determines
the probability law of the random variable $X$. In particular, if
$\MX(\gs) = \MY (\gs)$ for all $\gs$, then the random variables
$X$ and $Y$ have the same probability law.
\end{theorem}
This property is a rather profound mathematical fact that is used
frequently in  classical probability theory. In  light of Equation
\eqref{eq:mgf3}, i.e. the realization property of $\MX(\gs)$, this
phantom property is derived directly from the known result for the
standard (real) \mgf , applied to each comportment,  in classical
probability theory \cite{pCHU01a}.

Transform methods are particularly convenient when dealing with a
sum of \prv 's, since it covers addition of independent \prv \ to
multiplication of transforms, as we now show. Let $X$ and $Y$ be
independent \prv 's, and let $W = X + Y$. The transform associated
with W is, by definition,
$$\MW(\gs) \ = \ \mean{e^{\gs W}} \ = \ \mean{e^{\gs (X+Y )}} \ =
\ \mean{e^{\gs X}e^{ \gs Y}};$$
the last equality is due to Equation \eqref{eq:exp}.

Consider a fixed value of the parameter $\gs \in \hf$. Since $X$
and $Y$ are independent, $e^{\gs X}$ and $e^{\gs Y}$ are also
independent \prv 's. Hence, the expectation of their product is
the product of the expectations, and thus
$$
\MW(\gs) \ = \ \mean{e^{\gs X}}\mean{e^{\gs Y }} \ = \ \MX(\gs) \,
\MY (\gs).$$
 By the
same argument, if $X_1, \dots , X_n$ is a collection of
independent \prv 's, and $W = X_1 + \cdots + X_n$, then
$$\MW(\gs) \ = \ \MXi{1}(\gs) \ \cdots \ \MXi{n}(\gs).$$

\subsection{Examples of moment generating functions}

\begin{example}[\bfem{The transform of a linear function of a random variable}]
Let $\MX(\gs)$ be the transform associated with a \prv \  $X$.
Consider a new \prv \  $Y = uX + v$ for some  $u,v \in \hf$. We
then have
$$ \MY (\gs) \ = \ \mean{e^{\gs(uX+v)}} \ = \ e^{v \gs }\mean{e^{u \gs X}} \
= \ e^{v\gs }\MX(u \gs ).$$
\end{example}

\begin{example}[\bfem{The transform of the binomial}]
Let $X_1, \dots , X_n$ be independent Bernoulli \prv's, cf.
Example \ref{exmp:BernoulliMean}, with a common parameter $\gp$,
assigned to probability.   Then,
$$\MXi{i} (\gs) = (1 - \gp)e^{0 \gs} + \gp e^{1 \gs} = 1-
\gp + \gp e^\gs, \qquad \text{for all} \ \ i.$$ The \prv \ $Y =
X_1 + \cdots + X_n$ is phantom binomial with parameters $n \in
\Net$ and $\gp \in \hf$. Its transform is given by
$$\MY (\gs) = (1 - \gp + \gp e^\gs) ^n .$$
\end{example}

\begin{example}[\bfem{The sum of independent Poisson random variables is Poisson}] Let
$X$ and $Y$ be independent Poisson \prv 's with means $\muX$ and
$\muY$, respectively, and let $W = X + Y$. Then,
$$\MX(\gs) =
e^{\muX(e^\gs-1)}, \qquad \MY (\gs) = e^{ \muY(e^\gs-1)},$$
  and $$
\MW(\gs) \ = \ \MX(\gs) \MY (\gs) \ = \ e^{\muX(e^\gs-1)} e^{
\muY(e^\gs-1)} \ = \ e^{(\muX + \muY)(e^\gs-1)}.$$
Thus, $W$ has the same transform as a Poisson \prv \  with mean $
\muX + \muY$. By the uniqueness property of transforms, $W$ is
Poisson with mean $\muX + \muY$.
\end{example}


\section{Limit theorems}\label{sec:limitThm}

\subsection{Some useful inequalities} Before getting to
probability inequalities, we need furhter results about the weak
order $\oleq$ on $\hf$, including its relations with the phantom
absolute value as defined in Equation \eqref{eq:abs}. We recall
that $\oleq$  assumed satisfying Properties
\ref{prp:orderProperties}.

\begin{remark}\label{rmk:orderAbs} The classical relation $ z^2 = \abs{z}^2$ does not always
hold phantomly; we might have
 $ z^2 \ole \abs{z}^2$ or  $ z^2 \oge \abs{z}^2$. Note that
 $\abs{z}^2$ is real while $z^2$ is phantom. It is easy to verify
 that $ z^2 = \abs{z}^2$ holds   iff $b= -2a$.

 Moreover, from a metric point of view, there are numbers that are ``close'' in the sense of the weak order
 $\oleq$, but very far in the sense of $\abs{ \ }$; for example
 assuming $\oleq$ is the lexicographic order, a small increasing of $\ep$ makes
  $z_1 = \ep + \hi \ep$ greater than $z_2 = \ep + \hi b$, but sill
  $\abs{z_2} > \abs{z_1}$.
\end{remark}

The mismatch between $\oleq$ and the $\abs{ \ }$, as addressed in
Remark \ref{rmk:orderAbs}, yields different versions for phantom
Markov inequalities, aiming to provide later a phantom version of
the Chebyshev inequality.

\begin{proposition}[\bfem{Markov phantom inequalities}]\label{prop:markov} Given  a \prv \ $X$
that takes only values $\ogeq 0 $. Then
\begin{enumerate} \eroman
    \item $\hpf(X \ogeq  z) \oleq  \frac{\mean{X}}{ z}$ , \qquad for any
pseudo positive $z \in
    \hf$,\pSkip

    \item $\hpf(\abs{X} \geq  \abs{z}) \oleq  \frac{\mean{\abs{X}}}{\abs{z}}$ , \qquad for any   $z \in
    \hf$,
\pSkip

      \item $\abs{\hpf(\abs{X} \geq  \abs{z})} \leq  \abs{\frac{\mean{\abs{X}}}{\abs{z}}}$ , \qquad for any   $z \in
      \hf$.

\end{enumerate}
\end{proposition}

\begin{proof} \eroman
\begin{enumerate}
    \item

Fix a pseudo positive $z \in \hf$ and consider the random variable
$Y_z$
defined by $$Y_z = \left\{%
\begin{array}{ll}
    0, & \hbox{X $ \oleq$ z;} \\
    z, & \hbox{X $ \oge$ z.} \\
\end{array}%
\right.     $$ It is seen that the relation $Y_z \oleq X$ always
holds and therefore, using Properties \ref{prp:orderProperties}
(ii) for sums and products, $\mean{Y_z} \oleq \mean{X}$. On the
other hand, $\mean{Y_z} = zP(Y_z \oeq  z) = zP(X \ogeq z)$, from
which we obtain $z\hpf (X \ogeq z)  \oleq \mean{X}.$ The proof is
then completed by \ref{prp:orderProperties} (ii) for division.
\pSkip
%

 \item Apply part (i) to $\abs{X}$ and
$\abs{z}$, since both are positive. \pSkip

    \item We need to prove that $\abs{\mean{\abs{Y_z}}} \leq
    \abs{\mean{\abs{X}}}$, or equivalently that $0\leq
    \abs{\mean{\abs{X}}}^2 - \abs{\mean{\abs{Y_z}}}^2$;  then the required  inequality is obtained by part (ii).
We prove the assertion for a discrete \prv ; the continues version
is received similarly.

     Let $y_x \in \{ 0, \abs{z}\}$ for the value of
    $Y_z$ apply to  $x \in X$;  accordingly $\pX(x) = \pY(y_x)$ for each $x \in
    X$. We write $x$ and $y_x$ for $\abs{x}$  and $\abs{y_x}$,
    respectively, assuming both are real nonnegatives.
    Then, denoting $\pXr$ and $\pXp$ the real and the phantom \comp \ of $\pX$, respectively,  \\
$ \begin{array}{lll}
     \abs{\mean{\abs{Y_z}}}^2 & = &  \meanr{\abs{Y_z}}^2 + \meanr{\abs{Y_z}}
\meanp{\abs{Y_z}} +  \meanp{\abs{Y_z}}^2 / 2 \\ [1mm]
     & = &  \sum_{x', x''} y_{x'} y_{x''} \pXr(x') \pXr(x'') +
     \sum_{x', x''} y_{x'} y_{x''} \pXr(x') \pXp(x'')  \\ [1mm]
& &    + \sum_{x', x''} y_{x'} y_{x''} \pXp(x') \pXp(x'')/2 \\
[1mm]& = & \sum_{x', x''} y_{x'} y_{x''}(\pXr(x') \pXr(x'') +
\pXr(x') \pXp(x'') + \pXp(x') \pXp(x'')/2),
   \end{array}
$
\\ and $\abs{\mean{\abs{X}}}^2$ is expressed in the same way.

Letting  $\gX(x',x'' ) = \pXr(x') \pXr(x'') + \pXr(x') \pXp(x'') +
\pXp(x') \pXp(x'')/2$, as it is derived from the absolute value,
one observes that $\gX(x',x'' ) \geq 0$.

Putting all together, and considering the difference, we have
$$ \abs{\mean{\abs{X}}}^2 - \abs{\mean{\abs{Y_z}}}^2 = \sum_{x', x''} (x' x'' -  y_{x'} y_{x''}) \gX(x',x'' ),$$
in which all components are $\geq 0$. Since $x' \geq y_{x'}$ and
$x'' \geq y_{x''}$ then $x' x'' -  y_{x'} y_{x''} \geq 0$, and
thus the sum is $\geq 0$ as desired.
\end{enumerate}
\end{proof}

We write $\muXabs$ and $\sigXabs^2$ for $\mean{\abs{X}}$ and
$\var{\abs{X}}$, respectively, then have the phantom analogously
to the Chebyshev inequality.

\begin{proposition}[\bfem{Chebyshev phantom inequality}] If $X$ is a random variable with mean $\muXabs$ and
variance $\sigXabs ^2$, then
$$\abs{\hpf\( \abs{\abs{X} - \muXabs} \geq \abs{z} \)} \ \leq
\ \abs{\frac{\sigXabs^2}{\abs{z}^2}} , \qquad \text{for all }  z
\neq 0.$$
\end{proposition}
\begin{proof}
Consider the nonnegative random variable $(\abs{X} - \muXabs )^2$
and apply the Markov inequality (iii) with $z = \abs{w}^2$ to
obtain
$$\abs{\hpf\( \abs{(\abs{X} - \muXabs )^2} \geq \abs{w}^2 \)} \
\leq \
 \abs{ \frac{\mean{ \abs{(\abs{X} - \muXabs )^2} }}{\abs{w}^2}}.$$
Since, $(\abs{X} - \muXabs )^2$ is a real nonnegative number,
$\abs{(\abs{X} - \muXabs )^2} = (\abs{X} - \muXabs )^2$, and thus
$$\abs{\hpf\( {(\abs{X} - \muXabs )^2} \geq \abs{w}^2 \)} \ \leq \
\abs{ \frac{\mean{ {(\abs{X} - \muXabs )^2} }}{\abs{w}^2}} \ = \
\abs{\frac{\sigXabs^2}{ \abs{w}^2}} .$$

 The derivation is
completed by observing that the event  ${(\abs{X} - \muXabs )^2}
\geq \abs{w}^2$ is identical to the event $\abs{\abs{X} - \muXabs
} \geq \abs{w}$ and
$$
\abs{\hpf\( \abs{\abs{X} - \muXabs } \geq \abs{w} \)} \ = \
\abs{\hpf\( {(\abs{X} - \muXabs )^2} \geq \abs{w}^2 \)} \ \leq \
\abs{\frac{\sigXabs^2}{ \abs{w}^2}}.
$$
\end{proof}

An alternative form of the Chebyshev inequality is obtained by
letting $\abs{w} = c \sigXabs.$, where $c$ is a real positive,
which yields
$$
\abs{\hpf\( \abs{\abs{X} - \muXabs } \geq c \sigXabs \)} \ = \
\abs{\frac{\sigXabs^2}{ (c \sigXabs)^2}} \ = \ \frac{1}{ c^2}.
$$
Thus, the probability that a random variable $\abs{X}$ takes a
value more than $c$ times the standard deviations away from the
mean $\muXabs$ is at most $1/c^2$.

The Chebyshev inequality is generally more powerful than the
Markov inequality (the bounds that it provides are more accurate),
because it also makes use of information on the variance of $X$.
Still, as usual, the mean and the variance of a random variable
are only a rough summary of the properties of its distribution,
and we cannot expect the bounds to be close approximations of the
exact probabilities.
\subsection{The weak law of large numbers}

Consider a sequence $X_1,X_2, \dots$ of independent identically
distributed \prv's, each with mean $\mu$ and variance $\sig^2$.
Let
$$ S_n = X_1 + \cdots + X_n$$ be the sum of the first $n$ of them.
As in  classical theory, phantom limit theorems are mostly
concerned with the properties of $S_n$ and related \prv's, as $n$
becomes very large. In fact, the realization property of phantoms
provides the phantom analogues to these theorems in a trivial way.

Because of the independence of $X_i$'s, we have
$$\var{S_n} = \var{X_1} + \cdots
+ \var{X_n} = n \sig^2.$$
Thus, the distribution of $S_n$ spreads out as $n$ increases, and
does not have a meaningful limit. The situation is different if we
consider the \textbf{sample mean}
$$
M_n = \frac{X_1 + \cdots + X_n}{ n} = \frac{S_n}{ n},$$
which can also be written as
\begin{equation}\label{eq:Mn}
    M_n = M_{n,\pv} + \hi (\Su{M_n} - M_{n,\pv}).
\end{equation}
 A quick
calculation, together with the independence,  shows that
$$\mean{M_n} = \mu, \qquad  \var{M_n} = \frac{\sig^2}{ n} .$$

We apply Chebyshev inequality and obtain
\begin{equation}\label{eq:Cheb2} \abs{\hpf\( \abs{\abs{M_n} - \muMnabs} \geq \ep \)} \ \leq
\ \abs{\frac{\sigMnabs^2}{n \ep^2}}
 \qquad \text{for any real } \ep >0.
\end{equation}
We observe that for any real fixed $\ep  > 0$, the right-hand side
of this inequality goes to zero as $n$ increases. This form gives
one way to approach phantom limit theorems. However, in the
sequel, we focus on the  way established by the realization
property. This means that we consider phantom probability for
abstract events, or random variables.

Next we consider the phantom weak law of large numbers,  stated
below. It turns out that this law remains true even if the $X_i$
have infinite variance, but a much more elaborate argument is
needed, which we omit. The only assumption needed is that
$\mean{X_i}$ is well-defined and finite.

\begin{theorem}[\bfem{The weak law of large numbers (WLLN)}]\label{thm:WLLN} Let $X_1,X_2, \dots$ be
independent identically distributed \prv's  with mean $\mu$. For
every real $\ep \ge 0$, we have
$$\hpf( \abs{M_n - \mu} \geq \ep )
\To \  0, \qquad \text{as } \  n \to \infty, $$
or equivalently
$$\hpf( \abs{M_n - \mu}< \ep )
\To \  1, \qquad \text{as } \  n \to \infty. $$
\end{theorem}
\begin{proof} Recall that  $\abs{M_n - \mu} \to 0 $ iff both, $\re{M_n -
\mu} \to 0 $ and $\phf{M_n - \mu} \to 0$, cf. Lemma
\ref{lem:converge}, and that $\hpf_\pv$ is a standard (real)
probability measure for any phantom probability measure $\hpf =
\hpf_\pv + \hi \hpf_\ph$. Then, since $\abs{M_n - \mu} < \ep$ is
an inequality of real random variables, by the known WLLN for real
probabilities, $\hpf_\pv( \abs{M_n - \mu} < \ep ) \to 1$ as $n \to
\infty$, which means $\hpf_\ph( \abs{M_n - \mu} < \ep ) \to 0$,
since $\hpf$ is a phantom probability measure.
\end{proof}

As in classical theory, the phantom WLLN states that for a large
$n$, the ``bulk" of the distribution of $M_n$ is concentrated near
$\mu$. That is, if we consider a neighborhood around $\mu$, which
here is 2-dimensional, then there is a high probability that $M_n$
will fall in that neighborhood; as $n \to \infty$, this
probability converges to $1$. Of course, if $\ep$ is very small,
we may have to wait longer (i.e., need a larger value of $n$)
before we can assert that $M_n$ is highly likely to fall in that
neighborhood.

\begin{corollary} Let $X_1,X_2, \dots$ be
independent identically distributed \prv's  with mean $\mu$. For
every real $\ep \ge 0$, we have
$\abs{\hpf( \abs{M_n - \mu}< \ep )} \to   1$, as $n \to \infty$.
\end{corollary}


\subsection{The central limit theorem}

We can interpret the WLLN as stating that "$M_n$ converges to
$\mu$." However, since $M_1,M_2, \dots$ is a sequence of phantom
random variables, not a sequence of phantom  numbers, the meaning
of convergence, in the phantom sense, has to be precise. A
particular definition is provided below. To facilitate the
comparison with the ordinary notion of convergence, we also
include the definition of the latter.

\begin{definition}

 Let $X_1, X_2, \dots$ be a sequence of
\prv's (not necessarily independent), and let $z$ be a phantom
number. We say that the sequence $X_n$ \textbf{converges to $\bf
z$ in probability}, if for every real $\ep > 0$, we have
$$\tLim{n}{\infty}\hpf\( \abs{X_n - z} \geq \ep \) = 0,$$
or equivalently, for every real $\dl > 0$ and for every real $\ep
> 0$, there exists some $n_0$ such that $$\abs{\hpf\( \abs{X_n - z} \geq
\ep \)} \ \leq \ \dl,$$ for all $n > n_0$.
\end{definition}

According to the weak law of large numbers, the distribution of
the sample mean $M_n$  is increasingly concentrated in the near
vicinity of the true mean $\mu$.  In particular, its variance
tends to zero. On the other hand, the variance of the sum $S_n =
X_1 + \cdots + X_n = n M_n$ is unbounded, and the distribution of
$S_n$ cannot be said to converge to anything meaningful.

An intermediate view is obtained by considering the deviation $S_n
- n\mu$ of $S_n$ from its mean $n\mu$, and scaling it by a (real)
factor proportional to $1/\sqrt{n}$. What is special about this
particular scaling is that it keeps the variance, even though it
is phantom, at a constant level. The central limit theorem asserts
that the distribution of this scaled phantom random variable
approaches a normal phantom distribution.

More specifically, let $X_1,X_2, \dots$  be a sequence of
independent identically distributed \prv 's with mean $\mu$ and
variance $\sig^2$. We define
\begin{equation}\label{eq:z1}
W_n \ = \ \frac{S_n - n \mu}{\sig \sqrt{n}} \ = \  \frac{X_1 + \
\cdots \ + X_n - n \mu}{\sig \sqrt{n}}.
\end{equation}
 An easy calculation yields:
$$\mean{W_n} \  = \ \frac{\mean{X_1 + \ \cdots \ + X_n}  - n \mu}{\sig \sqrt{n}} \
= \  0,$$ and
$$
\begin{array}{lllllll}
  \var{W_n} & = & \frac{\var{X_1 + \ \cdots \ + X_n}} {\sig^2 \sqrt{n}}
  && &&
  \\[1mm]
  & = & \frac{ \var{X_1} + \ \cdots \ + \var{X_n}}{\sig^2 \sqrt{n}} &= &
  \frac{\sig^2 \sqrt{n}}{\sig^2 \sqrt{n}} & = & 1  \\
\end{array}
$$

\begin{theorem}[\bfem{The phantom central limit theorem}] Let $X_1,X_2, \dots$ be a sequence of
independent identically distributed \prv 's with common mean $\mu$
and a finite variance $\sig^2$, and let $W_n$ be defined as in
Equation \eqref{eq:z1}. Then, the \cpdf \  of $W_n$ converges to
the standard normal $\cpdf$; that is, for a given phantom value $z
\in \hf$,
$$ \Phi(z) \ = \ \frac{1}{\sqrt{2\pi}} \int_{S} \frac{1}{w'} e^{-w^2/2} dw,$$
with  $S = \{w \in \gm_{_{W_n}} \ : \   w \oleq z   \}$ assumed
piecewise continuous and differentiable,  in the sense that
$$
 \tLim{n}{\infty} \hpf(W_n \oleq z) \ = \ \Phi(z), \qquad  \text{for every } z \in \hf.$$
\end{theorem}

\begin{proof} The proof is established on the standard central limit
theorem, known for real distributions, cf. \cite{pCHU01a,pDUR96a}.
 We also use the fact that if $z_\pv \to z_{0, \pv} $  and $z_\ph \to
z_{0,\ph} $ as reals, then $z \to z_0$, and the properties of the
standard phantom normal distribution are as  addressed in
Proposition \ref{prop:Normal1}.

 By phantom computations, that are already familiar
to  the reader, we have
\begin{equation}\label{eq:z2}
W_n = \
\frac{S_{n,\pv} - n \mu_\pv}{\sig_\pv \sqrt{n}} + \hi \(
\frac{\Su{S_{n}} - n \Su{\mu}}{\su{\sig} \sqrt{n}}  -
\frac{S_{n,\pv} - n \mu_\pv}{\sig_\pv \sqrt{n}}  \).
\end{equation}
Suppose $z \in W_n$, then by the classical central limit theorem,
each component converges to the standard normal cumulative
distribution function, and thus using the realization property of
Equation \eqref{eq:normDomp} we get the desired.

When $z \notin W_n$, apply the same argument to $\supWnv{z}$, cf.
Equation \eqref{eq:supX}, for which $\Phi(z) = \Phi(\supWnv{z}) $
 by definition.
\end{proof}

The central phantom limit theorem is surprisingly general, maybe
even more general than the known classical one, which is a private
case of the phantom theorem. (Note that here  the integration is
performed along a path.) Besides independence, and the implicit
assumption that the mean and variance are well-defined and finite,
it places no other requirement on the distribution of the $X_i$,
even though they are phantoms, which could be discrete,
continuous, or mixed random variables.

This is of tremendous importance for several reasons, both
conceptual and practical. On the conceptual side, it indicates
that the sum of a large number of independent \prv 's is
approximately phantom normal. As such, it applies to many
situations in which a random effect is the sum of a large number
of small but independent random factors. Noise in many natural or
engineered systems has this property.

In a wide array of contexts, it has been found empirically that
the statistics of noise are well-described by (real) normal
distributions, and the central limit theorem provides a convincing
explanation of this phenomenon. Here, we add another argument,
recorded by the phantom \term \ which might provide more
information about the behavior of the noise.

On the practical side, the phantom central   limit theorem
eliminates the need for detailed probabilistic models and for
tedious manipulations of \pmf 's  and \pdf 's. Rather, it reduces
all the computations to a real familiar framework, and allows the
calculation of specific  probabilities by simply referring to the
table of the standard normal distribution. Furthermore, these
calculations only require  knowledge about the phantom means and
phantom variances.

\subsection{The strong law of large numbers}

\begin{theorem}[\bfem{The strong law of large numbers (SLLN)}] Let $X_1,X_2, \dots$ be a
sequence of independent identically distributed \prv's with mean
$\mu$. Then, the sequence of sample means $M_n = (X_1 + \cdots +
X_n)/n$ converges to $\mu$, with probability $1$, in the sense
that $$ \hpf\( \tLim{n}{\infty}  \frac{ X_1 + \cdots + X_n}{ n} =
\mu \)  = 1.$$
\end{theorem}

\begin{proof} Using the same argument as in the proof of  Theorem \ref{thm:WLLN}, $\hpf_\pv$ is a standard probability measure.  By the classical
SLLN, $\hpf_\pv \( \tLim{n}{\infty}  \frac{ X_1 + \cdots + X_n}{
n} = \mu \) =1$, and thus, since $\hpf$ is a phantom probability
measure, $\hpf_\ph \( \tLim{n}{\infty}  \frac{ X_1 + \cdots +
X_n}{ n} = \mu \) = 0  $. (Note that, for this purpose, the fact
that the random variable may take phantom values does not play a
role; equivalently, the $X_i$ can be viewed as random variables
that take values in $\Real^2$.)
\end{proof}

Consider a sequence of \prv's,  $X_1, X_2, \dots $,  (not
necessarily independent) associated with the same probability
model. Let $z$ be a phantom number. We say that $X_n$ converges to
$z_0$ \textbf{with probability} $1$ (or almost surely) if $$
\hpf\( \tLim{n}{\infty}
 X_n = z_0 \)  = 1.$$

In order to interpret the SSLN, one needs to use probabilistic
phantom models in terms of sample spaces. The contemplated
experiment is infinitely long and generates experimental values
for each one of the \prv's  in the sequence $X_1,X_2, \dots$.
Thus, one should rather  think of the sample space $\Om$ as a set
of infinite sequences $ \om = (x_1, x_2, \dots)$ of phantom
numbers: any such sequence is a possible outcome of the
experiment. Let us now define the subset $A$ of $\Om$ consisting
of those sequences $(x_1, x_2, \dots)$ whose long-term average is
$\mu$, i.e.,
$$(x_1, x_2, \dots) \in A \iff \tLim{n}{\infty} \frac{ x_1 +
\cdots + x_n}{ n} = \mu.$$ The SLLN states that most of the
phantom probabilities are concentrated on this particular subset
of $\Om$. Equivalently, the collection of outcomes that do not
belong to $A$ (infinite sequences whose long-term average  $ \neq
\mu$) has probability zero.

This means that the initial distortions of the probabilities
become meaningless as $n \to \infty$, as well as their phantom
\term s. (The latter have a special meaning when dealing with
Markov chains and stochastic processes, cf.  \cite{youMarkov}.)
Moreover, in the long term, the contribution of the phantom \term
\ lessens  and tends to zero.

The difference between the weak and the strong law is subtle and
deserves close scrutiny. The weak law states that the probability
$\hpf(\abs{M_n -\mu} \geq \ep)$ of a significant deviation of
$M_n$ from $\mu$ goes to zero as $n \to \infty $. Still, for any
finite $n$, this probability can be positive and it is conceivable
that once in a while, even if infrequently, $M_n$ deviates
significantly from $\mu$. The weak law provides no conclusive
information on the number of such deviations, but the strong law
does. According to the strong law, and with probability $1$, $M_n$
converges to $\mu$. This implies that for any given $\ep > 0$, the
difference $\abs{M_n - \mu}$ will exceed $\ep$ only a finite
number of times.

%


\end{document}